\documentclass{article}
\usepackage{color}
\usepackage{amssymb}
\usepackage{amstext}
\usepackage[dvips]{graphicx}

\newtheorem{theo}{Th\'{e}or\`{e}me}
\newtheorem{defi}{D\'{e}finition}
\newtheorem{propo}{Proposition}
\newtheorem{lemm}{Lemme}

\newtheorem{proposition}{Proposition}

\newtheorem{coro}{Corollaire}
\newcommand{\C}{\mathbb C}
\newcommand{\g}{\frak{g}}

\newcommand{\n}{\frak{n}}
\newcommand{\m}{\frak{t}}
\newcommand{\s}{\frak{s}}
\newcommand{\h}{\frak{h}}
\begin{document}

\title{Alg\`ebres de Lie : classifications, d\'{e}formations et rigidit\'{e},
g\'{e}om\'{e}trie diff\'{e}rentielle.}
\author{Michel GOZE}
\maketitle
\title{Cours donn\'{e} \`{a} l'ENSET d'Oran en novembre 2006 durant la cinqui%
\`{e}me Ecole de G\'{e}om\'{e}trie diff\'{e}rentielle et Syst\`{e}mes
Dynamiques}
\tableofcontents

{\noindent }

\bigskip

\part{Alg\`{e}bres de Lie : g\'{e}n\'{e}ralit\'{e}s et classifications}

\section{Alg\`{e}bres de Lie : d\'{e}finitions, exemples.}

\medskip

\subsection{ D\'{e}finition et exemples}

\medskip

Dans tout ce travail, les alg\`{e}bres de Lie consid\'{e}r\'{e}es seront
complexes.\ Lorsque nous serons int\'{e}ress\'{e}s par le cas r\'{e}el ou par des alg\`ebres de 
Lie sur un anneau quelconque,
nous le pr\'{e}ciserons alors.

\begin{defi}
Une alg\`ebre de Lie  est un couple $(\g,\mu)$ o\`u $\g$ est un espace vectoriel complexe et $\mu$ 
une application bilin\'eaire
$$
\mu :\g \times \g \rightarrow \g
$$
satisfaisant :
$$
\mu (X,Y)=-\mu (Y,X),\quad \forall X,Y \in \g,$$
$$
 \mu (X,\mu (Y,Z))+\mu (Y,\mu (Z,X))+\mu (Z,\mu (X,Y))=0,\quad \forall
X,Y,Z\in \g.$$
\end{defi}

\bigskip

\noindent Cette derni\`{e}re identit\'{e} est appel\'{e}e l'identit\'{e} de
Jacobi. Si $\mathfrak{g}$ est un espace vectoriel de $\dim $ension finie $n$%
, on dira que $(\mathfrak{g},\mu )$ est une alg\`{e}bre de Lie de dimension $%
n$. Sinon on dira que l'alg\`{e}bre de Lie $\mathfrak{g}$ est de dimension
infinie.\ Par soucis de simplification de langage, lorsque cela n'entraine
aucune cons\'{e}quence, on parlera de $\mathfrak{g}$ alg\`{e}bre de Lie au
lieu du couple $(\mathfrak{g,\mu ).}$

\medskip

\subsection{Exemples}

\medskip

1. Soit $V$ un espace vectoriel de dimension finie ou non.\ Si $\mu =0,$
alors $\mathfrak{g}=(V,0)$ est une alg\`{e}bre de Lie appel\'{e}e dans ce
cas ab\'{e}lienne.

\medskip

\noindent 2. Soit $\mathfrak{g}$ une espace vectoriel de dimension $2$.
Alors, pour toute application bilin\'{e}aire antisym\'{e}trique $\mu $ sur $%
\mathfrak{g}$ \`{a} valeurs dans $\mathfrak{g}$, le couple $\mathfrak{g}%
=(V,\mu )$ est une alg\`{e}bre de Lie. En effet la condition de Jacobi est
toujours, dans ce cas, satisfaite.

\noindent 3. Soit $sl(2,\mathbb{C})$ l'espace vectoriel des matrices d'ordre 
$2$ de trace nulle. Le produit 
\[
\mu (A,B)=AB-BA 
\]%
est bien d\'{e}fini sur $sl(2,\mathbb{C})$ car $tr(AB-BA)=0$ d\`{e}s que $%
A,B\in sl(2,\mathbb{C})$. Comme cette multiplication v\'{e}rifie l'identit%
\'{e} de Jacobi, $sl(2,\mathbb{C})$ est une alg\`{e}bre de Lie complexe de
dimension $3$.

\subsection{Sous-alg\`{e}bres de Lie, morphismes.}

Un sous-espace vectoriel $\mathfrak{h}$ de $\mathfrak{g}$ est une sous-alg%
\`{e}bre de Lie de $(\mathfrak{g,\mu )}$ \ si pour tout $X,Y\in \mathfrak{h}$
on a $\mu (X,Y)\in \mathfrak{h}$. Bien \'{e}videmment, une sous-alg\`{e}bre
de Lie est une alg\`{e}bre de Lie. Une type sp\'{e}cial de sous-alg\`{e}bres
est l'id\'{e}al.\ Un id\'{e}al $\mathfrak{I}$ de $(\mathfrak{g,\mu )}$ est
une sous-alg\`{e}bre de Lie v\'{e}rifiant%
\[
\forall X\in \mathfrak{I},\forall Y \in g,\mu (X,Y)\in \mathfrak{I}. 
\]

Soient ($\mathfrak{g}_{1},\mu _{1}$ ) et ($\mathfrak{g}_{2},\mu _{2}$) deux
alg\`{e}bres de Lie. Une application lin\'{e}aire%
\[
\varphi :\mathfrak{g}_{1}\longrightarrow \mathfrak{g}_{2} 
\]%
est un homomorphisme d'alg\`{e}bres de Lie si%
\[
\varphi (\mu _{1}(X,Y))=\mu _{2}(\varphi (X),\varphi (Y)) 
\]%
pour tout $X,Y\in \mathfrak{g}_{1}$. Si de plus $\varphi $ est un
isomorphisme lin\'{e}aire, on dira que c'est un isomorphisme d'alg\`{e}bres
de Lie. Enfin pr\'{e}cisons la notion bien utile d'alg\`{e}bres de Lie
quotient.\ Si $\mathfrak{I}$ est un id\'{e}al de $\ $l'alg\`{e}bre de Lie $%
\mathfrak{g,}$ il existe une unique structure d'alg\`{e}bre de Lie sur
l'espace vectoriel quotient $\mathfrak{g/I}$ pour laquelle la projection
canonique%
\[
\pi :\mathfrak{g\longrightarrow g/I} 
\]%
soit un homorphisme surjectif d'alg\`{e}bres de Lie.

\subsection{Alg\`{e}bres Lie-admissibles}

Soit $\mathcal{A}$ une alg\`{e}bre associative complexe dont la
multiplication est not\'{e}e $A.B$ avec $A,B\in \mathcal{A}$. On voit
facilement que le crochet suivant 
\[
\lbrack A,B]=AB-BA 
\]%
d\'{e}finit une structure d'alg\`{e}bre de Lie sur l'espace vectoriel
sous-jacent \`{a} $\mathcal{A}$. On note par $\mathcal{A}^{L}$ cette alg\`{e}%
bre de Lie. Par exemple, si $M(n,\mathbb{C})$ est l'espace vectoriel des
matrices d'ordre $n$, c'est une alg\`{e}bre associative pour le produit
usuel des matrices.\ Ainsi $[A,B]=AB-BA$ d\'{e}finit une structure d'alg\`{e}%
bre de Lie sur $M(n,\mathbb{C})$. Elle est not\'{e}e dans ce cas $gl(n,%
\mathbb{C})$. Pr\'{e}cisons toutefois qu'il existe des alg\`{e}bres de Lie
qui ne sont pas d\'{e}finies \`{a} partir d'alg\`{e}bres associatives.

\begin{defi}
Une alg\`ebre Lie-admissible est une alg\`ebre (non-associative) $\mathcal{A}$ dont le produit $A.B$ est tel que
$$[A,B]=AB-BA$$
est un produit d'alg\`ebre de Lie.  
\end{defi}

\bigskip

\noindent Ceci est \'{e}quivalent \`{a} dire que $A.B$ v\'{e}rifie: 
\[
\begin{array}{l}
(A.B).C-A.(B.C)-(B.A).C+B.(A.C)-(A.C).B+A.(C.B)-(C.B).A+C.(B.A) \\ 
+(B.C).A-B.(C.A)+(C.A).B-C.(A.B)=0%
\end{array}%
\]%
pour tout $A,B,C\in \mathcal{A}.$

\bigskip

\noindent\textbf{Exemple : Les alg\`{e}bres sym\'etriques gauche.}

\bigskip

Les alg\`{e}bres sym\'etriques gauche sont des exemples int\'{e}ressants, pour
diverses raisons, d'alg\`{e}bres Lie-admissibles.\ Une alg\`{e}bre sym\'etrique gauche
 est une alg\`{e}bre non-associative (on entend par l\`{a} une alg\`{e}%
bre non n\'{e}cessairement associative) dont le produit $A.B$ v\'{e}rifie 
\[
(A.B).C-A.(B.C)-(B.A).C+B.(A.C)=0. 
\]%
Il est clair qu'une alg\`{e}bre sym\'etrique gauche est Lie-admissible. Elles
sont \'{e}tudi\'{e}es par exemple, dans la recherche des connexions affines
invariantes \`{a} gauche sur un groupe de Lie sans courbure ni torsion.\ En
effet si un groupe de Lie $G$ admet une telle connexion, son alg\`{e}bre de
Lie $\mathfrak{g}$ provient d'une alg\`{e}bre sym\'etrique gauche, c'est-\`{a}-dire,
il existe sur l'espace vectoriel $\mathfrak{g}$ une structure d'alg\`{e}bre
sym\'etrique gauche de produit $A.B$ tel que le produit d'alg\`{e}bre de Lie de $%
\mathfrak{g}$ soit donn\'{e} par 
\[
\mu (A,B)=A.B-B.A 
\]%
Dans ce cas aussi, il existe des alg\`{e}bres de Lie qui ne sont donn\'{e}es
par aucune alg\`{e}bre sym\'etrique gauche, ce qui signifie que pour les groupes de
Lie correspondants, toute connexion affine invariante \`{a} gauche sans
torsion a une courbure non triviale. Nous donnerons un exemple dans la derni%
\`{e}re partie de ce cours.\ A titre d'illustration regardons le cas
particulier o\`{u} $\mathfrak{g}$ est une alg\`{e}bre de Lie ab\'{e}lienne
de dimension $2$. Dans ce cas toute alg\`{e}bre sym\'etrique gauche dont
l'anticommutateur donne le crochet (trivial) de $\mathfrak{g}$ est
commutative.\ C'est donc une alg\`{e}bre associative commutative. Comme on
connait la classification de ces alg\`{e}bres associatives, on en d\'{e}duit
la classification des structures affines dans le plan d\'{e}finies par une
connexion affine sans courbure ni torsion dans le groupe de Lie $\mathbb{R}%
^{2}$ (dans le cas r\'{e}el) ou $\mathbb{C}^{2}$ (dans le cas complexe). Par
exemple, consid\'{e}rons l'alg\`{e}bre associative commutative de dimension $%
2$ donn\'{e}e dans une base $\left\{ e_{1},e_{2}\right\} $ par 
\[
e_{1}.e_{1}=e_{1},\ e_{1}.e_{2}=e_{2}.e_{1}=e_{2},\ e_{2}.e_{2}=e_{2}. 
\]%
Soit $X=ae_{1}+be_{2}.$ La translation \`{a} gauche $l_{X}$ associ\'{e}e 
\`{a} ce produit est donn\'{e}e par 
\[
l_{X}(e_{1})=ae_{1}+be_{2},\ l_{X}(e_{2})=(a+b)e_{2}. 
\]%
L'alg\`{e}bre de Lie ab\'{e}lienne$\mathfrak{\ g}$ d\'{e}finie par $\left[
e_{1},e_{2}\right] =e_{1}.e_{2}-e_{2}.e_{1}=0$ \ se repr\'{e}sente comme une
sous alg\`{e}bre de l'alg\`{e}bre de Lie $Aff(\mathbb{R}^{2})$ du 
groupe de Lie affine du plan (r\'{e}el ou complexe) de la mani\`{e}re suivante
: 
\[
\mathfrak{g=}\left\{ \left( 
\begin{array}{ll}
l_{X} & X \\ 
0 & 0%
\end{array}%
\right) \ =\left( 
\begin{array}{lll}
a & 0 & a \\ 
b & a+b & b \\ 
0 & 0 & 0%
\end{array}%
\right) \right\} . 
\]%
Le groupe de Lie correspondant s'\'{e}crit en consid\'{e}rant
l'exponentielle de ces matrices (voir paragraphe suivant) 
\[
G=\left\{ \left( 
\begin{array}{lll}
e^{a} & 0 & e^{a}-1 \\ 
e^{a}(e^{b}-1) & e^{a}e^{b} & e^{a}(e^{b}-1) \\ 
0 & 0 & 1%
\end{array}%
\right) .\right\} \ 
\]%
Ceci revient \`{a} dire que $G$ est le groupe des transformations affines du
plan : 
\[
\left\{ 
\begin{array}{l}
x\longrightarrow e^{a}x+e^{a}-1 \\ 
y\longrightarrow e^{a}(e^{b}-1)x+e^{a}e^{b}y+e^{a}(e^{b}-1)%
\end{array}%
\right. 
\]%
Notons que tous ces groupes sont class\'{e}s, et que le r\'{e}sultat est 
\'{e}galement connu pour $n=3$.

\medskip

\noindent Notons \'egalement la notion d'alg\`ebres sym\'etriques droite, encore appel\'ees alg\`ebres pr\'e-Lie 
qui v\'erifient
$$(A.B).C-A.(B.C)=(A.C).B-A.(C.B).$$
Elles jouent un r\^ole important dans l'\'etude des alg\`ebres de Gerstenhaber, de la cohomologie 
d'Hochschild ou m\^eme dans l'\'etude des alg\`ebres de Rota-Baxter.

\subsection{Alg\`{e}bres de Lie de dimension infinie}

La th\'{e}orie des alg\`{e}bres de Lie de dimension infinie, c'est-\`{a}-dire dont l'espace vectoriel sous-jacent est de dimension infinie est
assez diff\'{e}rente de celle de la dimension finie.\ Nous ne l'aborderons
pas trop dans ce cours.\ Dans ce cas, la structure topologique de l'espace
vectoriel joue un r\^{o}le pr\'{e}pond\'{e}rant.\ Par exemple, afin de
donner des contre exemples au troisi\`{e}me th\'{e}or\`{e}me de Lie-Cartan,
W.T\ Van Est a \'{e}tudi\'{e} les alg\`{e}bres de Lie de Banach. Certaines
alg\`{e}bres de Lie infinies sont de nos jours fortement \'{e}tudi\'{e}es.\
Par exemple, les alg\`{e}bres de Kac-Moody sont des alg\`{e}bres de Lie
infinies gradu\'{e}es et d\'{e}finies par g\'{e}n\'{e}rateurs et relations.
Elles sont construites d'une mani\`{e}re analogue \`{a} celle des alg\`{e}%
bres de Lie simples. Un autre exemple est donn\'{e}e par l'alg\`{e}bre de Lie
des champs de vecteurs sur une vari\'{e}t\'{e} diff\'{e}rentiable $M$. Le
produit de Lie est alors le crochet de Lie des champs de vecteurs.\ La
structure d'une telle alg\`{e}bre est tr\`{e}s compliqu\'{e}e.\ Elle a \'{e}t%
\'{e} \'{e}tudi\'{e}e lorsque $M=\mathbb{R}$ ou $M=S^{1}$. Dans ce cas l'alg%
\`{e}bre de Lie est associ\'{e}e (voir paragraphe suivant) au groupe de Lie
des diff\'{e}omorphismes de $\mathbb{R}$ ou $S^{1}$. Une autre classe int\'eressante d'alg\`ebres
de Lie infinies concerne les alg\`ebres de Lie-Cartan. Elles sont d%
\'{e}finies comme les alg\`{e}bres de Lie des transformations infinit\'{e}%
simales (des champs de vecteurs) qui laissent invariants une structure donn%
\'{e}e, comme une structure symplectique ou une structure de contact. Par
exemple si $(M,\Omega )$ est une vari\'{e}t\'{e} symplectique, c'est-\`{a}%
-dire $\Omega $ est une forme symplectique sur $M$, on consid\`{e}re alors 
\[
L(M,\Omega )=\{X\ {\mbox{\rm champs \ de  \ vecteurs \ sur}}\ M,L_{X}\Omega =0\} 
\]%
o\`{u} 
\[
L_{X}\Omega =i(X)d\Omega +d(i(X)\Omega )=d(i(X)\Omega ) 
\]%
est la d\'{e}riv\'{e}e de Lie. Alors $L(M,\Omega )$ est une alg\`{e}bre de
Lie r\'{e}elle de dimension infinie. Elle admet une sous-alg\`{e}bre $L_{0}$
constitut\'{e}e des champs de vecteurs de $L(M,\Omega )$ \`{a} support
compact . Andr\'{e} Lichnerowicz prouva que toute sous-alg\`{e}bre de Lie de
dimension finie de $L_{0}$ est r\'{e}ductive (c'est-\`{a}-dire produit
direct d'une sous-alg\`{e}bre semi simple par un centre ab\'{e}lien) et que
tout id\'{e}al non nul est de dimension infinie. \medskip

\section{Alg\`{e}bres de Lie et groupes de Lie}

\subsection{L'alg\`{e}bre de Lie d'un groupe de Lie}

Soit $G$ un groupe de Lie complexe de dimension $n$. Pour tout $g\in G$, \
notons par $L_{g}$ l'automorphisme de $G$ donn\'{e} par 
\[
L_{g}(x)=gx. 
\]%
Cet automorphisme est appel\'{e} la translation \`{a} gauche par $\ g$. Sa
diff\'{e}rentielle $(L_{g})_{x}^{\ast }$ en un point $x\in G,$ est
l'isomorphisme vectoriel 
\[
(L_{g})_{x}^{\ast }:T_{x}(G)\longrightarrow T_{gx}(G) 
\]%
o\`{u} $T_{x}(G)$ d\'{e}signe l'espace tangent au point $x$ \`{a} $G$. Un
champ de vecteurs $X$ sur $G$ est dit invariant \`{a} gauche s'il v\'{e}%
rifie 
\[
(L_{g})_{x}^{\ast }(X(x))=X(gx) 
\]%
pour tout $x$ et $g$ dans $G$. On montre que si $[X,Y]$ d\'{e}signe le
crochet de Lie des champs de vecteurs, alors si $X$ et $Y$ sont des champs
invariants \`{a} gauche sur $G$, le crochet de Lie $[X,Y]$ est aussi un
champ invariant \`{a} gauche. Ceci montre que l'espace vectoriel des champs
invariants \`{a} gauche sur $G$, muni du crochet de Lie, est une alg\`{e}bre
de Lie. On la note $L(G)$ et elle est appel\'{e}e l'alg\`ebre
 de Lie du groupe de Lie $G.$ Comme un champ invariant \`{a}
gauche $X$ est enti\`{e}rement d\'{e}fini par sa valeur $X(e)$ en l'\'{e}l%
\'{e}ment neutre $e$ de $G$, l'espace vectoriel $L(G)$ s'identifie
naturellement \`{a} l'espace tangent T$_{e}(G)$ de $G$ en $e$. Si $u,v\in
T_{e}(G),$ il existe $X,Y\in L(G)$ tels que $u=X(e)$ et $v=Y(e)$. Posons%
\[
\mu (u,v)=[X,Y](e). 
\]%
Alors $(T_{e}(G),\mu )$ est une alg\`{e}bre de Lie isomorphe \`{a} $L(G)$.\
On confondra souvent les deux. Par exemple l'alg\`{e}bre de Lie du groupe de
Lie alg\'{e}brique $SL(2,\mathbb{C})$ est isomorphe \`{a} $sl(2,\mathbb{C})$.

\medskip

Ainsi cette construction permet de d\'{e}finir pour chaque groupe de Lie,
une et une seule (classe d'isomorphie d') alg\`{e}bre de Lie.\ Mais
l'inverse n'est pas vrai.\ En dimension finie, plusieurs groupes de Lie
peuvent avoir la m\^{e}me alg\`{e}bre de Lie.\ En dimension infinie, il
existe des alg\`{e}bres de Lie qui ne sont des alg\`{e}bres d'aucun groupe de
Lie. Pr\'{e}cisons bri\`{e}vement ces deux remarques. Supposons tout d'abord
que $\mathfrak{\ g}$ soit une alg\`{e}bre de Lie de dimension finie. A
partir de sa multiplication $\mu $, on paut d\'{e}finir un groupe local dont
le produit est donn\'{e} par la formule de Campbell-Hausdorff : 
\[
X.Y=X+Y+\displaystyle\frac{1}{2}\mu (X,Y)+\displaystyle\frac{1}{12}\mu (\mu (X,Y),Y)-\displaystyle\frac{1}{12}\mu (\mu (X,Y),X)+.... 
\]%
qui est une somme infinie de termes exprim\'{e}s \`{a} l'aide de la
multiplication $\mu $. Cette structure locale peut \^{e}tre \'{e}tendue en
une structure globale de groupe de Lie en imposant une condition topologique
de simple connexit\'{e} et de connexit\'{e}. On a donc une correspondance
biunivoque entre l'ensemble des alg\`{e}bres de Lie de dimension finie
complexes (mais ceci reste vrai dans le \ cas r\'{e}el) et l'ensemble des
groupes de Lie connexes et simplement connexes dont la dimension en tant que
vari\'{e}t\'{e} diff\'{e}rentielle est la dimension de l'alg\`{e}bre de Lie.
Deux groupes de Lie de dimension finie ayant des alg\`{e}bres de Lie
isomorphes sont donc localement isomorphes.

\bigskip

Dans le cas de la dimension infinie, la situation est diff\'{e}rente.\
W.T.Van Est a montr\'{e} l'existence d'alg\`{e}bres de Lie de Banach de
dimension infinie qui n'\'{e}taient des alg\`{e}bres de Lie d'aucun groupes
de Lie de dimension infinie.\ Ainsi, dans ce cas, il n'y a pas d'\'{e}%
quivalent au troisi\`{e}me th\'{e}or\`{e}me de Lie-Cartan.

\subsection{Relation entre un groupe de Lie et son alg\`{e}bre de Lie}

\begin{defi}
Un groupe de Lie est dit lin\'eaire si c'est un sous-groupe de Lie du groupe 
$GL(n,\C)$ des matrices inversibles d'ordre $n$.
\end{defi}

\bigskip

Si $G$ est un groupe de Lie lin\'{e}aire, son alg\`{e}bre de Lie $\mathfrak{g%
}$ est une sous-alg\`{e}bre de Lie de $gl(n,\mathbb{C}),$ l'alg\`ebre de Lie des matrices complexes d'ordre $n$.
Dans ce cas l'application Exponentielle d\'{e}finit une correspondance entre
l'alg\`{e}bre de Lie $\mathfrak{g}$ et le groupe de Lie $G$.\ Rappelons que
l'application exponentielle 
\[
Exp:gl(n,\mathbb{C})\longrightarrow Gl(n,\mathbb{C}) 
\]%
est donn\'{e}e par la s\'{e}rie enti\`{e}re convergente 
\[
Exp(A)=Id+A+\frac{A^{2}}{2!}+\frac{A^{3}}{3!}+...=\sum_{n=0}^{\infty
}\frac{A^{n}}{n!}. 
\]%
Si $G$ est un sous-groupe de Lie de $GL(n,\mathbb{C}),$ 
l'application exponentielle envoie l'alg\`{e}bre de Lie $\mathfrak{g}$ de $%
G$ dans $G$.\ Ainsi nous d\'{e}finissons une application exponentielle pour
tous les groupes lin\'{e}aires et leurs alg\`{e}bres de Lie.

La d\'{e}finition ci-dessus est d'un usage tr\`{e}s pratique mais n'est pas
donn\'{e}e pour les groupes de Lie qui ne sont pas des groupes de matrices.
La difficult\'{e} pourrait \^{e}tre contourn\'{e}e en s'appuyant sur le th%
\'{e}or\`{e}me d'Ado. Ce th\'{e}or\`{e}me pr\'{e}cise que toute alg\`{e}bre
de Lie de dimension finie peut se repr\'{e}senter comme une alg\`{e}bre de
matrices, plus pr\'{e}cis\'{e}ment il existe un entier $N$ tel que l'alg\`{e}%
bre de Lie donn\'{e}e $\mathfrak{g}$ soit isomorphe \`{a} une sous-alg\`{e}%
bre de dimension $n$ de $gl(N,\mathbb{C)}.\ $\ Mais l'application du th\'{e}%
or\`{e}me d'Ado est parfois difficile car nous ne connaissons pas pr\'{e}cis%
\'{e}ment l'entier $N$ et l'application exponentielle est dans ce contexte
non seulement difficile \`{a} \'{e}crire mais d\'{e}pend aussi du choix de
la repr\'{e}sentation de $\mathfrak{g}$ dans $gl(N,\mathbb{C)}.$ Nous allons
donc d\'{e}finir directement, pour un groupe de Lie abstrait $G$, cette
application exponentielle et v\'{e}rifier qu'elle co\"{\i}ncide avec
l'application exponentielle des groupes de Lie lin\'{e}aires. Chaque vecteur 
$X$ de $\mathfrak{g}$ d\'{e}termine une application lin\'{e}aire de $\mathbb{%
R}$ dans $\mathfrak{g}$ ayant $X$ comme image de $1$ et qui soit
un homomorphisme d'alg\`{e}bre de Lie. Comme $\mathbb{R}$ est l'alg\`{e}bre
de Lie r\'{e}elle du groupe de Lie connexe et simplement connexe $\mathbb{R}$%
, cette application induit un homomorphisme de groupes de Lie 
\[
c:\mathbb{R}\longrightarrow G 
\]%
tel que 
\[
c(s+t)=c(s)+c(t) 
\]%
pour tout $s$ et $t$. L'op\'{e}ration dans la partie droite de la formule
correspond \`{a} la multiplication dans $G$.\ Compte tenu de la ressemblance
de cette formule avec la propri\'{e}t\'{e} caract\'{e}ristique de l'application exponentielle des
matrices, on est conduit \`{a} poser la d\'{e}finition suivante: 
\[
Exp(X)=c(1). 
\]%
Cette application est appel\'{e}e l'application exponentielle et envoie bien
l'alg\`{e}bre de Lie $\mathfrak{g}$ dans le groupe de Lie $G$. Elle d\'{e}%
termine un diff\'{e}omorphisme entre un voisinage de $\ 0$ dans $\mathfrak{g}
$
et un voisinage de l'\'el\'ement neutre dans $G$. L'application exponentielle n'est pas
toujours surjective m\^{e}me si le groupe $G$ est suppos\'{e} connexe. Par
exemple, on montre que l'application exponentielle 
\[
Exp:sl(2,\mathbb{C})\longrightarrow Sl(2,\mathbb{C}) 
\]%
n'est pas surjective. Mais si l'alg\`{e}bre de Lie $\mathfrak{g}$ est
nilpotente (voir la d\'{e}finition dans les paragraphes suivants), alors $%
Exp $ est bijective.

\medskip

\section{Classifications des alg\`{e}bres de Lie complexes}

\medskip

\subsection{Alg\`{e}bres de Lie isomorphes}

\noindent 

\begin{defi}
Deux alg\`ebres de Lie $\g$ and $\g'$ de dimension $n$ de multiplication  $\mu $ et $\mu ^{\prime }$ sont dites isomorphes s'il
existe $f\in Gl(n,\mathbb{C)}$ tel que
$$
\mu _f(X,Y)=f*\mu (X,Y)=f^{-1}(\mu (f(X),f(Y)))
$$
pour tout $X,Y\in \g.$ \end{defi}

\bigskip

Par exemple, toute alg\`{e}bre de Lie de dimension $2$ a une multiplication
qui v\'{e}rifie 
\[
\mu (e_{1},e_{2})=ae_{1}+be_{2}. 
\]%
Un tel produit v\'{e}rifie toujours l'identit\'{e} de Jacobi. Si $a$ ou $b$
est non nul, par exemple $b$, le changement de base 
\[
\left\{ 
\begin{array}{l}
f(e_{1})=1/b\ e_{1} \\ 
f(e_{2})=a/be_{1}+e_{2}%
\end{array}%
\right. 
\]%
d\'{e}finit une multiplication isomorphe donn\'{e}e par 
\[
\mu _{f}(e_{1},e_{2})=e_{2}. 
\]%
On en d\'{e}duit que toute alg\`{e}bre de Lie de dimension $2$ est soit ab%
\'{e}lienne, soit isomorphe \`{a} l'alg\`{e}bre de Lie dont le produit v\'{e}%
rifie 
\[
\mu (e_{1},e_{2})=e_{2}. 
\]%
La classification des alg\`{e}bres de Lie de dimension $n$ consiste \`{a} d%
\'{e}crire un repr\'{e}sentant de chacune des classes d'isomorphie.\ Le r%
\'{e}sultat pr\'{e}c\'{e}dent donne la classification des alg\`{e}bres de
Lie complexes (et r\'{e}elles) de dimension $2$. Notons que la
classification g\'{e}n\'{e}rale est de nos jours encore un probl\`{e}me
ouvert.\ Elle est parfaitement connue jusqu'en dimension $5$.\ Au del\`{a},
seules des classifications partielles sont \'{e}tablies, ou bien des familles
particuli\`{e}res d'alg\`{e}bres de Lie ont \'{e}t\'{e} class\'{e}es.\ Nous
allons pr\'{e}senter ces familles.

\subsection{Alg\`{e}bres de Lie simples et semi-simples}

\noindent 

\begin{defi}
Une alg\`ebre de Lie $\g$ est appel\'ee simple si sa dimension est sup\'erieure ou \'egale \`a 
2
 et si elle ne contient pas d'id\'eaux propres (autre que $\{0\}$ et $\g$).
\end{defi}

\bigskip

La classification des alg\`{e}bres simples complexes (et r\'{e}elles) est
bien connue.\ Elle est due essentiellement aux travaux d'Elie Cartan, de
Dynkin et de Killing. Elle se r\'{e}sume au r\'{e}sultat suivant:

\begin{propo}
\noindent Toute alg\`{e}bre simple complexe de dimension finie est

i) soit isomorphe \`{a} une alg\`{e}bre de type classique, c'est-\`{a}-dire 
\`{a} l'une des alg\`{e}bres suivantes :$su(n,\mathbb{C})$ (type $A_{n}$), $%
so(2n+1,\mathbb{C})$ (type $B_{n}$), $sp(n,\mathbb{C})$ (type $C_{n}$), $%
so(2n,\mathbb{C})$ (type $D_{n}$)

ii) soit isomorphe \`{a} une alg\`{e}bre exeptionnelle $E_{6},\ E_{7},\
E_{8},\ F_{4},\ G_{2}$.
\end{propo}

\noindent On peut trouver les d\'{e}finitions pr\'{e}cises de ces alg\`{e}%
bres dans le livre de J.P. Serre intitul\'{e} SemiSimple Lie algebras. Notons
que l'alg\`{e}bre de Lie $E_{8}$ a eu les honneurs  de tous les
medias ces derniers mois.

\begin{defi}
\medskip Une alg\`{e}bre de Lie est appel\'{e} semi-simple si elle est non
nulle et si elle n'admet pas d'id\'{e}aux ab\'{e}liens non nuls.
\end{defi}

Ces alg\`{e}bres se caract\'{e}risent aussi par le fait qu'elles sont des
produits directs d'alg\`{e}bres simples.\ Rappelons que si $\mathfrak{g}_{1} 
$ et $\mathfrak{g}_{2}$ sont deux alg\`{e}bres de Lie de multiplications
respectives $\mu _{1}$ et $\mu _{2},$ alors le produit direct (ou la somme
directe externe $\mathfrak{g}_{1}\oplus \mathfrak{g}_{2}$ )\ est aussi une
alg\`{e}bre de Lie pour le produit 
\[
\mu (X_{1}+X_{2},Y_{1}+Y_{2})=\mu _{1}(X_{1},Y_{1})+\mu _{2}(X_{2},Y_{2}) 
\]%
pour tout $X_{1},Y_{1}\in \mathfrak{g}_{1}$ et $X_{2},Y_{2}\in \mathfrak{g}%
_{2}$. La classification des alg\`{e}bres simples implique celle des
semi-simples. Notons \'{e}galement que les alg\`{e}bres de Lie semi-simples
se caract\'{e}risent par le fait que la forme bilin\'{e}aire de
Killing-Cartan%
\[
K(X,Y)=Tr(adX\circ adY) 
\]%
est non d\'{e}g\'{e}n\'{e}r\'{e}e. Cette forme d\'{e}finit donc un produit
scalaire invariant au sens suivant :%
\[
K(ad(Y)(X),Z)+K(X,ad(Y)(Z))=0 
\]%
pour tout $X,Y,Z\in \mathfrak{g}$ o\`{u} $adY$ est l'endomorphisme donn\'{e}%
e par $ad(Y)(X)=\mu (Y,X),$ appel\'{e} application adjointe.

\subsection{Alg\`{e}bres de Lie nilpotentes}

Soit $\mathfrak{g}$ une alg\`{e}bre de Lie de multiplication $\mu $ et consid%
\'{e}rons la suite d\'{e}croissante suivante d'id\'{e}aux de $\mathfrak{g}$ 
\[
\left\{ 
\begin{array}{l}
{\mathcal{C}}^{0}(\mathfrak{g})=\mathfrak{g} \\ 
{\mathcal{C}}^{1}(\mathfrak{g})=\mu (\mathfrak{g},\mathfrak{g}) \\ 
{\mathcal{C}}^{p}(\mathfrak{g})=\mu ({\mathcal{C}}^{p-1}(\mathfrak{g}),%
\mathfrak{g})\ for\ p>2.%
\end{array}%
\right. 
\]%
o\`{u} $\mu ({\mathcal{C}}^{p-1}(\mathfrak{g}),\mathfrak{g})$ est la sous-alg%
\`{e}bre de Lie de $\mathfrak{g}$ engendr\'{e}e par les produits $\mu (X,Y)$
avec $X\in {\mathcal{C}}^{p-1}(\mathfrak{g})$ et $Y\in \mathfrak{g}$.

\medskip

\noindent 

\begin{defi}
Une alg\`ebre de Lie  $\g$ est appel\'ee nilpotente s'il existe un entier $k$ tel que
$${\mathcal{C}}^k(\g)=\{0\}.$$
Si un tel entier existe, le plus petit $k$ tel que ${\mathcal{C}}^k(\g)=\{0\}$ est appel\'e l'indice de nilpotence ou nilindex de
 $\g$.
\end{defi}

\medskip Pour toute alg\`{e}bre de Lie nilpotente de dimension $n$, son
nilindex est inf\'{e}rieur ou \'{e}gal \`{a} $n-1.$

\noindent \textbf{Exemples.}

1. Toute alg\`{e}bre de Lie ab\'{e}lienne v\'{e}rifie ${\mathcal{C}}^{1}(%
\mathfrak{g})=\{0\}$. Elle est donc nilpotente d'indice $1$. Bien entendu,
toute alg\`{e}bre nilpotente d'indice $1$ est ab\'{e}lienne.

2. L'alg\`{e}bre de Heisenberg de dimension $(2p+1)$ est l'alg\`{e}bre de
Lie $\mathfrak{h}_{2p+1}$ dont le produit dans une base $%
\{e_{1},...,e_{2p+1}\}$ est donn\'{e} par 
\[
\mu (e_{2i+1},e_{2i+2})=e_{2p+1} 
\]%
pour $i=0,...,p-1$, les autres produits \'{e}tant nuls. Cette alg\`{e}bre de
Lie est nilpotente d'indice $2$.\ La sous-alg\`{e}bre d\'{e}riv\'{e}e ${%
\mathcal{C}}^{1}(\mathfrak{g})$ est engendr\'{e}e par $\{e_{2p+1}\}$ et co%
\"{\i}ncide avec le centre $Z(\mathfrak{g)}$.

\bigskip

Le r\'{e}sultat suivant est fondamental dans l'\'{e}tude des alg\`{e}bres
nilpotentes

\medskip

\noindent 

\begin{theo}Th\'eor\`eme de Engel. Soit $(\g,\mu)$ une alg\`ebre de Lie complexe de dimension finie. Alors $\g$
est nilpotente si et seulement si  pour tout $X \in \g$, l'endomorphisme
$$ad X : \g \longrightarrow \g$$
donn\'e par $adX(Y)=\mu(X,Y)$ est nilpotent.
\end{theo}

\medskip On peut s'int\'{e}resser \`{a} des sous familles d'alg\`{e}bres
nilpotentes.\ En particulier

\noindent 

\begin{defi}
Une alg\`ebre de Lie nilpotente de dimension 
$n$
 est appel\'ee filiforme si son indice de nilpotence est \'egal \`a  $n-1$.
\end{defi}

\medskip

\noindent Par exemple, l'alg\`{e}bre de Lie de dimension $4$, d\'{e}finie
dans la base $\left\{ e_{1},e_{2},e_{3},e_{4}\right\} $ par 
\[
\left\{ 
\begin{array}{l}
\mu (e_{1},e_{2})=e_{3} \\ 
\mu (e_{1},e_{3})=e_{4}%
\end{array}%
\right. 
\]%
est filiforme.

\medskip

\bigskip La classification des alg\`{e}bres de Lie nilpotentes complexes (ou
r\'{e}elles) n'est connue, de nos jours, que jusqu'en dimension $7$.

\begin{itemize}
\item En dimension $1$ et $2$, toute alg\`{e}bre nilpotente est ab\'{e}%
lienne.

\item En dimension $3$, toute alg\`{e}bre nilpotente non ab\'{e}lienne est
isomorphe \`{a} l'alg\`{e}bre de Heisenberg $\mathfrak{h}_{3}.\ $Rappelons
que le produit est d\'{e}fini par 
\[
\mu (e_{1},e_{2})=e_{3} 
\]%
o\`{u} $\left\{ e_{1},e_{2},e_{3}\right\} $ est une base de $\mathfrak{g.}$

\item En dimension $4$, l'alg\`{e}bre est soit ab\'{e}lienne, soit isomorphe 
\`{a} une somme directe externe d'une alg\`{e}bre de Lie ab\'{e}lienne de
dimension $1$ avec l'alg\`{e}bre de Heisenberg de dimension $3$, soit est
filiforme et donn\'{e}e par le produit%
\[
\left\{ 
\begin{array}{l}
\mu (e_{1},e_{2})=e_{3} \\ 
\mu (e_{1},e_{3})=e_{4}%
\end{array}%
\right. 
\]%
dans la base $\left\{ e_{1},e_{2},e_{3},e_{4}\right\} .$

\item En dimensions $5$ et $6$, il n'existe toujours qu'un nombre fini de
classes d'isomorphie d'alg\`{e}bres de Lie nilpotentes. On pourra consulter
cette liste sur le site%
\[
\frame{ \ \ $%
\begin{array}{c}
\\ 
\ \ \ \ \ \text{http://www.math.uha.fr/\symbol{126}algebre\ \ }\ \  \\ 
\\ 
\end{array}%
$\ \ } 
\]

\item En dimension $7$, (et au del\`{a}), il existe une infinit\'{e} de
classes d'isomorphie d'alg\`{e}bres de Lie nilpotentes.\ En dimension $7$,
nous avons $6$ familles \`{a} un param\`{e}tre d'alg\`{e}bres non
isomorphes ainsi qu' une grande famille discr\`{e}te.\ Cette liste peut
aussi \^{e}tre consult\'{e}ee sur le site

\[
\frame{ \ \ $%
\begin{array}{c}
\\ 
\ \ \ \ \ \text{http://www.math.uha.fr/\symbol{126}algebre\ \ }\ \  \\ 
\\ 
\end{array}%
$\ \ } 
\]

\item Pour les dimensions sup\'{e}rieures ou \'{e}gales \`{a} $8$, seules
des classifications partielles sont connues. Par exemple les alg\`{e}bres
filiformes sont class\'{e}es jusqu'en dimension $\ 11$. Une alg\`ebre filiforme est
dite gradu\'ee si elle est munie d'une d\'erivation diagonalisable. La graduation
dans ce cas est donn\'ee par les espaces raciniens. Ces alg\`{e}bres
filiformes gradu\'{e}es sont class\'{e}es (ou plut\^ot d\'{e}crites) en toute
dimension.\ Pour approcher une classification g\'{e}n\'{e}rale, nous pouvons
utiliser l'invariant suivant, appel\'{e}e suite caract\'{e}ristique (ou
parfois invariant de Goze dans certains travaux trop \'{e}logieux pour
l'auteur) ainsi d\'{e}finie
\end{itemize}

\noindent%

\begin{defi}
Soit g une alg\`ebre de Lie nilpotente complexe de dimension $n$. Pour tout $X\in \g$, soit $c(X)$
 la suite d\'ecroissante des invariants de similitude de l'op\'erateur nilpotent $adX$
 (la suite d\'ecroissante des dimensions des blocs de Jordan). La suite caract\'eristique de $\g$ est la suite

$$c(\g)=max\{c(X), \ X \in \g-C^1(\g)\}$$

l'ordre sur les suites \'etant l'ordre lexicographique.
\end{defi}

\bigskip

En particulier, la suite caract\'{e}ristique d'une alg\`{e}bre filiforme de
dimension $n$ est $(n-1,1)$ et cette suite caract\'{e}rise la classe des alg%
\`{e}bres filiformes, la suite caract\'{e}ristique d'une alg\`{e}bre ab\'{e}%
lienne est $(1,\cdots ,1)$ et celle de \medskip l'alg\`{e}bre de Heisenberg $%
\mathfrak{h}_{2p+1}$ de dimension $(2p+1)$\ est $(2,1,\cdots ,1)$ cette
suite caract\'{e}risant aussi cette alg\`{e}bre. \ La classification des alg%
\`{e}bres de Lie nilpotentes gradu\'{e}es dont la suite caract\'{e}ristique
est $\ (n-2,1,1)$ est \'{e}galement connue \ \cite{G.K}, \cite{Ga} et \cite{Gom}.

\medskip

\subsection{Alg\`{e}bres de Lie r\'{e}solubles}

Soit $\mathfrak{g}$ une alg\`{e}bre de Lie.\ Consid\'{e}rons la suite d'id%
\'{e}aux suivante : 
\[
\left\{ 
\begin{array}{l}
{\mathcal{D}}^{0}(\mathfrak{g})=\mathfrak{g} \\ 
{\mathcal{D}}^{1}(\mathfrak{g})=\mu (\mathfrak{g},\mathfrak{g})=\mathcal{C}%
^{1}(\mathfrak{g}) \\ 
{\mathcal{D}}^{p}(\mathfrak{g})=\mu ({\mathcal{D}}^{p-1}(\mathfrak{g}),{%
\mathcal{D}}^{p-1}(\mathfrak{g}))\ \text{ pour }\ p>2.%
\end{array}%
\right. 
\]%
Comme 
\[
{\mathcal{D}}^{p}(\mathfrak{g})\subseteq  {\mathcal{D}}^{p-1}(\mathfrak{g}),\
p>0 
\]%
cette suite est d\'{e}croissante.

\medskip

\noindent \noindent 

\begin{defi}
Une alg\`ebre de Lie  $\g$ est dite r\'esoluble s'il existe un entier $k$ tel que
$${\mathcal{D}}^k(\g)=\{0\}.$$
\end{defi}

\medskip

\noindent \textbf{Exemples.}

1. Toute alg\`{e}bre nilpotente est r\'{e}soluble car 
\[
{\mathcal{D}}^{p}(\mathfrak{g})\subset {\mathcal{C}}^{p}(\mathfrak{g}). 
\]

2. Toute alg\`{e}bre de Lie de dimension $2$ est r\'{e}soluble.\ En effet,
soit elle est ab\'{e}lienne, soit isomorphe \`{a} l'alg\`{e}bre donn\'{e}e
par%
\[
\mu (e_{1},e_{2})=e_{2} 
\]%
et cette alg\`{e}bre est r\'{e}soluble.

\bigskip

Nous pouvons construire des alg\`{e}bres de Lie r\'{e}solubles \`{a} partir
d'alg\`{e}bres nilpotentes par extension par d\'{e}rivations.\ Ce proc\'{e}d%
\'{e} est ainsi d\'{e}fini:

\medskip

\noindent \noindent 

\begin{defi}
Une d\'erivation de l'alg\`ebre de Lie $\g$ est un endomorphisme lin\'eaire $f$ satisfaisant
$$\mu(f(X),Y)+\mu(X,f(Y))=f(\mu(X,Y))$$
pour tout $X,Y \in \g$.
\end{defi}

\medskip

\noindent Par exemple, les endomorphismes $adX$ donn\'{e}s par 
\[
adX(Y)=\mu (X,Y) 
\]%
sont des d\'{e}rivations (appel\'{e}es d\'{e}rivations int\'{e}rieures).
Dans une alg\`{e}bre de Lie semi-simple, toutes les d\'{e}rivations sont int%
\'{e}rieures. Notons que toute d\'{e}rivation int\'{e}rieure est singuli\`{e}%
re c'est-\`a-dire non inversible car $X\in Ker(adX)$. 

\noindent 

\begin{proposition}
Une alg\`ebre de Lie munie d'une d\'erivation r\'eguli\`ere (ou inversible) est nilpotente.
\end{proposition}

Ceci \'{e}tant, soit $\mathfrak{g}$ une alg\`{e}bre de Lie nilpotente de
dimension $n$ et soit $f$ une d\'{e}rivation non int\'{e}rieure. Consid\'{e}%
rons l'espace vectoriel $\mathfrak{g}^{\prime }=\mathfrak{g}\oplus \mathbb{C}
$ de dimension $n+1$ et notons par $e_{n+1}$ une base du compl\'{e}mentaire $%
\mathbb{C}$ dans cette extension. D\'{e}finissons une multiplication dans $%
\mathfrak{g}^{\prime }$ par 
\[
\left\{ 
\begin{array}{l}
\mu ^{\prime }(X,Y)=\mu (X,Y),\ X,Y\in \mathfrak{g} \\ 
\mu ^{\prime }(X,e_{n+1})=f(X),\ X\in \mathfrak{g}%
\end{array}%
\right. 
\]%
Ce produit munit $\mathfrak{g}^{\prime }$ d'une structure d'alg\`{e}bre de
Lie de dimension $n+1$.\ Cette alg\`{e}bre est r\'{e}soluble d\`{e}s que $f$
n'est pas une d\'{e}rivation nilpotente.

\bigskip \newpage

\part{D\'{e}formations et rigidit\'{e} des alg\`{e}bres de Lie complexes de
dimension finie}

\section{\protect\medskip Sur les origines du probl\`{e}me}

L'alg\`{e}bre de Lie de Poincar\'{e}, qui se d\'{e}duit de l'\'{e}tude des
sym\'{e}tries d'un espace vectoriel pseudo-euclidien, incluant les
rotations, pseudo-rotations et translations est un produit semi-direct de
l'alg\`{e}bre de type $so(p,q)$ qui correspondant \`{a} l'alg\`{e}bre des
matrices antisym\'{e}triques par rapport \`{a} une forme quadratique non d%
\'{e}g\'{e}n\'{e}r\'{e}e de signature $(p,q)$ par une alg\`{e}bre ab\'{e}%
lienne (engendr\'{e}e par les op\'{e}rateurs moments). Cette alg\`{e}bre
n'est pas semi-simple ni m\^{e}me r\'{e}ductive, c'est-\`{a}-dire produit
direct d'une alg\`{e}bre semi-simple par un centre ab\'{e}lien. Toute la th%
\'{e}orie classique des repr\'{e}sentations des alg\`{e}bres simples ou
semi-simples ne peut donc s'appliquer aux alg\`{e}bres de Poincar\'{e}. Il
existe un outil tr\`{e}s utile, appel\'{e} contraction, qui permet de relier
ces alg\`{e}bres, ou plus g\'{e}n\'{e}ralement des alg\`{e}bres de Lie non r%
\'{e}ductives \`{a} des alg\`{e}bres de Lie r\'{e}ductives et exploiter la th%
\'{e}orie des repr\'{e}sentations de ces derni\`{e}res. C'est ainsi que l'alg%
\`{e}bre de Poincar\'{e} correspondant aux sym\'{e}tries d'un espace \`{a} $%
4 $ dimensions, et \`{a} une forme quadratique de signature $(3,1)$ sera reli%
\'{e}e \`{a} l'alg\`{e}bre simple $so(5)$. Cette proc\'{e}dure de
contraction permet donc de remonter des propri\'{e}t\'{e}s \`{a} des alg\`{e}%
bres \`{a} partir d'autres qui ne leur sont pas isomorphes. Dans le cas de
l'alg\`{e}bre de Poincar\'{e} classique, la contraction est obtenue en consid%
\'{e}rant la c\'{e}l\'{e}rit\'{e} comme un param\`{e}tre et en faisant
tendre ce param\`{e}tre vers l'infini.\ Conceptuellement, ceci signifie que
dans l'ensemble des constantes de structure des alg\`{e}bres de Lie d'une
dimension donn\'{e}e, on consid\`{e}re une suite de familles de constantes
chacune de ces familles correspondant \`{a} des alg\`{e}bres de Lie
isomorphes et examiner la limite, si elle existe \'{e}ventuellement. Il nous
faut donc \'{e}tablir un cadre formel et topologique pour cette op\'{e}%
ration. On se fixe une dimension $n$, une base commune \`{a} tous les epaces
vectoriels sous-jacents aux alg\`{e}bres de Lie de dimension $n,$ et sans
restriction aucune supposer que tous ces espaces vectoriels sont $\mathbb{C}%
^{n}$. Dans la base fix\'{e}e, chaque alg\`{e}bre de Lie est d\'{e}termin%
\'{e}e par ses constantes de structures et l'ensemble des alg\`{e}bres de
Lie est donc assimil\'{e} \`{a} l'ensemble des constantes de structures.\
Cet ensemble est naturellement muni d'une structure de vari\'{e}t\'{e} alg%
\'{e}brique plong\'{e}e dans un espace vectoriel. Deux topologies
apparaissent naturellement sur cet ensemble, la topologie de Zariski, la
plus naturelle, et la topologie m\'{e}trique induite par l'espace vectoriel,
qui est plus fine. Cette vari\'{e}t\'{e} admet une partition en classes,
chaque classe correspond \`{a} des alg\`{e}bres isomorphes.\ Dans ce cadre,
une contraction appara\^{\i}t comme un point adh\'{e}rent \`{a} une classe.
Cette approche am\`{e}ne \`{a} une \'{e}tude topologique plus pr\'{e}cise.\
En particulier la description d'un voisinage (pour les deux topologies)
d'une alg\`{e}bre donn\'{e}e.\ Ici appara\^{\i}t la notion de d\'{e}%
formations. Pour la topologie de Zariski, on peut consid\'{e}rer une d\'{e}%
formation comme un point proche g\'{e}n\'{e}rique.\ Il existe une approche
classique pour \'{e}tudier les points g\'{e}n\'{e}riques, utiliser une
extension non archim\'{e}dienne du corps de base et consid\'{e}rer les alg%
\`{e}bres de Lie \`{a} coefficients dans un anneau local ayant l'extension
comme corps de fractions. C'est ainsi que si nous prenons comme corps, le
corps des fractions de l'anneau des s\'eries formelles, nous obtenons les d%
\'{e}formations de Gerstenhaber. Nous allons donc d\'{e}finir cette notion
de d\'{e}formation en ayant comme arri\`{e}re pens\'{e}e perp\'{e}tuelle
l'objectif topologique.

\section{La vari\'{e}t\'{e} alg\'{e}brique des alg\`{e}bres de Lie complexes
de dimension $n$}

\subsection{D\'{e}finition de $L^{n}$}

Une alg\`{e}bre de Lie complexe de dimension $n$ peut \^{e}tre consid\'{e}r%
\'{e}e comme une paire $\mathfrak{g}=(\mathbb{C}^{n},\mu )$ o\`{u} $\mu $
est le produit de l'alg\`{e}bre de Lie et $\mathbb{C}^{n}$ l'espace
vectoriel sosu-jacent \`{a} $\mathfrak{g}.$ On peut donc identifier une alg%
\`{e}bre de Lie \`{a} son produit $\mu $ et l'ensemble des alg\`{e}bres de
Lie complexes de dimension $n$ s'identifie \`{a} l'ensemble des applications
bilin\'{e}aires antisym\'{e}triques sur $\mathbb{C}^{n}$ \`{a} valeurs dans $%
\mathbb{C}^{n}$ et v\'{e}rifiant l'identit\'{e} de Jacobi. Fixons, une fois
pour toute, une base $\{e_{1},e_{2},\cdots ,e_{n}\}$ de $\mathbb{C}^{n}.$
Cela n'a que peu d'importance car nous allons \'{e}tudier les alg\`{e}bres
de Lie \`{a} isomorphisme pr\`{e}s. Les constantes de structures  de $\mu $
sont les nombres complexes $C_{ij}^{k}$ donn\'{e}s par 
\[
\mu (e_{i},e_{j})=\sum_{k=1}^{n}C_{ij}^{k}e_{k}. 
\]%
Comme la base est fix\'{e}e, le produit $\ \mu $ s'identifie \`{a} ses
constantes de structure.\ Celles-ci satisfont les \'{e}quations polynomiales
qui se d\'{e}duisent des conditions de Jacobi et de l'antisym\'{e}trie: 
\[
(1)\left\{ 
\begin{array}{l}
C_{ij}^{k}=-C_{ji}^{k}\ ,\ \ \ 1\leq i<j\leq n\ ,\ \ \ 1\leq k\leq n \\ 
\\ 
\sum_{l=1}^{n}C_{ij}^{l}C_{lk}^{s}+C_{jk}^{l}C_{li}^{s}+C_{ki}^{l}C_{jl}^{s}=0\ ,\ \ \ 1\leq i<j<k\leq n\ ,\ \ \ 1\leq s\leq n.%
\end{array}%
\right. 
\]%
Si nous notons par $L^{n}$ l'ensemble des alg\`{e}bres de Lie complexes de
dimension $n$, cet ensemble apparait donc comme une vari\'{e}t\'{e} alg\'{e}%
brique plong\'{e}e dans $\mathbb{C}^{n^{3}}$ et param\'{e}tr\'{e}e par les $%
C_{ij}^{k}$. Comme il n'y a aucune confusion possible, on pourra toujours
noter par $\mu $ les \'{e}l\'{e}ments de $L^{n}.$

\subsection{Orbite d'une alg\`{e}bre de Lie dans $L^{n}$}

Si $\mu \in L^{n}$ et si $f\in Gl(n,\mathbb{C)}$ est un isomorphisme lin\'{e}%
aire, nous avons d\'{e}fini le produit\ $\mu _{f\text{ }}$ isomorphe \`{a} $%
\mu $ par%
\[
\mu _{f}(X,Y)=f\ast \mu (X,Y)=f^{-1}(\mu (f(X),f(Y))) 
\]%
pour tout $X,Y\in \mathbb{C}^{n}.$ On d\'{e}finit ainsi une action du groupe
alg\'{e}brique $Gl(n,\mathbb{C)}$ sur $L^{n}$. On notera $\mathcal{O}(\mu )$
l'orbite du point $\mu $ relative \`{a} cette action. Consid\'{e}rons le
sous-groupe $G_{\mu }$ de $Gl(n,\mathbb{C)}$ d\'{e}fini par 
\[
G_{\mu }=\{f\in Gl(n,\mathbb{C)}\mid \text{\quad }f\ast \mu =\mu \} 
\]%
Son alg\`{e}bre de Lie est l'alg\`{e}bre de Lie des d\'{e}rivations $Der(%
\mathfrak{g})$ des d\'{e}rivations $\mathfrak{g=(\mu ,}\mathbb{C}^{n})$%
. L'orbite $\mathcal{O}(\mu )$ est donc isomorphe \`{a} l'espace homog\`{e}ne 
$Gl(n,\mathbb{C)}$/$G_{\mu }.\,$\ Cet espace est une vari\'{e}t\'{e} diff%
\'{e}rentiable de classe $\mathcal{C}^{\infty }$ et de dimension \'{e}gale 
\`{a} 
\[
\dim \mathcal{O}(\mu )=n^{2}-\dim Der(\mu ). 
\]%
On a donc montr\'{e}

\medskip

\noindent 

\begin{proposition}
L' orbite $\mathcal{O}(\mu )$ de l'alg\`ebre de Lie $\g=(\mu,\C^n)$ 
est une vari\'et\'e diff\'erentiable homog\`ene de dimension
$$
\dim \mathcal{O}(\mu )=n^{2}-\dim Der(\mu ).
$$
\end{proposition}

\subsection{Rappel sur la cohomologie de Chevalley d'une alg\`{e}bre de Lie}

\bigskip Soit $\mathfrak{g=(\mu ,}\mathbb{C}^{n})$ une alg\`{e}bre de Lie
complexe de dimension finie. Notons par $\mathcal{C}^{p}(\mathfrak{g},%
\mathfrak{g})$ l'espace vectoriel des applications $p$-lin\'{e}aires altern%
\'{e}es sur l'espace vectoriel $\mathfrak{g}$ et \`{a} valeurs dans $%
\mathfrak{g}$. Pour tout $p$, on d\'{e}finit l'endomorphisme%
\[
\delta _{p}:\mathcal{C}^{p}(\mathfrak{g},\mathfrak{g})\longrightarrow 
\mathcal{C}^{p+1}(\mathfrak{g},\mathfrak{g}) 
\]%
par%
\[
\begin{array}{ll}
\delta _{p}\Phi (Y_{1},\ldots ,Y_{p+1})= & \sum_{l=1}^{p+1}(-1)^{l+1}\mu (Y_{l},\Phi (Y_{1},\ldots ,\hat{Y}%
_{l},\ldots ,Y_{p+1})) \\ 
& +\sum_{r<s}(-1)^{r+s}\Phi (\mu (Y_{r},Y_{s}),Y_{1},\ldots ,\hat{Y}%
_{r},\ldots ,\hat{Y}_{s},\ldots ,Y_{p+1})%
\end{array}%
\]%
avec $Y_{i}\in \mathbb{C}^{n}$ et o\`{u} $\hat{Y}$ signifie que ce vecteur
n'appara\^{\i}t pas. Ces homomorphismes satisfont%
\[
\delta _{p+1}\circ \delta _{p}=0 
\]%
et la cohomologie de Chevalley de $\mathfrak{g}$ est la cohomologie $H^{\ast
}(\mathfrak{g},\mathfrak{g})$ associ\'{e}e au complexe $(\mathcal{C}^{p}(%
\mathfrak{g},\mathfrak{g}),\delta _{p})_{p\geq 0}$ c'est-\`{a}-dire%
\[
H^{p}(\mathfrak{g},\mathfrak{g})=\frac{Z^{p}(\mathfrak{g},\mathfrak{g})}{%
B^{p}(\mathfrak{g},\mathfrak{g})} 
\]%
o\`{u} $Z^{p}(\mathfrak{g},\mathfrak{g})=\left\{ \Phi \in \mathcal{C}^{p}(%
\mathfrak{g},\mathfrak{g}),\delta _{p}\Phi =0\right\} =Ker$ $\delta _{p}$
et\ $B^{p}(\mathfrak{g},\mathfrak{g})=Im\delta _{p-1}$.\ On a, en
particulier

\begin{itemize}
\item Si $X\in \mathfrak{g,}$ $\delta _{0}X=adX.$

\item Si $f\in End(\mathfrak{g),}$ $\delta _{1}f(X,Y)=\mu (f(X),Y)+\mu
(X,f(Y))-f(\mu (X,Y))$

\item Si $\varphi \in \mathcal{C}^{2}(\mathfrak{g},\mathfrak{g})$%
\[
\begin{array}{ll}
\delta _{2}\varphi (X,Y,Z)= & \mu (\varphi (X,Y),Z)+\mu (\varphi (Y,Z),X)+\mu
(\varphi (Z,X),Y)+\varphi (\mu (X,Y),Z) \\ 
& +\varphi (\mu (Y,Z),X)+\varphi (\mu (Z,X),Y).%
\end{array}%
\]
\end{itemize}

\bigskip

\subsection{Sur la g\'{e}om\'{e}trie tangente au point $\protect\mu $ \`{a} $%
L^{n}$ et $\mathcal{O(\protect\mu )}$}

Nous avons vu que l'orbite $\mathcal{O}$($\mu $) de $\mu $ est une vari\'{e}t%
\'{e} diff\'{e}rentiable plong\'{e}e dans $L^{n}$ d\'{e}finie par 
\[
\mathcal{O}(\mu )=\frac{Gl(n,\mathbb{C)}}{G_{\mu }} 
\]%
Soit $\mu ^{\prime }$ un point proche de $\mu $ dans $\mathcal{O}$($\mu $)
(pour la topologie de la vari\'{e}t\'{e}). Il existe $f\in Gl(n,\mathbb{C)}$
tel que $\mu ^{\prime }=f\ast \mu .$ Supposons que $f$ soit proche de
l'identit\'{e}, $f=Id+\varepsilon g,$ avec $g\in gl(n)$. Alors 
\begin{eqnarray*}
\mu ^{\prime }(X,Y) &=&\mu (X,Y)+\varepsilon \lbrack -g(\mu (X,Y))+\mu
(g(X),Y)+\mu (X,g(Y))] \\
&&+\varepsilon ^{2}[\mu (g(X),g(Y))-g(\mu (g(X),Y)+\mu (X,g(Y))-g\mu (X,Y)]
\end{eqnarray*}%
et 
\[
lim_{\varepsilon \rightarrow 0}\frac{\mu ^{\prime }(X,Y)-\mu (X,Y)}{%
\varepsilon }=\delta ^{1}g(X,Y) 
\]%
o\`{u} $\delta ^{1}g$ est l'op\'{e}rateur cobord associ\'{e} \`{a} la
cohomologie de Chevalley de $\mathfrak{g}$.\ 

\medskip

\noindent 

\begin{proposition}
L'espace tangent \`a l'orbite $\mathcal{O}$($\mu )$ au point  $\mu $ est l'espace
 $B^{2}(\mu ,\mu )$ des  2-cobords de la cohomologie de Chevalley de $\g=(\mu,\C^n)$.
\end{proposition}

L'espace tangent au point $\mu $ \`{a} $L^{n}$ est plus d\'{e}licat \`{a} d%
\'{e}finir car $L^{n}$ est une vari\'{e}t\'{e} alg\'{e}brique avec de
nombreux points singuliers.\ Nous pouvons toutefois d\'{e}finir les droites
tangentes. On consid\`{e}re pour cela un point $\mu _{t}=\mu +t\varphi $ $%
\in L^{n}$ o\`{u} $t$ est un param\`{e}tre complexe. Mais $\mu _{t}$ v\'{e}%
rifie la condition de Jacobi pour tout $t$ si et seulement si%
\[
\left\{ 
\begin{array}{c}
\delta ^{2}\varphi =0, \\ 
\varphi \in L^{n}.%
\end{array}%
\right. 
\]%
On a ainsi

\medskip

\noindent 

\begin{proposition}
Un droite affine $\Delta $ passant par $\mu $ est tangente en $\mu $
\`a $L^{n}$ si son vecteur directeur appartient \`a $Z^{2}(\mu ,\mu ).$
\end{proposition}

\medskip

\medskip Cette \'{e}tude g\'{e}om\'{e}trique montre qu'une \'{e}tude plus alg%
\'{e}brique de $L^{n}$ est n\'{e}cessaire.\ Ceci nous conduit \`{a} \'{e}%
tudier le sch\'{e}ma affine associ\'{e}.\ Rappelons que la notion de vari%
\'{e}t\'{e} alg\'{e}brique et de sch\'{e}ma co\"{\i}ncident lorsque ce dernier
est r\'{e}duit.\ Mais nous allons voir que pour $L^{n},$ except\'{e}es pour
les petites valeurs de $n$, le sch\'{e}ma n'est pas r\'{e}duit et donc bon
nombre de propri\'{e}t\'{e}s va s'exprimer en terme de sch\'{e}mas associ%
\'{e}s.

\noindent

\subsection{Le sch\'{e}ma $\mathcal{L}^{n}.$}

\medskip

On consid\`{e}re les $C_{ij}^{k}$ avec $1\leq i<j\leq n$ et $1\leq k\leq n$
comme des variables formelles et soit $\mathbb{C}[C_{ij}^{k}]$ l'anneau des
polyn\^{o}mes complexes en les $\frac{n^{3}-n^{2}}{2}$ variables $%
C_{ij}^{k}. $ Soit $I_{J}$ l'id\'{e}al de $\mathbb{C}[C_{ij}^{k}]$ engendr%
\'{e} par les polyn\^omes de Jacobi : 
\[
\sum_{l=1}^{n}C_{ij}^{l}C_{lk}^{s}+C_{jk}^{l}C_{li}^{s}+C_{ki}^{l}C_{jl}^{s} 
\]%
avec $1\leq i<j<k\leq n,$ et $1\leq s\leq n.$ Ces polyn\^omes sont au nombre de $n(n-1)(n-2)/6$. 
La vari\'{e}t\'{e} alg\'{e}%
brique $L^{n}$ est l'ensemble alg\'{e}brique associ\'{e} \`{a} l'id\'{e}al $%
I_{J}:$%
\[
L^{n}=V(I_{J})=\left\{ x\in \mathbb{C}^{\frac{n^{3}-n^{2}}{2}},\text{ tels
que }P(x)=0,\forall P\in I_{J}\right\} . 
\]%
Soit $rad(I_{J})$ le radical de l'id\'{e}al $I_{J}$ c'est-\`{a}-dire%
\[
rad(I_{J})=\{P\in \mathbb{C}[C_{jk}^{i}],\text{ tels que }\exists r\in 
\mathbb{N}^{\ast }\text{, }P^{r}\in I_{J}\} 
\]%
En g\'{e}n\'{e}ral, $radI_{J}\neq I_{J}$ (Rappelons que si $I$ est premier
alors $rad(I)=I$). Si $M$ est un sous-ensemble de $\mathbb{C}^{\frac{%
n^{3}-n^{2}}{2}},$ on note par $\mathfrak{i}(M)$ l'id\'{e}al de $\mathbb{C}%
[C_{jk}^{i}]$ d\'{e}fini par 
\[
\mathfrak{i}(M)=\{P\in \mathbb{C}[C_{jk}^{i}],P(x)=0\quad \forall x\in M\}. 
\]%
On a alors 
\[
\mathfrak{i}(L^{n})=\mathfrak{i}(V(I_{J}))=rad(I_{J}). 
\]%
Consid\'{e}rons l'anneau 
\[
A(L^{n})=\frac{\mathbb{C}[C_{jk}^{i}]}{I_{J}} 
\]%
qui est aussi une $\mathbb{C}-$alg\`{e}bre de type fini. Il correspond \`{a}
l'anneau des fonctions r\'{e}guli\`{e}res sur $L^{n}.$ Rappelons que si un id%
\'{e}al $I$ d'un anneau $A$ est premier, l'anneau quotient $A/I$ est un
anneau int\`{e}gre et que tout id\'{e}al maximal est premier. L'anneau
quotient est dit r\'{e}duit s'il ne contient pas d'\'{e}l\'{e}ments
nilpotents.\ Comme en g\'{e}n\'{e}ral $radI_{J}\neq I_{J}$, l'alg\`{e}bre $%
A(L^{n})$ n'est pas r\'{e}duite.

L'alg\`{e}bre affine $\Gamma (L^{n})$ de la vari\'{e}t\'{e} alg\'{e}brique $%
L^{n}$ est l'anneau quotient 
\[
\Gamma (L^{n})=\frac{\mathbb{C}[C_{jk}^{i}]}{\mathfrak{i}(L^{n})}. 
\]%
Comme on a 
\[
\mathfrak{i}(L^{n})=rad\mathfrak{i}(L^{n}), 
\]%
on en d\'{e}duit que $\Gamma (L^{n})$ est toujours r\'{e}duite.

Le spectre maximal $Spm(A(L^{n}))$ \ de $A(L^{n})$ est l'ensemble des id\'{e}%
aux maximaux de l'alg\`{e}bre $A(L^{n}).$ Pour chaque $x\in L^{n},$
l'ensemble%
\[
\mathcal{M}_{x}=\left\{ f\in A(L^{n}),f(x)=0\right\} 
\]%
est un id\'{e}al maximal de $A(L^{n}).$ L'application $x\in L^{n}\rightarrow 
$ $\mathcal{M}_{x}$ d\'{e}finit une bijection entre $L^{n}$ et $%
Spm(A(L^{n})) $ 
\[
Spm(A(L^{n}))\sim L^{n}. 
\]%
On rajoute en quelque sorte des nouveaux points \`{a} $L^{n}$ en consid\'{e}%
rant une extension de $Spm(A(L^{n}))$ \`{a} savoir%
\[
Spec(A(L^{n}))=\left\{ I,\text{ o\`{u} }I\text{ est un id\'{e}al premier
propre de }A(L^{n})\right\} 
\]%
Comme tout id\'{e}al maximal est premier, on a bien $Spm(A(L^{n}))\subset
Spec(A(L^{n})).$ On munit $Spec(A(L^{n}))$ de la topologie de Zariski : Les
ensembles ferm\'{e}s sont d\'{e}finis ainsi : soit $\mathfrak{a}$ un id\'{e}%
al de $A(L^{n}),$ on consid\`{e}re le sous ensemble $V(\mathfrak{a})$ de $%
Spec(A(L^{n}))$ constitu\'{e} des id\'{e}aux premiers de $A(L^{n})$
contenant $\mathfrak{a}$. Ces ensembles forment la famille de ferm\'{e}s de
la topologie$.$ On en d\'{e}duit que les ensembles%
\[
Spec(A(L^{n}))_{f}=\left\{ I\in Spec(A(L^{n})),f\notin I\right\} 
\]%
o\`{u} $f\in A(L^{n})$ forment une base d'ouverts. Il existe un faisceau de
fonctions $\mathcal{O}_{Spec(A(L^{n}))}$ sur $Spec(A(L^{n}))$ \`{a} valeurs
dans $\mathbb{C}$ tel que pour tout $f\in A(L^{n})$, on ait 
\[
\Gamma (D(f),\mathcal{O}_{Spec(A(L^{n}))})=A(L^{n})_{f} 
\]%
o\`{u} $A(L^{n})_{f}$ est l'anneau des fonctions 
\[
x\rightarrow \frac{g(x)}{f(x)^{n}} 
\]%
pour $x\in D(f)=\left\{ y,f(y)\neq 0\right\} $ et $g$ $\in A(L^{n}).$ En
particulier on a 
\[
\Gamma (Spec(A(L^{n})),\mathcal{O}_{Spec(A(L^{n}))})=A(L^{n}) 
\]%
Le sch\'{e}ma affine de l'anneau $A(L^{n})$ est l'espace $(Spec(A(L^{n})),%
\mathcal{O}_{Spec(A(L^{n}))})$ not\'e aussi $Spec(A(L^{n}))$ ou $\mathcal{L}%
^{n}.$

\bigskip

Comme tout point $x$ de $L^{n}$ correspond bijectivement \`{a} un point de $%
Spec(A(L^{n}))$ donn\'{e} par un id\'{e}al maximal, nous pouvons d\'{e}finir
un espace tangent en $x$ \`{a} $Spec(A(L^{n})).$ On consid\`{e}re pour cela
une d\'{e}formation infinit\'{e}simale de l'alg\`{e}bre $\Gamma
(Spm(A(L^{n})),\mathcal{O}_{Spm(A(L^{n}))})$ en ce point. Si $%
F_{1},...,F_{N} $ avec $N=\frac{1}{6}n(n-1)(n-2)$ sont les polyn\^{o}mes de
Jacobi, alors l'espace tangent en $x$ au sch\'{e}ma $\mathcal{L}^{n}$ est 
\[
T_{x}(\mathcal{L}^{n})=Ker\,d_{x}(F_{1},..,F_{N}) 
\]%
o\`{u} $d_{x}$ d\'{e}signe la matrice jacobienne des $F_{i}$ au point $x$.
D'apr\`{e}s la d\'{e}finition de la cohomologie de Chevalley de l'alg\`{e}%
bre de Lie $\mathfrak{g}$ correspondant au point $x,$ on a :

\medskip

\noindent 

\begin{theo}
$T_x(\mathcal{L}^{n})=Z^{2}(\frak{g},\frak{g})$
\end{theo}

\bigskip

\section{Contractions d' alg\`{e}bres de Lie}

\medskip

Parfois le mot contraction est remplac\'{e} par d\'{e}g\'{e}n\'{e}r\'{e}%
scence ou d\'{e}g\'{e}n\'{e}ration pour les anglophiles.\ Dans toute cette
partie, nous confondrons sans retenue une alg\`{e}bre de Lie complexe de
dimension $n$, sa multiplication et le point correspondant de la vari\'{e}t%
\'{e} $L^{n}$ voire m\^{e}me du sch\'{e}ma affine. Rappelons que $L^{n}$ est
fibr\'{e}e par les orbites relatives \`{a} l'action du groupe alg\'{e}brique 
$GL(n,\mathbb{C})$ et nous avons not\'{e}e par $\mathcal{O(\mu )}$ l'orbite
du point $\mu .$

\subsection{D\'{e}finition d'une contraction d'alg\`{e}bre de Lie et exemples%
}

Soit $g=(\mu ,\mathbb{C}^{n})$ une alg\`{e}bre de Lie complexe de dimension $%
n$.

\medskip

\noindent \noindent%

\begin{defi}
Une alg\`ebre de Lie $\g _0=(\mu _{0}, \C^n)$, $\mu _0 \in L^{n}$  est appel\'ee contraction de 
$\g =(\mu, \C^n)$ 
si $\mu _{0}\in \overline{\mathcal{O}(\mu )}.$
\end{defi}

\medskip

La notion historique de contraction est celle de Segal.\ Elle est ainsi d%
\'{e}finie : Soit $\{f_{p}\}$ une suite dans $GL(n,\mathbb{C})$. On d\'{e}%
duit la suite $\{\mu _{p}\}$ dans $L^{n}$ en posant%
\[
\mu _{p}=f_{p}\ast \mu 
\]

Si cette suite admet une limite$\ \mu _{0}$ dans l'espace vectoriel des
applications bilin\'{e}aires altern\'{e}es, alors $\mu _{0}\in L^{n}$ et $%
\mu _{0}$ est appel\'{e}e une contraction de $\mu $. Le lien entre ces deux
notions est donn\'{e} par la proposition suivante

\medskip

\noindent 

\begin{proposition}
Pour tout $\mu \in L^n$, l'adh\'erence de Zariski  $\overline{\mathcal{O}(\mu )}$ de l'orbite ${\mathcal{O}(\mu )}$ 
est \'equivalente \`a l'adh\'erence pour la topologie m\'etrique induite
$$\overline{\mathcal{O}(\mu )}=\overline{\mathcal{O}(\mu )}^d.$$
\end{proposition}

\medskip

\noindent \textit{D\'{e}monstration.} Ceci repose essentiellement sur le
fait que le corps de base est $\mathbb{C}$.

\medskip

\noindent 

\begin{coro}
Toute contraction de $\mu \in L^n$ est obtenue  par une contraction de Segal.
\end{coro}

\medskip

\textbf{Exemples.}

\begin{itemize}
\item Le cas ab\'{e}lien.\ Toute alg\`{e}bre de Lie de dimension $n$ se
contracte sur l'alg\`{e}bre de Lie ab\'{e}lienne de dimension $n.\ $En effet
si $\mu $ est donn\'{e}e dans la base $\{e_{1},\ldots ,e_{n}\}$ par $\mu
(e_{i},e_{j})=\sum C_{ij}^{k}e_{k},$ on consid\`{e}re alors l'isomorphisme
donn\'{e} par $f_{\varepsilon }(X_{i})=\varepsilon X_{i},\varepsilon \neq 0.$
Alors la multiplication $\mu _{\varepsilon }=f_{\varepsilon }\ast \mu $ v%
\'{e}rifie $\mu _{\varepsilon }(X_{i},X_{j})=\sum \varepsilon
C_{ij}^{k}X_{k} $ et $lim_{\varepsilon \rightarrow 0}\mu _{\varepsilon }$
existe et co\"{\i}ncide avec le produit de l'alg\`{e}bre ab\'{e}lienne.

\item L'ag\`{e}bre de Poincar\'{e}. Le groupe de Lie de Poincar\'{e} est le
groupe des isom\'{e}tries de l'espace de Minkowski. C'est un groupe de Lie
non compact de dimension $10$. C'est en fait le groupe affine associ\'{e} au
groupe de Lorentz.\ C'est donc un produit semi-direct des translations et
des transformations de Lorentz. L'alg\`{e}bre de Poincar\'{e} est l'alg\`{e}%
bre de Lie du groupe de Poincar\'{e}. On note classiquement par $P$ les g%
\'{e}n\'{e}rateurs des translations et par $M$ ceux des transformations de
Lorentz.\ La multiplication de cette alg\`{e}bre, que nous allons noter
aussi classiquement par le crochet est donn\'{e}e par%
\[
\left\{ 
\begin{array}{l}
\lbrack P_{\mu },P_{\nu }]=0, \\ 
\lbrack M_{\mu \nu },P_{\rho }]=\eta _{\mu \rho }P_{\nu }-\eta _{\nu \rho
}P_{\mu } \\ 
\lbrack M_{\mu \nu },M_{\rho \sigma }]=\eta _{\mu \rho }M_{\nu \sigma }-\eta
_{\mu \sigma }M_{\nu \rho }-\eta _{\nu \rho }M_{\mu \sigma }+\eta _{\nu \rho
}M_{\mu \sigma }%
\end{array}%
\right. 
\]%
o\`{u} $\eta $ est la m\'{e}trique de Minkowski et dans cette notation $%
M_{ab}=-M_{ba}.$ Ceci \'{e}tant consid\'{e}rons l'alg\`{e}bre de Lie simple $%
so(5)$ dont les \'{e}l\'{e}ments sont les matrices antisym\'{e}triques
d'ordre $5$.\ Elle est de dimension $10.~$Nous allons montrer que l$^{\prime
}$a$\lg $\`{e}bre de Poincar\'{e} est une contraction de $so(5)$. L'alg\`{e}%
bre $so(5)$ a trois formes r\'{e}elles non \'{e}quivalentes (une forme r\'{e}%
elle est une alg\`{e}bre de Lie r\'{e}elle dont la complexifi\'{e}e est
isomorphe \`{a} l'alg\`{e}bre complexe $so(5)).$ Ces formes sont l'alg\`{e}%
bre de de Sitter, l'alg\`{e}bre anti de Sitter et la forme compacte, l'alg%
\`{e}bre r\'{e}elle $so(5,\mathbb{R)}$. Elle correspondent respectivement
aux signatures de la forme de Killing Cartan suivantes : $(-++++),(-+++-)$
et $(+++++)$. Si $M_{ab}$ sont les g\'{e}n\'{e}rateurs de l'alg\`{e}bre de
de Sitter (et anti de Sitter), $0\leq a,b\leq 4$, alors les \'{e}l\'{e}ments 
$M_{\mu \nu }$ pour $0\leq \mu ,\nu \leq 3$ engendrent une sous-alg\`{e}bre
isomorphe \`{a} l'alg\`{e}bre de Lorentz $so(3,1).$ Soit $\varepsilon $ un r%
\'{e}el positif.\ Le changement de base%
\[
P_{m}^{\varepsilon }=\varepsilon M_{m4} 
\]%
implique%
\[
\lbrack P_{m}^{\varepsilon },P_{n}^{\varepsilon }]=\varepsilon
^{2}[M_{m4},M_{n4}] 
\]%
et%
\[
\lbrack M_{lm},P_{n}^{\varepsilon }]=\eta _{mn}P_{l}^{\varepsilon }-\eta
_{\ln }P_{m}^{\varepsilon } 
\]%
Lorsque $\varepsilon \rightarrow 0$, alors%
\[
\lbrack P_{m}^{\varepsilon },P_{n}^{\varepsilon }]=\varepsilon
^{2}[M_{m4},M_{n4}]\rightarrow 0 
\]%
et on obtient les crochets de l'alg\`{e}bre de Poincar\'{e}.
\end{itemize}

\medskip

\subsection{Alg\`{e}bres de Lie de contact, Alg\`{e}bres de Lie frob\'{e}%
niusiennes.}

Consid\'{e}rons l'ouvert (de Zariski) $\mathcal{C}_{2p+1}$ de $L^{2p+1\text{ 
}}$ dont les \'{e}l\'{e}ments sont les alg\`{e}bres de Lie de dimension $%
(2p+1)$ munie d'une forme de contact.\ Rappelons qu'une forme de contact est
un \'{e}l\'{e}ment $\omega \in \mathfrak{g}^{\ast }$ (le dual vectoriel de $%
\mathfrak{g}$) satisfaisant 
\[
\omega \wedge (d\omega )^{p}\neq 0 
\]%
o\`u $d \omega$ est la $2$-forme ext\'erieure sur $\mathfrak{g}$ donn\'ee par
$$ d\omega (X,Y)=\omega[X,Y]$$
pour tous $X,Y \in \mathfrak{g}$.
 Il existe alors une base $\left\{ X_{1},X_{2},...,X_{2p+1}\right\} $ de $%
\mathfrak{g}$ \ telle que la base duale $\left\{ \omega _{1},\omega
_{2}...,\omega _{2p+1}\right\} $ satisfasse $\omega =\omega _{1}$ et 
\[
d\omega _{1}=\omega _{2}\wedge \omega _{3}+\omega _{4}\wedge \omega
_{5}+...+\omega _{2p}\wedge \omega _{2p+1}. 
\]%
En particulier les constantes de structure de $\mathfrak{g}$ relatives \`{a}
cette base v\'{e}rifent 
\[
C_{23}^{1}=C_{34}^{1}=...=C_{2p2p+1}^{1}=1. 
\]%
Consid\'{e}rons l'isomorphisme $f_{\varepsilon }$ de $\mathbb{C}^{2p+1}$ donn%
\'{e} par 
\[
f_{\varepsilon }(X_{1})=\varepsilon ^{2}X_{1},\text{ }f_{\varepsilon
}(X_{i})=\varepsilon X_{i}\quad i=2,...,2p+1. 
\]%
Les structures de constantes $D_{ij}^{k}$ de $\mu _{\varepsilon
}=f_{\varepsilon }\ast \mu $ par rapport \`{a} la base $\left\{
X_{1},X_{2},...,X_{2p+1}\right\} $ sont 
\[
\left\{ 
\begin{array}{l}
D_{23}^{1}=D_{34}^{1}=...=D_{2p2p+1}^{1}=1 \\ 
D_{ij}^{k}=\varepsilon C_{ij}^{k}\text{ pour tous les autres indices.}%
\end{array}%
\right. 
\]%
Ainsi $lim_{\varepsilon \rightarrow 0}\mu _{\varepsilon }$ existe et cette
limite correspond \`{a} la multiplication de l'alg\`{e}bre de Heisenberg $%
\mathfrak{h}_{2p+1}$ \ de dimension $2p+1.$ Ainsi $\mathfrak{h}_{2p+1}\in 
\overline{\mathcal{C}_{2p+1}}.$

\medskip

\noindent 

\begin{proposition}
Toute alg\`ebre de Lie  de dimension $2p+1$ munie d'une forme de contact se contracte sur l'alg\`ebre de Lie de Heisenberg
$\frak{h}_{2p+1}$ de dimension $2p+1$. De plus, toute alg\`ebre de Lie qui se contracte sur $\frak{h}_p$
admet une forme de contact.
\end{proposition}

\medskip

\bigskip

Le cas des alg\`{e}bres frob\'{e}niusiennes est l'analogue en dimension
paire. Soit $\g$ une alg\`{e}bre de Lie de dimension $2p$. Elle est appel\'{e}%
e frob\'{e}niusienne s'il existe une forme lin\'{e}aire non nulle $\omega
\in {\mathfrak{g}}^{\ast }$ telle que%
\[
\lbrack d\omega ]^{p}\neq 0. 
\]%
Dans ce cas, la $2$-forme $\theta =d\omega $ est une forme symplectique
exacte sur $\mathfrak{g}.$

\medskip

\noindent 

\begin{theo}
\cite{Go2} Soit $\{F_{\varphi}\;|\;\varphi  \in\mathbb{C}^{p-1}\}$
la famille \`a $\left(  p-1\right)  $-param\`etres d'alg\`ebres de Lie de dimension $2p$ donn\'ee par
\[
\left\{
\begin{array}
[c]{l}d\omega_{1}=\omega_{1}\wedge\omega_{2}+\sum_{k=1}^{p-1}\omega_{2k+1}\wedge\omega_{2k+2}\\
d\omega_{2}=0\\
d\omega_{2k+1}=\varphi_{k}\omega_{2}\wedge\omega_{2k+1},\;1\leq k\leq p-1\\
d\omega_{2k+2}=-\left(  1+\varphi_{k}\right)  \omega_{2}\wedge\omega
_{2k+2},\;1\leq k\leq p-1
\end{array}
\right.
\]
o\`u $\left\{  \omega_{1},..,\omega_{2p}\right\}  $ est une base $\left(
\mathbb{C}^{2p}\right)  ^{\ast}$. La famille $\{F_{\varphi}\}$ est un mod\`ele irr\'eductible complexe pour la propriét\'e
``Il existe une forme symplectique exacte'', c'est-\`a-dire toute alg\`ebre de Lie de dimension $2p$ frob\'eniusienne se 
contracte sur un \'el\'ement de la famille, et il n'existe aucune contraction entre deux \'el\'ements distincts
de  $\{F_{\varphi}\}$.
\end{theo}

\medskip

Notons que la famille $F_{\varphi }$ est gradu\'{e}e : si $\left\{
X_{1},..,X_{2p}\right\} $ est la base duale de $\left\{ \omega _{1},..,\omega
_{2p}\right\} $, alors $F_{\varphi }=\left( F_{\varphi }\right) _{0}\oplus
\left( F_{\varphi }\right) _{1}\oplus \left( F_{\varphi }\right) _{2}$, o%
\`{u} $\left( F_{\varphi }\right) _{0}=\mathbb{C}{X_{2}},\left( F_{\varphi
}\right) _{1}=\sum_{k=3}^{2p}\mathbb{C}{X_{k}}$ et $\left( F_{\varphi
}\right) _{2}=\mathbb{C}{X_{1}}$. Cette d\'{e}composition se r\'{e}v\`{e}le
int\'{e}ressante lors de calcul de cohomologie des alg\`{e}bres frob\'{e}%
niusiennes.

\bigskip

\noindent \textit{D\'{e}monstration.} Soit $\mathfrak{g}$ une alg\`{e}bre de
Lie frob\'{e}niusienne de dimension $2p.$ Il existe une base $%
\{X_{1},...,X_{2p}\}$ de $\mathfrak{g}$ telle que la base duale $\{\omega
_{1},...,\omega _{2p}\}$ v\'{e}rifie 
\[
d\omega _{1}=\omega _{1}\wedge \omega _{2}+...+\omega _{2p-1}\wedge \omega
_{2p}, 
\]%
c'est-\`{a}-dire nous pouvons supposer $\omega _{1}$ frob\'{e}niusienne.
Consid\'{e}rons la famille de changement de bases : 
\[
f_{\epsilon }(X_{1})=\epsilon ^{2}X_{1},\ f_{\epsilon }(X_{2})=X_{2},\
f_{\epsilon }(X_{i})=\epsilon X_{i},\ i=3,...,2p 
\]%
o\`{u} $\varepsilon $ est un param\`{e}tre complexe. Lorsque $\varepsilon
\rightarrow 0,$ nous obtenons des contractions de $\mathfrak{g}$ dont les 
\'{e}quations de Maurer-Cartan sont 
\[
\left\{ 
\begin{array}{l}
d\omega _{1}=\omega _{1}\wedge \omega _{2}+...+\omega _{2p-1}\wedge \omega
_{2p}, \\ 
d\omega _{2}=0, \\ 
d\omega _{3}=C_{23}^{3}\omega _{2}\wedge \omega _{3}+C_{24}^{3}\omega
_{2}\wedge \omega _{4}+...+C_{22p-1}^{3}\omega _{2}\wedge \omega
_{2p-1}+C_{22p}^{3}\omega _{2}\wedge \omega _{2p}, \\ 
d\omega _{4}=C_{23}^{4}\omega _{2}\wedge \omega _{3}+(-1-C_{23}^{3})\omega
_{2}\wedge \omega _{4}+...+C_{22p-1}^{4}\omega _{2}\wedge \omega
_{2p-1}+C_{22p}^{4}\omega _{2}\wedge \omega _{2p}, \\ 
.... \\ 
d\omega _{2p-1}=C_{22p}^{4}\omega _{2}\wedge \omega _{3}+C_{2p}^{3}\omega
_{2}\wedge \omega _{4}+...+C_{22p-1}^{2p-1}\omega _{2}\wedge \omega
_{2p-1}+C_{22p}^{2p-1}\omega _{2}\wedge \omega _{2p}, \\ 
d\omega _{2p}=C_{22p-1}^{4}\omega _{2}\wedge \omega _{3}+C_{22p-1}^{3}\omega
_{2}\wedge \omega _{4}+...+C_{22p-1}^{2p}\omega _{2}\wedge \omega
_{2p-1}+(-1-C_{22p-1}^{2p-1})\omega _{2}\wedge \omega _{2p}.%
\end{array}%
\right. 
\]%
Le reste de la d\'{e}monstration consiste \`{a} r\'{e}duire l'op\'{e}rateur $%
\psi $ d\'{e}fini comme la restriction de l'op\'{e}rateur adjoint $adX_{2}$
au sous-espace invariant $F$ engendr\'{e} par $\{X_{3},...,X_{2p}\}$. On v%
\'{e}rifie directement les propri\'{e}t\'{e}s suivantes :

- Si $\alpha $ et $\beta $ sont des valeurs propres de $\psi $ telles que $%
\alpha \neq -1-\beta $, \ alors les espaces propres correspondants $%
F_{\alpha }$ et $F_{\beta }$ v\'{e}rifient $[F_{\alpha },F_{\beta }]=0.$

- Si $\alpha $ est une valeur propre de of $\psi $ diff\'{e}rente de $-\frac{%
1}{2}$, alors pour tout $X$ et $Y$ $\in F_{\alpha }$, on a $[X,Y]=0.$

- Si $\alpha $ est valeur propre de $\psi $, alors $-1-\alpha $ est aussi
une valeur propre de $\psi .$

- Les multiplicit\'{e}s des valeurs propres $\alpha $ \ et $-1-\alpha $ sont 
\'{e}gales.

- La suite ordonn\'{e}e des blocs de Jordan des valeurs propres $\alpha $ et 
$-1-\alpha $ sont les m\^{e}mes.

\noindent A partir de ces remarques, on peut trouver une base de Jordan de $%
\psi $ telle que la matrice $\psi $ restreint au sous-espace invariant $%
C_{\alpha }\oplus C_{-1-\alpha }$ o\`{u} $C_{\lambda }$ est le sous-espace
caract\'{e}ristique associ\'{e} \`{a} la valeur propre $\lambda $ soit de la
forme : 
\[
\left( 
\begin{array}{lllllll}
\alpha & 0 & 1 & 0 & 0 & 0 & ... \\ 
0 & -1-\alpha & 0 & 0 & 0 & 0 & ... \\ 
0 & 0 & \alpha & 0 & 1 & 0 & ... \\ 
0 & -1 & 0 & -1-\alpha & 0 & 0 & ... \\ 
0 & 0 & 0 & 0 & \alpha & 0 & ... \\ 
0 & 0 & 0 & -1 & 0 & -1-\alpha & ... \\ 
&  &  &  &  &  & 
\end{array}%
\right) . 
\]%
Ainsi les valeurs propres et la dimension de blocs de Jordan correspondants
classent les \'{e}l\'{e}ments de la famille du mod\`{e}le frob\'{e}niusien.

\medskip

\textbf{Remarque : }Dans \cite{Go2} on trouve \'{e}galement le mod\`{e}le frob%
\'{e}niusien dans le cas r\'{e}el.

\noindent

\subsection{Contractions d'In\"{o}n\"{u}-Wigner}

Pour d\'{e}finir les contractions d'In\"{o}n\"{u}-Wigner, on consid\`{e}re la
famille \`{a} un param\`{e}tre d'isomorphismes $\{f_{\epsilon }\}$ de $GL(n,%
\mathbb{C})$ de la forme%
\[
f_{\epsilon }=f_{1}+\epsilon f_{2} 
\]%
o\`{u} $f_{1}\in gl(n,\mathbb{C})$ est un op\'{e}rateur singulier, c'est-%
\`{a}-dire $det(f_{1})=0$ et $f_{2}\in GL(n,\mathbb{C})$. Les matrices de ces
applications lin\'{e}aires peuvent se r\'{e}duire simultan\'{e}ment sous la
forme 
\[
f_{1}=\left( 
\begin{array}{ll}
Id_{r} & 0 \\ 
0 & 0%
\end{array}%
\right) \ ,\ f_{2}=\left( 
\begin{array}{ll}
v & 0 \\ 
0 & Id_{n-r}%
\end{array}%
\right) 
\]%
o\`{u} $rang(f_{1})=rang(v)=r$.

\noindent Les contractions associ\'{e}es \`{a} de telles familles
d'isomorphismes sont appel\'{e}es les contractions d'In\"{o}n\"{u}-Wigner.\
Elles permettent de contracter une alg\`{e}bre de Lie donn\'{e}e $\mathfrak{g%
}$ dans une alg\`{e}bre de Lie $\mathfrak{g}_{0}$ en laissant invariant une
sous-alg\`{e}bre $\mathfrak{h}$ of $\mathfrak{g}$ ce qui signifie que $%
\mathfrak{h}$ est encore une sous-alg\`{e}bre de $\mathfrak{g}_{0}$. Par
exemple, l'alg\`{e}bre de Lorentz homog\`{e}ne peut se contracter, via une
contraction d'In\"{o}n\"{u}-Wigner dans l'alg\`{e}bre homog\`{e}ne de Galil\'{e}%
e. De m\^{e}me, l'alg\`{e}bre de De Sitter peut se contracter dans l'alg\`{e}%
bre de Lorentz non homog\`{e}ne.

Donnons \`{a} pr\'{e}sent une br\`{e}ve description des contractions d'Inon%
\"{u}-Wigner. Soit $\mathfrak{g}=(\mu ,\mathbb{C}^{n})$ une alg\`{e}bre de
Lie et $\mathfrak{h}$ une sous-alg\`{e}bre de Lie de $\mathfrak{g}$. Fixons
une base $\{e_{1},...,e_{n}\}$ de $\mathbb{C}^{n}$ telle que $%
\{e_{1},...,e_{p}\}$ soit une base de $\mathfrak{h}$. Alors 
\[
\mu (e_{i},e_{j})=\sum_{k=1}^{p}C_{ij}^{k}e_{k},\ \ i,j=1,...,p. 
\]%
Consid\'{e}rons la famille d'isomorphismes d'In\"{o}n\"{u}-Wigner donn\'{e}e par 
\[
\left\{ 
\begin{array}{l}
f_{\epsilon }(e_{i})=(1+\epsilon )e_{i},\ \ i=1,...,p \\ 
\\ 
f_{\epsilon }(e_{l})=\epsilon e_{l},\ \ l=p+1,...,n.%
\end{array}%
\right. 
\]%
On a ici $f_{\epsilon }=f_{1}+\epsilon f_{2}$ avec 
\[
f_{1}=\left( 
\begin{array}{ll}
Id_{p} & 0 \\ 
0 & 0%
\end{array}%
\right) \ ,\ f_{2}=\left( 
\begin{array}{ll}
Id_{p} & 0 \\ 
0 & Id_{n-p}%
\end{array}%
\right) . 
\]%
La multiplication $\mu _{\epsilon }=f_{\epsilon }\ast \mu $ s'\'{e}crit 
\[
\left\{ 
\begin{array}{l}
\mu _{\epsilon }(e_{i},e_{j})=(1+\epsilon )^{-1}\mu (e_{i},e_{j}),\ \
i,j=1,...,p \\ 
\\ 
\mu _{\epsilon }(e_{i},e_{l})=\epsilon (1+\epsilon
)^{-1}\sum_{k=1}^{p}C_{ij}^{k}e_{k}+(1+\epsilon
)^{-1}\sum_{k=p+1}^{n}C_{il}^{k}e_{k},\ \ i=1,..,p,\ l=p+1,...,n \\ 
\\ 
\mu _{\epsilon }(e_{l},e_{m})=\epsilon ^{2}(1+\epsilon
)^{-1}\sum_{k=1}^{p}C_{lm}^{k}e_{k}+\epsilon
\sum_{k=p+1}^{n}C_{lm}^{k}e_{k},\ \ l,m=p+1,...,n%
\end{array}%
\right. 
\]%
Si $\epsilon \rightarrow 0$, la suite $\{\mu _{\epsilon }\}$ a pour limite $%
\mu _{0}$ donn\'{e}e par 
\[
\left\{ 
\begin{array}{l}
\mu _{0}(e_{i},e_{j})=\mu (e_{i},e_{j}),\ \ i,j=1,...,p \\ 
\\ 
\mu _{0}(e_{i},e_{l})=\sum_{k=p+1}^{n}C_{il}^{k}e_{k},\ \ i=1,..,p,\
l=p+1,...,n \\ 
\\ 
\mu _{0}(e_{l},e_{m})=0,\ \ l,m=p+1,...,n.%
\end{array}%
\right. 
\]%
L'alg\`{e}bre $\mathfrak{g}_{0}=(\mu _{0},\mathbb{C}^{n})$ est une
contraction d'In\"{o}n\"{u}-Wigner de $\mathfrak{g}$ et $\mathfrak{h}$ est bien
une sous-alg\`{e}bre de $\mathfrak{g}_{0}$. Notons \'{e}galement que $%
\mathbb{C}\{e_{p+1},...,e_{n}\}$ est une sous-alg\`{e}bre ab\'{e}lienne de $%
\mathfrak{g}_{0}$.

\medskip

\noindent 
\begin{proposition}
Si $\g_0$ est une contraction d' In\"on\"u-Wigner de $\g$ qui laisse invariante la sous-alg\`ebre $\frak{h}$ de
 $\g$ alors
$$g_0=\frak{h} \oplus \frak{a}$$
o\`u $\frak{a}$ est un id\'eal ab\'elien de $\g_0$.
\end{proposition}

\medskip

\noindent \textbf{Remarques}

1. Si $\mathfrak{h}$ est un id\'{e}al de $\mathfrak{g}$ alors $\mathfrak{g}%
_{0}=\mathfrak{h}\oplus \mathfrak{a}$ avec%
\[
\lbrack \mathfrak{h},\mathfrak{a}]=0. 
\]

\medskip

2. Il existe des contractions d'alg\`{e}bre de Lie qui ne soient pas des
contractions d'In\"{o}n\"{u}-Wigner. Par exemple si nous consid\'{e}rons
l'alg\`{e}bre de Lie r\'{e}soluble de dimension $4$ donn\'{e}e par 
\[
\left[ e_{1},e_{2}\right] =e_{2},\ \left[ e_{3},e_{4}\right] =e_{4}. 
\]%
cette alg\`{e}bre peut se contracter sur l'alg\`{e}bre filiforme suivante : 
\[
\left[ e_{1},e_{2}\right] =e_{3},\ \ \left[ e_{1},e_{3}\right] =e_{4}. 
\]%
D'apr\`{e}s la proposition pr\'{e}c\'{e}dente, une telle contraction ne peut 
\^{e}tre une contraction d'In\"{o}n\"{u}-Wigner.

\medskip

3. Il existe une notion de contraction d'In\"{o}n\"{u}-Wigner pour les
groupes de Lie.\ Elle est subordonn\'{e}e \`{a} celle des alg\`{e}bres de
Lie. Ainsi tout groupe de Lie peut se contracter, en ce sens, sur un de ses
sous-groupes \`{a} un param\`{e}tre. Le groupe des rotations de dimension $3$
se contracte sur le groupe Euclidien \`{a} $2$ dimensions. Une contraction
du groupe de Lorentz homog\`{e}ne donne le groupe de Galil\'{e}e en laissant
invariant le sous-groupe d'invariance des coordonn\'{e}es temporelles.\ De m%
\^{e}me une contraction du groupe de Lorentz non-homog\`{e}ne donne le
groupe de Galil\'{e}e, en laissant invariant le sous groupe \ engendr\'{e}
par les rotations spatiales et les d\'{e}placements sur le temps. Tous ces
exemples sont d\'{e}crits dans le papier historique d'In\"{o}n\"{u}-Wigner 
\cite{I.N}.

\bigskip

\subsection{Les contractions de Weimar-Woods}

Saletan et Levy-Nahas g\'{e}n\'{e}ralis\`{e}rent la notion de contraction
d'In\"{o}n\"{u}-Wigner en consid\'{e}rant

\noindent a) Pour Saletan, des familles d'isomorphismes de la forme 
\[
f_{\varepsilon }=f_{1}+\varepsilon f_{2} 
\]%
avec $det(f_{2})\neq 0$. Ces isomorphismes se r\'{e}duisent \`{a} 
\[
f_{\varepsilon }=\varepsilon Id+(1-\varepsilon )g 
\]%
avec $det(g)=0$. Si $q$ est le nilindex de la aprtie nilpotente de $g$ dans
sa d\'{e}composition de Jordan, alors e partant d'une alg\`{e}%
bre de Lie donn\'{e}e $\mathfrak{g}$, on la contracte via $f_{\varepsilon }$ pour
obtenir une nouvelle alg\`{e}bre de Lie (si elle existe) $\mathfrak{g}_{1}$,
on contracte \`{a} nouveau $\mathfrak{g}_{1}$ via $f_{\varepsilon }$ et on
obtient une nouvelle alg\`{e}bre de Lie $\mathfrak{g}_{2}$ ainsi de suite.
On d\'{e}finit de la sorte une famille de contractions. Cette suite
stationne \`{a} l'ordre $q$, on obtient alors la contract\'{e}e de Saletan.
Notons que la contraction d'In\"{o}n\"{u}-Wigner correspond \`{a} $q=1$.

\noindent b) Pour Levi-Nahas, qui g\'{e}n\'{e}ralise la construction de
Saletan , la famille $f_{\varepsilon }$ \ d'isomorphismes a la forme
suivante: 
\[
f=\varepsilon f_{1}+(\varepsilon )^{2}f_{2}, 
\]%
o\`{u} $f_{1}$ et $f_{2}$ v\'{e}rifient les hypoth\`{e}ses de Saletan.

\bigskip

Les contractions de Weimar-Woods sont plus g\'{e}n\'{e}rales et permettent
une graduation dans la contraction.\ On consid\`{e}re dans ce cas des
isomorphismes du type%
\[
f(e_{i})=\epsilon ^{n_{i}}e_{i} 
\]%
o\`{u} $n_{i}\in \mathbb{Z}$ et la contraction est donn\'{e}e lorsque $%
\epsilon \rightarrow 0.$ Ce type de contractions peut \^{e}tre consid\'{e}%
rer comme des contractions g\'{e}n\'{e}ralis\'{e}es d'In\"{o}n\"{u}-Wigner
avec des exponents entiers.\ L'int\'{e}r\^{e}t est de pouvoir construire
toutes les contractions et de faire le lien avec la notion de d\'{e}%
formation que l'on pr\'{e}sente dans le paragraphe suivant.

\bigskip

\textbf{Remarque. }R.\ Hermann introduit aussi une notion de contraction,
mais dans ce cas la notion de dimension n'est plus un param\`{e}tre fixe.

\bigskip

\section{D\'{e}formations d'alg\`{e}bres de Lie}

\medskip La notion de d\'{e}formation se veut une notion duale de la notion
de contraction.\ En gros, \'{e}tant donn\'{e}e une alg\`{e}bre de Lie,
peut-on d\'{e}terminer toutes les alg\`{e}bres de Lie se contractant sur
cette alg\`{e}bre donn\'{e}e. Ceci revient \`{a} d\'{e}terminer, \'{e}tant
donn\'{e} un point de $L^{n}$ toutes les orbites ayant ce point comme point
adh\'{e}rent. Alors que pour la notion de contraction, le fait de pouvoir
utiliser le crit\`{e}re de S\'{e}gal permettait de regarder le probl\`{e}me
dans l'espace vectoriel des constantes de structure muni de la topologie m%
\'{e}trique, pour le probl\`{e}me inverse, nous sommes oblig\'{e}s de
travailler dans $L^{n}$ avec la topologie de Zariski. Nous allons donc d\'{e}%
finir une d\'{e}formation d'un point de $L^{n}$ comme un point proche au
sens de Zariski de ce point. La notion de point g\'{e}n\'{e}rique va donc
jouer un r\^{o}le pr\'{e}pond\'{e}rant pour param\'{e}trer ces points.

\subsection{Points g\'{e}n\'{e}riques dans $L^{n}$}

Soit $Spec(A(L^{n}))$ le spectre de l'anneau $A(L^{n})$ muni de sa topologie
de Zariski$\ $\ Nous savons qu'un point $P\in Spec(A(L^{n}))$ est ferm\'{e}
si et seulement si l'id\'{e}al $P$ est maximal. Dans ce cas l'adh\'{e}rence
du point $P$ est $P$ lui m\^{e}me.\ On veut g\'{e}n\'{e}raliser cette
situation.

\begin{defi}
Soit $Z$ un sous-ensemble ferm\'{e} irr\'{e}ductible de $Spec(A(L^{n})).$ Un
point $P\in Z$ est dit point g\'{e}n\'{e}rique de $Z$ si $\ Z$ est \'{e}gal 
\`{a} l'adh\'{e}rence du point $P$. Ceci revient \`{a} dire que tout ouvert
de $Z$ contient $P$.
\end{defi}

En particulier si $P~$est un id\'{e}al premier, c'est un point g\'{e}n\'{e}%
rique de l'ensemble $V(P)$ suppos\'{e} irr\'{e}ductible$.$ Nous allons nous
int\'{e}resser aux points g\'{e}n\'{e}riques de $Spec(A(L^{n}))$ qui
appartiennent \`{a} des composantes irr\'{e}ductibles de $L^{n}$ passant par
un point repr\'{e}sentant une alg\`{e}bre de Lie donn\'{e}e et poss\'{e}dant
des points g\'{e}n\'{e}riques.\ Ces points g\'{e}n\'{e}riques vont repr\'{e}%
senter les d\'{e}formations de l'alg\`{e}bre donn\'{e}e. Il existe une mani%
\`{e}re classique de construire ces points g\'{e}n\'{e}riques bas\'{e}e sur
des notions d'extension. Rappelons tout d'abord le r\'{e}sultat suivant

\begin{lemm}
Soit $\mu $ $\ =(C_{ij}^{k})$ un point de la vari\'{e}t\'{e} affine complexe 
$L^{n}$. Alors pour des extensions finies $K$ de $\mathbb{C}$, \ il existe
des s\'{e}ries enti\`{e}res $C_{ij}^{k}(t)\in K[[t]]$ telles que $%
C_{ij}^{k}(0)=C_{ij}^{k}$ et le point $(C_{ij}^{k}(t))$ est un point g\'{e}n%
\'{e}rique de $L^{n}.$
\end{lemm}

En particulier, nous pouvons consid\'{e}rer comme point g\'{e}n\'{e}rique
des points \`{a} coefficients dans l'anneau $\mathbb{C}[[t]].$

\bigskip

Une autre fa\c{c}on de d\'{e}terminer les points g\'{e}n\'{e}riques et
d'utiliser une extension non archim\'{e}dienne de $\mathbb{C}$, appel\'{e}e
l'extension de Robinson. \ Les adeptes de l'analyse non standard utilise
cette extension mais dans un contexte particulier.\ Ici nous sommes int\'{e}%
ress\'{e}s par cette extension, not\'{e}e $\mathbb{C}^{\ast }$, car elle est
munie d'une valuation et son anneau de valuation a des propri\'{e}t\'{e}s
analogues \`{a} l'anneau des s\'{e}ries formelles. De plus cette extension b%
\'{e}n\'{e}ficie d'un principe de transfert ce qui est d'un b\'{e}n\'{e}fice
appr\'{e}ciable. Soit donc $\mathbb{C}^{\ast }$ une extension non archim\'{e}%
dienne de Robinson.\ C'est un corps valu\'{e} contenant $\mathbb{C}$. Les \'el\'ements de 
$\mathbb{C}^*- \mathbb{C}$  sont dits non-standard. Il existe un principe de transfert de $\mathbb{C}$
\`a $\mathbb{C}^*$  qui peut se r\'esumer en disant que pour v\'erifier qu'une 
formule usuelle d\'ependant de param\`etres standard est vraie pour tout $x$, il suffit de la v\'erifier 
pour tout $x$ standard. On notera par $A$ l'anneau
de valuation de $\mathbb{C}^*$ qui est muni d'une structure d'alg\`{e}bre. Comme $A$ est un
anneau de valuation, c'est un anneau local (nous reviendrons sur cette d\'{e}%
finition) et poss\`{e}de donc un unique id\'{e}al maximal $\mathfrak{m}$ qui
a la propri\'{e}t\'{e} suivante :%
\[
\forall x\in \mathbb{C}^{\ast }-A,x^{-1}\in \mathfrak{m.} 
\]%
Posons $N_1=\frac{n^{3}-n^{2}}{2}.$ Pour tout $x\in \mathbb{C}^{N_1},$ soit%
\[
I_{x}=\left\{ f\in A(L^{n}),f(x)=0.\right\} 
\]%
C'est un id\'{e}al premier de $A(L^{n})$ donc un \'{e}l\'{e}ment de $%
Spec(A(L^{n})).$

\begin{defi}
Soit $P$ un \'{e}l\'{e}ment de $Spec(A(L^{n})).$ Un \'{e}l\'{e}ment $x\in 
\mathbb{C}^{\ast N_1}$ est dit g\'{e}n\'{e}rique pour $P$ si%
\[
f\in P\Longleftrightarrow f(x)=0. 
\]
\end{defi}

Ceci signifie que $P=I_{x}.$ Il est clair que l'\'{e}l\'{e}ment $x\in 
\mathbb{C}^{n}$ est g\'{e}n\'{e}rique si et seulement si $P=I_{x}$ est maximal.

\subsection{$B$-D\'{e}formations et d\'{e}formations g\'{e}n\'{e}riques
d'une alg\`{e}bre de Lie}

\medskip Soit $B$ une $\mathbb{C}$-alg\`{e}bre associative unitaire.\ C'est
en particulier un anneau unitaire. Il existe alors un homomorphisme d'anneau
unitaire%
\[
\omega :\mathbb{C\rightarrow B} 
\]%
donn\'{e} par $\omega (1)=e$ o\`{u} $e$ est l'unit\'{e} de $B$. Un
homomorphisme 
\[
\epsilon :B\rightarrow \mathbb{C} 
\]%
est appel\'{e} une augmentation s'il v\'{e}rifie%
\[
\epsilon (\omega (a))=a 
\]%
pour tout $a\in \mathbb{C}$. $\ $\ Dans ce cas $\overline{B}=Ker\epsilon $
est un id\'{e}al de $B$ et les \'{e}l\'{e}ments de l'anneau quotient $B/\overline{B}$ $%
\ $sont appel\'{e}s les ind\'{e}composables de $B$.\ Par exemple l'anneau $%
\mathbb{C}[[t]]$ est muni d'une augmentation, elle est d\'{e}finie par $%
\epsilon (\sum_{n\geq 0}a_{n}t^{n})=a_{0}.$ De m\^{e}me l'anneau de
valuation $A$ de l'extension de Robinson est muni d'une augmentation.\ En
effet, par construction, il existe un homomorphisme, appel\'{e} dans ce cas
ombre,%
\[
\alpha \in A\rightarrow ^{\circ }a\in \mathbb{C} 
\]%
ce qui se traduit, dans le langage nonstandard, en disant que tout \'{e}l%
\'{e}ment limit\'{e}, c'est-\`{a}-dire non infiniment grand, est infiniment
proche d'un unique \'{e}l\'{e}ment standard appel\'{e} son ombre. Rappelons
qu'un anneau est dit local s'il poss\`{e}de un unique id\'{e}al maximal.

\bigskip \noindent \noindent%

\begin{defi}
Soit $\g$ une alg\`ebre de Lie complexe et $B$ une $\C$-alg\`ebre munie d'une augmentation. Une $B$-d\'eformation de $\g$
est une alg\`ebre de Lie $\h$ sur $B$ telle que $\g$ soit $\C$-isomorphe \`a l'alg\`ebre de Lie complexe 
$\overline \h = \C \otimes _B \h$.
\end{defi}

\medskip

Par alg\`ebre de Lie sur l'alg\`ebre $B$, on entend un $B$-module muni d'un produit de Lie.
Notons que, comme $B$ admet une augmentation, alors $\mathbb{C}$ admet une
structure de $B$-module.\ Ainsi le produit tensoriel $\mathbb{C}\otimes _{B}%
\mathfrak{h}$ est bien d\'{e}fini.\ La multiplication externe dans $%
\overline{\mathfrak{h}}$ est donn\'{e}e par $\alpha (\alpha ^{\prime
}\otimes x)=\alpha \alpha ^{\prime }\otimes x.$ Si on suppose de plus que $%
\mathfrak{h}$ soit un $B$-module libre alors $\mathfrak{h}=B\otimes 
\mathfrak{g.}$ Notons par $\varphi _{\mathfrak{h}}$ l'isomorphisme entre $%
\overline{\mathfrak{h}}$ et $\mathfrak{g}$. \ Deux $B$-d\'{e}formations $(%
\mathfrak{h,}\varphi _{\mathfrak{h}})$ et $(\mathfrak{h}_{1}\mathfrak{,}%
\varphi _{\mathfrak{h}_{1}})$ sont dites \'{e}quivalentes s'il existe un
isomorphisme de $B$-alg\`{e}bres de Lie $\phi :\mathfrak{h\rightarrow h}_{1}$
tel que $\overline{\phi }=$ $\varphi _{\mathfrak{h}_{1}}^{-1}\circ \varphi _{%
\mathfrak{h}}$ o\`{u} $\overline{\phi }:$ $\overline{\mathfrak{h}}%
\rightarrow $ $\overline{\mathfrak{h}}_{1}$ est l'isomorphisme de $\mathbb{C}
$-alg\`{e}bre de Lie qui se d\'{e}duit de $\phi .$ Nous pouvons noter que si 
$B$ est un anneau local complet , alors il suffit de supposer que $\phi $
soit un homomorphisme \ tel que $\overline{\phi }=$ $\varphi _{\mathfrak{h}%
_{1}}^{-1}\circ \varphi _{\mathfrak{h}}$ pour que ce soit une \'{e}%
quivalence de d\'{e}formations (c'est-\`{a}-dire, $\phi $ est inversible).

\begin{propo}
Soit $\mathfrak{h}$ une $B$-d\'{e}formation de $\mathfrak{g.}$ Alors $%
\mathfrak{h}$ est de type fini sur $B$ si et seulement si $\mathfrak{g}$ est
de dimension finie (sur $\mathbb{C}$)
\end{propo}

A ce stade l\`{a}, nous avons une d\'{e}finition g\'{e}n\'{e}rale d'une d%
\'{e}formation. Bien entendu, le comportement de la d\'{e}formation d\'{e}%
pend fortement des propri\'{e}t\'{e}s de l'anneau $B$. Ici, nous allons nous
int\'{e}resser maintenant \`{a} une classe particuli\`{e}re de $B$-d\'{e}%
formations, le choix de $B$ sera dict\'{e} pour qu'une d\'{e}formation soit
un point g\'{e}n\'{e}rique des composantes alg\'{e}briques de $L^{n}$
passant par le point correspondant \`{a} l'alg\`{e}bre de Lie donn\'{e}e.

\medskip

\noindent 

\begin{defi}
Un anneau $B$ est dit un anneau de d\'eformation si $B$ est un anneau de valuation 
dont le corps r\'esiduel est isomorphe \`a $\C$. On appelle d\'eformation (valu\'ee) d'une alg\`ebre de Lie complexe, 
une 
$B$-d\'eformation o\`u $B$ est un anneau de d\'eformation. La d\'eformation sera dite g\'en\'erique si elle correspond \`a un point g\'en\'erique 
de la vari\'et\'e $L^n$.
\end{defi}

\medskip

Notons que tout anneau de d\'eformation est naturellement muni d'une structure de $\mathbb{C}$-alg\`ebre.
Soit $K$ le corps des fractions de $B$ et soit $v$ la valuation.\ Alors $B$
correspond aux \'{e}l\'{e}ments de $K$ qui ont une valuation positive et l'id%
\'{e}al maximal $\mathfrak{m}$ correspond aux \'{e}l\'{e}ments de valuation
strictement positive. Par hypoth\`{e}se le corps r\'{e}siduel $B/\mathfrak{m=%
}\mathbb{C}$. Ainsi $B$ est muni naturellement d'une augmentation.

\medskip

\textbf{Exemples.}

\begin{itemize}
\item L'anneau des s\'{e}ries formelles $\mathbb{C[[}t]]$ est un anneau de d%
\'{e}formations.\ Ici $\mathfrak{m}$ correspond aux s\'{e}ries formelles
sans terme constant. Les d\'{e}formations correspondantes sont appel\'{e}es
des d\'{e}formations formelles ou d\'{e}formations de Gerstenhaber.

\item L'anneau des \'{e}l\'{e}ments limit\'{e}s $\mathcal{L}$ dans une
extension non archim\'{e}dienne de Robinson.\ Dans ce cas l'id\'{e}al
maximal $\mathfrak{m}$ correspond aux infiniment petits. Les d\'{e}%
formations correspondantes sont appel\'{e}es des perturbations.
\end{itemize}

Ces deux exemples d\'{e}finissent donc des d\'{e}formations g\'{e}n\'{e}%
riques.\ Nous allons donc les \'{e}tudier en d\'{e}tail.

\subsection{D\'{e}formations formelles}

\noindent Rappelons que ces d\'{e}formations sont les $\mathbb{C[[}t]]-$d%
\'{e}formations. Les points correspondants dans la vari\'{e}t\'{e} $L^{n}$
sont, d'apr\`{e}s le d\'{e}but du paragraphe, des points g\'{e}n\'{e}riques.
Rappelons donc la d\'{e}finition de ces d\'{e}formations. \ Auparavant, d%
\'{e}finissons le produit $\circ $ qui \`{a} deux formes bilin\'{e}aires
altern\'{e}es fait correspondre une forme de degr\'{e} $3$:%
\[
\varphi \circ \psi (X,Y,Z)=\varphi (\psi (X,Y),Z)+\varphi (\psi
(Y,Z),X)+\varphi (\psi (Z,X),X). 
\]%
En particulier si $\mu $ est une multiplication d'alg\`{e}bre de Lie et $%
\varphi $ une application bilin\'{e}aire altern\'{e}e, alors%
\[
\delta _{\mu }\varphi =\mu \circ \varphi +\varphi \circ \mu 
\]%
et l'identit\'{e} de Jacobi se r\'{e}sume \`{a} $\mu \circ \mu =0.$ Ces
produits font partie de la classe des produits de Gerstenhaber.

\noindent%

\begin{theo}
Une d\'eformation formelle d'une alg\`ebre de Lie correspondant au point  $\mu _{0}\in L^{n}$ est donn\'ee  par
une famille $\{\varphi_i\}_{i\in \mathbb{N}}$ d'applications bilin\'eaires altern\'ees telles que $\varphi_0=\mu_0$  v\'erifiant
$$
\left\{
 \begin{array}{l}

\mu_0\circ \mu_0=0
\\
\mu_0\circ \varphi_1+\varphi_1 \circ \mu_0=\delta _{\mu_0}\varphi _1=0\\

\varphi_1\circ \varphi_1=-\mu_0\circ \varphi_2-\varphi_2 \circ \mu_0=-\delta _{\mu_0}\varphi _2\\
\vdots  \\ 
\varphi _{p}\circ \varphi _{p}+\sum_{1\leq i\leq p-1}\varphi _{i}\circ
\varphi _{2p-i}+\varphi _{2p-i}\circ \varphi _{i}=-\delta _{\mu _{0}}\varphi
_{2p} \\ 
\sum_{1\leq i\leq p}\varphi _{i}\circ \varphi _{2p+1-i}+\varphi
_{2p+1-i}\circ \varphi _{i}=-\delta _{\mu _{0}}\varphi _{2p+1} \\ 
\vdots 
\end{array}\right. .
$$
\end{theo}

\medskip \textit{D\'{e}monstration. }Montrons tout d'abord que la
multiplication $\mu $ de la $\mathbb{C[[}t]]-$d\'{e}formation de $\mu _{0}$
est enti\`{e}rement d\'{e}finie par sa restriction \`{a} $\mathbb{C}^{n}.$
La $\mathbb{C[[}t]]-$alg\`{e}bre de Lie de multiplication $\mu $ a pour
module sous-jacent le module $\mathbb{C[[}t]]\otimes \mathbb{C}^{n}.$ Chaque 
\'{e}l\'{e}ment de cette alg\`{e}bre s'\'{e}crit donc comme somme finie d'%
\'{e}l\'{e}ments du type $s(t)\otimes v,$ $s(t)\in \mathbb{C[[}t]]$ et $v\in 
$ $\mathbb{C}^{n}.$ On en d\'{e}duit $\mu (s(t)\otimes v,r(t)\otimes
w)=s(t)r(t)\mu (v,w)$. Ainsi $\mu $ est d\'{e}fini par sa restriction \`{a} $%
\mathbb{C}^{n}.$ On en d\'{e}duit, que si $v,w\in \mathbb{C}^{n},$ alors $%
\mu $ s'\'{e}crit%
\[
\mu (v,w)=\mu _{0}(v,w)+t\varphi _{1}(v,w)+\ldots +t^{p}\varphi
_{p}(v,w)+\ldots 
\]%
On \'{e}crira donc naturellement $\mu _{t}$ la d\'{e}formation de $\mu$. 
L'identit\'{e} de Jacobi pour la multiplication 
\[
\mu _{t}=\sum_{p\geq 0}t^{p}\varphi _{p} 
\]%
s'\'{e}crit 
\[
\mu _{t}\circ \mu _{t}=\mu _{0}\circ \mu _{0}+t\delta _{\mu _{0}}\varphi
_{1}+t^{2}(\varphi _{1}\circ \varphi _{1}+\delta _{\mu _{0}}\varphi
_{2})+t^{3}(\varphi _{1}\circ \varphi _{2}+\varphi _{2}\circ \varphi
_{1}+\delta _{\mu _{0}}\varphi _{3})+...=0 
\]%
ce qui est \'{e}quivalent au syst\`{e}me infini ci-dessus. On s'aper\c{c}oit
que le premier terme $\varphi _{1}$ de la d\'{e}formation $\mu _{t}$ de $\mu
_{0}$ appartient \`{a} $Z^{2}(\mu _{0},\mu _{0}).$ Ce terme est appel\'{e}
la partie infinit\'{e}simale de $\mu _{t}.$

\medskip

\noindent \noindent%

\begin{defi}
Une d\'eformation formelle de $\mu _{0}$ est appel\'ee une d\'eformation lin\'eaire si elle est de longueur 
1, c'est-\`a-dire si elle s'\'ecrit $\mu_t=\mu _{0}+t\varphi _{1}$  où $\varphi
_{1}\in Z^{2}(\mu _{0},\mu _{0}).$
\end{defi}

\medskip

Pour une telle d\'{e}formation, on a n\'{e}cessairement $\varphi
_{1}\circ \varphi _{1}=0$ soit $\varphi _{1}\in L^{n}.$

\bigskip

Consid\'{e}rons maintenant le probl\`{e}me suivant : soit $%
\varphi _{1}\in Z^{2}(\mu _{0},\mu _{0})$ avec $\mu _{0}\in L^{n}.$
Est-ce-que cette application est le premier terme d'une d\'{e}formation de $%
\mu _{0},$ autrement dit est-ce la partie infinit\'{e}simale d'une d\'{e}%
formation? Si c'est le cas, il existe une famille $\varphi _{i}\in C^{2}(\mu
_{0},\mu _{0}),$ $i\geq 2$, tel que le syst\`{e}me $(I)$ soit satisfait. Ce
probl\`{e}me d'existence est appel\'{e} le probl\`{e}me d'int\'{e}gration
formelle de $\varphi _{1}$ au point $\mu _{0}$. Comme le syst\`{e}me $(I)$
est infini, nous allons essayer de le r\'{e}soudre par induction. Pour $%
p\geq 2,$ soit $(I_{p})$ le sous-syst\`{e}me donn\'{e} par%
\[
(I_{p})\text{ }\left\{ 
\begin{array}{l}
\varphi _{1}\circ \varphi _{1}=-\delta _{\mu _{0}}\varphi _{2} \\ 
\varphi _{1}\circ \varphi _{2}+\varphi _{2}\circ \varphi _{1}=-\delta _{\mu
_{0}}\varphi _{3} \\ 
\vdots \\ 
\sum_{1\leq i\leq \lbrack p/2]}a_{i,p}(\varphi _{i}\circ \varphi
_{p-i}+\varphi _{p-i}\circ \varphi _{i})=-\delta _{\mu _{0}}\varphi _{p} \\ 
\end{array}%
\right. 
\]%
avec $a_{i,p}=1$ si $i\neq p/2$ $\ $et $a_{i,p}=1/2$ \ si $i=[p/2]$.

\medskip

\noindent \noindent%

\begin{defi}
On dit que $\varphi _{1}\in Z^{2}(\mu _{0},\mu _{0})$ est int\'egrable jusqu'\`a l'ordre $p$ s'il existe 
$\varphi _{i}\in C^{2}(\mu _{0},\mu _{0}),$ $i=2,\cdots,p$ tel que $(I_p)$ soit satisfait.
\end{defi}

\medskip

Supposons que $\varphi _{1}$ soit int\'{e}grable jusqu'\`{a} l'ordre $p$. On
montre directement que $$\sum_{1\leq i\leq \lbrack p+1/2]}a_{i,p+1}(\varphi
_{i}\circ \varphi _{p+1-i}+\varphi _{p+1-i}\circ \varphi _{i})\in Z^{3}(\mu
_{0},\mu _{0}).$$ Ainsi $\varphi _{1}$ est int\'{e}grable jusqu'\`{a} l'ordre 
$p+1$ si et seulement si cette $3$-cochaine est dans $B^{3}(\mu _{0},\mu
_{0})$. On en d\'{e}duit

\bigskip

\noindent 

\begin{proposition}
Si $H^{3}(\mu _{0},\mu _{0})=0$ alors tout $\varphi _{1}\in Z^{2}(\mu
_{0},\mu _{0})$ est la partie infinit\'esimale d'une d\'eformation de $\mu
_{0}.$
\end{proposition}

\medskip

La classe de cohomologie $[$$\sum_{1\leq i\leq \lbrack
p+1/2]}a_{i,p+1}(\varphi _{i}\circ \varphi _{p+1-i}+\varphi _{p+1-i}\circ
\varphi _{i})]$ est appel\'{e} l'obstruction d'ordre $p+1$. L'obstruction
d'ordre $2$ est donc donn\'{e}e par la classe de cohomologie de $\varphi
_{1}\circ \varphi _{1}.$ Elle peut \^{e}tre \'{e}crite en utilisant la forme
quadratique suivante:

\noindent%

\begin{defi}
La forme quadratique de Rim
 $$sq : H^2(\mu_0,\mu_0)\longrightarrow H^3(\mu_0,\mu_0)$$
est d\'efinie par
$$sq([\varphi _1])=[\varphi_1 \circ \varphi _1]$$
pour tout  $\varphi _{1}\in Z^{2}(\mu _{0},\mu _{0})$. 
\end{defi}

\medskip

\noindent En utilisant cette application, l'obstruction d'ordre $2$ s'\'{e}%
crit : $sq([\varphi _{1}])=0.$

\bigskip

\noindent \textbf{Remarque.} Dans le paragraphe suivant, nous allons \'{e}%
tudier les d\'{e}formations formelles dans le cadre des d\'{e}formations valu%
\'{e}es (ou g\'{e}n\'{e}riques). Nous verrons alors que le syst\`{e}me
infini $(I)$ est \'{e}quivalent \`{a} un syst\`{e}me fini ce qui signifie
qu'il n'existe qu'un nombre fini d'obstructions.

\bigskip

\noindent\textbf{D\'{e}formations formelles \'{e}quivalentes. }Consid\'{e}%
rons le groupe%
\[
GL_{t}(n)=GL(n,\mathbb{C)\otimes C}[[t]]=\left\{ u=Id_{n}+t\phi
_{1}+t^{2}\phi _{2}+\ldots ,\text{ }\phi _{i}\in gl(n,\mathbb{C)}\right\} 
\]%
la multiplication \'{e}tant induite par la composition des applications. L'%
\'{e}quivalence des d\'{e}formations se traduit donc ici de la mani\`{e}re
suivante : deux d\'{e}formations $\mu _{t}=\sum_{p\geq
0}t^{p}\varphi _{p}$ et $\mu _{t}^{\prime }=\sum_{p\geq
0}t^{p}\varphi _{p}^{\prime }$ de $\mu $\ sont \'{e}quivalentes si et
seulement si il existe $u\in GL_{t}(n)$ tel que%
\[
u\circ ( \sum_{p\geq 0}t^{p}\varphi _{p})=(\sum_{p\geq
0}t^{p}\varphi _{p}^{\prime })\circ (u\otimes u). 
\]

\bigskip

\noindent \textbf{Exemples}

1. Soit $\mu _{0}$ la multiplication de l'alg\`{e}bre de Lie nilpotente
nilpotente filiforme de dimension $n$ donn\'{e}e par 
\[
\mu (X_{1},X_{i})=X_{i+1} 
\]%
pour $i=2,..,n-1$. Alors toute alg\`{e}bre de Lie filiforme de dimension $n$
est isomorphe \`{a} une d\'{e}formation formelle lin\'{e}aire (D\'efinition 17).

2. Consid\'{e}rons une alg\`{e}bre de Lie frob\'{e}niusienne de dimension $%
2p $. Au paragraphe pr\'{e}c\'{e}dent, nous avons donn\'{e} la
classification de ces alg\`{e}bres de Lie \`{a} contraction pr\`{e}s. Dans 
\cite{A.C.1} on montre qu'une telle multiplication peut s'\'{e}crire, \`{a}
isomorphisme pr\`{e}s, sous la forme 
\[
\mu =\mu _{0}+t\varphi _{1} 
\]%
o\`{u} $\mu _{0}$ est une multiplication d'un mod\`{e}le frob\'{e}niusien.
Ainsi, toute multiplication d'une alg\`{e}bre de Lie frob\'{e}niusienne de
dimension $2p$ est formellement \'{e}quivalent \`{a} une d\'{e}formation lin%
\'{e}aire de la multiplication d'une alg\`{e}bre de Lie frob\'{e}niusienne
mod\`{e}le.

\bigskip

\noindent \noindent%

\begin{defi}
Une  d\'eformation formelle $\mu _{t}^{{}}$ de $\mu _{0}$ est dite triviale si elle est \'equivalente \`a
 $\mu _{0}.$
\end{defi}

\bigskip

Soit $\mu _{t}^{1}=\mu _{0}+\sum_{p=1}^{\infty }t^{p}\varphi _{p}$ et $\mu
_{t}^{2}=\mu _{0}+\sum_{t=p}^{\infty }t^{p}\psi _{p}$ deux d\'{e}formations 
\'{e}quivalentes de $\mu _{0}.$ Alors $u\circ (\sum_{p\geq
0}t^{p}\varphi _{p})=(\sum_{p\geq 0}t^{p}\psi _{p})\circ (u\otimes
u).$ Ceci implique en particulier que, apr\`{e}s avoir remplac\'{e} $u$ par $%
Id_{n}+t\phi _{1}+t^{2}\phi _{2}+\ldots \,$, 
\[
\varphi _{1}+\phi _{1}\circ \mu _{0}=\psi _{1}+\mu _{0}\circ (\phi
_{1}\otimes Id+Id\otimes \phi _{1}) 
\]%
ce qui s'\'{e}crit 
\[
\varphi _{1}-\psi _{1}\in B^{2}(\mu _{0},\mu _{0}). 
\]%
Comme\ $\varphi _{1}\in Z^{2}(\mu _{0},\mu _{0}),$ on en d\'{e}duit que $%
\psi _{1}$ appartient \`{a} la m\^{e}me classe de cohomologie que $\varphi
_{1}.\ $Notons que si $\varphi _{1}$ est nul, cette propri\'{e}t\'{e} porte
sur le premier terme non nul.

\medskip \noindent%
\begin{theo}
L'espace $H^{2}(\mu _{0},\mu _{0})$ param\'etrise les classes d'\'equivalence des d\'eformations lin\'eaires 
des d\'eformations formelles de
$\mu _0$.
\end{theo}

\bigskip

\noindent{\bf Remarque : D\'eformations et cohomologie}. Le r\'esultat pr\'ec\'edent fait un lien entre la th\'eorie des
d\'eformations et la cohomologie \`a valeur dans l'alg\`ebre. Cette cohomologie est la cohomologie de Chevalley pour les alg\`ebres de
Lie, la cohomologie d'Hochschild dans le cas associatif, etc. Mais il est faut de croire que la th\'eorie des d\'eformations permet de
d\'efinir une th\'eorie cohomologique (et pourtant ceci se lit assez souvent). La th\'eorie des d\'eformations permet 
de d\'efinir des espaces de 
degr\'e $1$ et $2$ qui correspondraient aux espaces $H^1$ et $H^2$ si la cohomologie existait. Sinon on est r\'eduit \`a d\'efinir un complexe
dont les $3$-cochaines sont triviales. Il existe des alg\`ebres ternaires pour lesquelles aucune 
cohomologie n'est d\'efinissable (voir l'article de Nicolas Goze et Elisabeth Remm sur arxiv 0803.0553 ), 
la classe des alg\`ebres de Jordan en est un autre exemple et pourtant on sait d\'eformer ces structures.

\medskip

On  d\'{e}duit du th\'eor\`eme ci-dessus :

\noindent \noindent%

\begin{coro}
If $H^{2}(\mu _{0},\mu _{0})=0$, alors toute d\'eformation formelle de $\mu _{0}$ est
triviale
\end{coro}

\medskip

\noindent \textbf{Remarque. }On peut r\'{e}duire la notion de d\'{e}%
formation formelle en consid\'{e}rant des d\'{e}formations convergentes.\
Dans ce cas une telle d\'{e}formation peut s'\'{e}crire ainsi

\begin{defi}
 Une d\'{e}formation formelle convergente $\mu _{t}$ de la
multiplication $\mu $ d'alg\`{e}bre de Lie est donn\'{e}e par 
\[
\mu _{t}=\mu _{0}+S_{1}(t)\varphi _{1}+S_{2}(t)\varphi _{2}+S_{3}(t)\varphi
_{3}+...+S_{p}(t)\varphi _{p} 
\]%
o\`{u}

1. $\{\varphi _{1},\varphi _{2},\varphi _{3},...,\varphi _{p}\}$ sont lin%
\'{e}airement ind\'ependantes dans $C^{2}(\mathbb{C}^{n},\mathbb{C}^{n})$

2. Les $S_{i}(t)$ sont des s\'{e}ries convergentes de rayon de convergence $%
r_{i}>0$

3. Les valuations $v_{i}$ de $S_{i}(t)$ satisfont $v_{i}<v_{j}$ pour $i<j.$

\end{defi}

Par exemple la d\'{e}formation formelle

\[
\widetilde{\mu _{t}}=\mu _{0}+\sum_{t=1}^{\infty }t^{i}\varphi _{1} 
\]%
s'\'{e}crit comme une d\'{e}formation convergente de longueur $1$%
\[
\mu _{t}=\mu _{0}+t\frac{1-t^{n}}{1-t}\varphi _{1} 
\]

Nous \'{e}tudierons plus g\'{e}n\'{e}ralement ces d\'{e}formations dans le
cadre des d\'{e}formations valu\'{e}es.

\subsection{Perturbations d'alg\`{e}bres de Lie}

Cette notion de d\'{e}formation est peut \^{e}tre la plus proche de ce
concept.\ En effet elle est bas\'{e}e sur une extension non archim\'{e}%
dienne de Robinson et permet d'utiliser le langage de la math\'{e}matique
non standard. Ainsi, sous ce point de vue une perturbation d'une alg\`{e}bre
de Lie sera une alg\`{e}bre de Lie de m\^eme dimension dont les
constantes de structure sont infiniment proche de l'alg\`{e}bre initiale.
Sans aucun doute, on a ici, un concept de d\'{e}formation qui traduit
parfaitement l'id\'{e}e m\^{e}me de d\'{e}formation. De plus, dans une telle
extension, une d\'{e}formation appara\^{\i}t comme un point g\'{e}n\'{e}%
rique de la composante alg\'{e}brique contenant le point associ\'{e}e \`{a}
l'alg\`{e}bre de Lie donn\'{e}e. Comme nous verrons que toutes ces notions
de d\'{e}formations g\'{e}n\'{e}riques sont conceptuellement \'{e}%
quivalentes, nous pouvons utiliser n'importe quelle parmi ces d\'{e}%
formations.\ L'int\'{e}r\^{e}t des perturbations est de pr\'{e}senter un
cadre calculatoire tr\`{e}s pratique. Enfin, dans une telle extension nous
avons un principe de transfert. Rappelons dans un premier temps les propri%
\'{e}t\'{e}s de cette extension.

\bigskip Soit $\mathbb{C}^{\ast }$ une extension de Robinson de $\mathbb{C}$%
. C'est un corps valu\'{e} non archim\'{e}dien.\ Si on note $A$ son alg\`{e}%
bre de valuation et $\mathfrak{m}$ l'id\'{e}al maximal de $A$, alors on a 
\[
x\in \mathbb{C}^{\ast }-A\Longleftrightarrow x^{-1}\in \mathfrak{m.} 
\]%
Ceci est en fait la d\'{e}finition d'une alg\`{e}bre de valuation.\ Dans ce
contexte, on appelle infiniment petits les \'{e}l\'{e}ments de $\mathfrak{m,}
$ limit\'{e}s les \'{e}lements de $A$ et infiniment grands les \'{e}l\'{e}%
ments de $\mathbb{C}^{\ast }-A.$ Les op\'{e}rations alg\'{e}briques sur ces 
\'{e}l\'{e}ments correspondent parfaitement \`{a} celles que l'on met
sur ce type d'\'{e}l\'{e}ments en analyse. Le corps r\'{e}siduel $A/%
\mathfrak{m}$ est isomorphe \`{a} $\mathbb{C}$. On note pour tout $x\in A$
sa classe par $^{\circ }x.$ On a un homomorphisme naturel d'anneaux $\mathbb{%
C}\rightarrow A$ et l'application $x\in A\rightarrow ^{\circ }x\in \mathbb{C}
$ est une augmentation. Si on appelle standard un \'{e}l\'{e}ment de $%
\mathbb{C}^{\ast }$ appartenant \`{a} $\mathbb{C}$, l'augmentation se
traduit en disant que tout \'{e}l\'{e}ment limit\'{e} est infiniment proche
d'un unique \'{e}l\'{e}ment standard et que tous les infiniment petis sont
infiniment proches de $0$. Nous pouvons \'{e}tendre toutes ces notions \`{a} 
$\mathbb{C}^{n}$ pour tout $n\in \mathbb{N}$. 

\bigskip

Soit $\mu _{0}$ un point de $L^{n}.$

\medskip

\noindent \noindent%

\begin{defi}
Une perturbation $\mu$ of $\mu _0$ est une $A$-d\'eformation de $\mu_0$ telle que
$$\mu(X,Y) - \mu _0 (X,Y) \in \frak{m}^n$$
pour tout $X,Y\in \Bbb{C}^n$.
\end{defi}

\medskip

Le r\'{e}sultat essentiel se r\'{e}sume en disant que toute perturbation
admet une d\'{e}composition canonique finie.\ Plus pr\'{e}cis\'{e}ment on a :

\noindent%

\begin{theo}
Soit  $\mu _0$ un point de $L^n$ et $\mu$ une perturbation of $\mu _0$. Il existe un entier $k$ et des \'el\'ements
$\epsilon _1,\epsilon _2, \cdots, \epsilon _k \in \frak{m}$ tels que pour tout $X,Y \in \C^n$ on ait
$$\mu(X,Y) = \mu _0(X,Y) + \epsilon_1 \phi_1(X,Y) + \epsilon_1\epsilon_2 \phi_2(X,Y)+\cdots+
\epsilon_1\epsilon_2\cdots \epsilon_k\phi_k(X,Y)$$
o\`u les  $\phi_i$ sont des applications bilin\'eaires altern\'ees sur $\C^n$ \`a valeurs dans $\C^n$
lin\'eairement ind\'ependantes.
\end{theo}

\bigskip

L'int\'{e}r\^{e}t d'une telle d\'{e}composition est appr\'{e}ciable.\ En
effet, alors que dans le cadre des d\'{e}formations formelles l'identit\'{e}
de Jacobi relative \`{a} la d\'{e}formation se traduit sur un syst\`{e}me
infini portant sur les applications $\varphi _{i}$, la m\^{e}me identit\'{e}
est \'{e}quivalente dans le cadre des pertubations \`{a} un syst\`{e}me fini
formellent r\'{e}solvable. Comme nous verrons qu'il y a \'{e}quivalence entre
les deux notions de d\'{e}formations, on en d\'{e}duit

\bigskip

\textit{Le syst\`{e}me infini (I) (th\'{e}or\`{e}me 4) de Gerstenhaber est 
\'{e}quivalent \`{a} un nombre fini de syst\`{e}mes finis.}

\bigskip

En particulier, le premier terme d'une perturbation est toujours dans
l'espace $Z^{2}(\mu _{0},\mu _{0}).$ Nous verrons tout cela dans le
paragraphe qui suit. Une \'{e}tude sp\'{e}cifique des perturbations est
faite dans \cite{Go}

\medskip

\medskip

\subsection{D\'{e}formations valu\'{e}es ou d\'{e}formations g\'{e}n\'{e}%
riques}

\subsubsection{Une d\'{e}composition canonique dans $\mathfrak{m}^{k}$}

Rappelons rapidement ce qu'est un anneau de valuation.\ Soit $\mathbb{F}$ un
corps commutatif et $A$ un sous anneau de $\mathbb{F}$. On dit que $A$ est
un anneau de valuation de $\mathbb{F}$ si $A$ est un anneau local int\`{e}%
gre tel que: 
\[
\mbox{ Si}\ x\in \mathbb{F}-A,\quad \mbox{alors}\quad x^{-1}\in \mathfrak{m}%
. 
\]%
o\`{u} $\mathfrak{m}$ est l'id\'{e}al maximal de $A$. Un anneau $A$ est appel%
\'{e} anneau de valuation si c'est un anneau de valuation de son corps de
fractions. Par exemple si $\mathbb{K}$ est un corps commutatif de caract\'{e}%
ristique $0,$ alors l'anneau des s\'{e}ries formelles $\mathbb{K}[[t]]$ est
un anneau de valuation alors que l'anneau $\mathbb{K}[[t_{1},t_{2}]]$ de
deux ou plus ind\'{e}termin\'{e}es ne l'est pas.

Soit $\mathfrak{m}^{2}$ le produit cart\'{e}sien $\mathfrak{m}\times 
\mathfrak{m}$ . Soit $(a_{1},a_{2})\in \mathfrak{m}^{2}$ avec $a_{i}\neq 0$
pour $i=1,2$.

\smallskip

\noindent i) Supposons $a_{1}.a_{2}^{-1}\in A$. Soit $\alpha =\pi
(a_{1}.a_{2}^{-1})$ o\`{u} $\pi $ est la projection canonique dans $A/$ $%
\mathfrak{m}$. Rappelons que nous avons un morphisme naturel\ $\omega :%
\mathbb{K}\rightarrow A$ qui permet d'identifier $\alpha $ avec $s(\alpha )$
dans $A$. Alors 
\[
a_{1}.a_{2}^{-1}=\alpha +a_{3} 
\]%
avec $a_{3}\in \mathfrak{m}$. Si $a_{3}\neq 0$, 
\[
(a_{1},a_{2})=(a_{2}(\alpha +a_{3}),a_{2})=a_{2}(\alpha ,1)+a_{2}a_{3}(0,1). 
\]%
Si $\alpha \neq 0$ on peut \'{e}crire aussi 
\[
(a_{1},a_{2})=aV_{1}+abV_{2} 
\]%
avec $a,b\in \mathfrak{m}$ et $V_{1},V_{2}$ sont lin\'{e}airement ind\'{e}%
pendants dans $\mathbb{K}^{2}$. Si $\alpha =0$, alors $a_{1}.a_{2}^{-1}\in 
\mathfrak{m}$ et $a_{1}=a_{2}a_{3}$. On a 
\[
(a_{1},a_{2})=(a_{2}a_{3},a_{2})=ab(1,0)+a(0,1). 
\]%
Ainsi dans ce cas, $V_{1}=(0,1)$ et $V_{2}=(1,0)$. Si $a_{3}=0$, alors $%
a_{1}a_{2}^{-1}=\alpha $ et $(a_{1},a_{2})=a_{2}(\alpha ,1)=aV_{1}.$Ceci
correspond \`{a} la d\'{e}composition pr\'{e}c\'{e}dent mais pour $b=0$.

\smallskip

\noindent ii) Si $a_{1}.a_{2}^{-1}\in \mathbb{F}-A$, alors $%
a_{2}.a_{1}^{-1}\in \mathfrak{m}$. Dans ce cas on pose $%
a_{2}.a_{1}^{-1}=a_{3}$ et on obtient 
\[
(a_{1},a_{2})=(a_{1},a_{1}.a_{3})=a_{1}(1,a_{3})=a_{1}(1,0)+a_{1}a_{3}(0,1) 
\]%
avec $a_{3}\in \mathfrak{m}$. On obtient dans ce cas la d\'{e}composition : 
\[
(a_{1},a_{2})=aV_{1}+abV_{2} 
\]%
avec $a,b\in \mathfrak{m}$ et $V_{1},V_{2}$ lin\'{e}airement ind\'{e}%
pendants dans $\mathbb{K}^{2}$. \ Cette d\'{e}composition se g\'{e}n\'{e}%
ralise facilement pour un point de $\mathfrak{m}^{k}$. Ceci s'\'{e}crit :

\medskip

\noindent 

\begin{theo}
Pour tout point $(a_1,a_2,...,a_k) \in \frak{m}^k$ il existe $h \ (h \leq k$) et $h$-vecteurs lin\'eairement ind\'ependants 
$V_1,V_2,..,V_h$
de l'espace  $\mathbb{K} ^k$  et
$b_1,b_2,..,b_h \in \frak{m}$ tels que
$$(a_1,a_2,...,a_k)=b_1V_1+b_1b_2V_2+...+b_1b_2...b_hV_h.$$
\end{theo}

\medskip

Le param\`{e}tre $h$ qui appara\^{\i}t dans cette d\'{e}composition est appel%
\'{e}e la longueur de la d\'{e}composition. Il peut \^{e}tre inf\'{e}rieur 
\`{a} $k$. Il correspond \`{a} la dimension du plus petit $\mathbb{K}$%
-espace vectoriel $V$ tel que $(a_{1},a_{2},...,a_{k})\in V\otimes \mathfrak{%
m}$. Si les coordonn\'{e}es $a_{i}$ du vecteur $(a_{1},a_{2},...,a_{k})$
sont dans $A$ et pas n\'{e}cessairement dans l'id\'{e}al maximal, on \'{e}%
crit alors $a_{i}=\alpha _{i}+a_{i}^{\prime }$ avec $\alpha _{i}\in \mathbb{K%
}$ et $a_{i}^{\prime }\in \mathfrak{m}$, et on d\'{e}compose 
\[
(a_{1},a_{2},...,a_{k})=(\alpha _{1},\alpha _{2},...,\alpha
_{k})+(a_{1}^{\prime },a_{2}^{\prime },...,a_{k}^{\prime }). 
\]%
On applique alors le th\'{e}or\`{e}me ci-dessus au vecteur $(a_{1}^{\prime
},a_{2}^{\prime },...,a_{k}^{\prime })$.

Cette d\'{e}composition est unique au sens suivant :

\begin{theo}
Soit $b_{1}V_{1}+b_{1}b_{2}V_{2}+...+b_{1}b_{2}...b_{h}V_{h}$ et $%
c_{1}W_{1}+c_{1}c_{2}W_{2}+...+c_{1}c_{2}...c_{s}W_{s}$ deux d\'{e}%
compositions du vecteur $(a_{1},a_{2},...,a_{k})$. Alors

\noindent i. \thinspace\ $h=s$,

\noindent ii. Le drapeau engendr\'{e} par la famille libre $%
(V_{1},V_{2},..,V_{h})$ est \'{e}gal au drapeau engendr\'{e} par la famille
libre $(W_{1},W_{2},...,W_{h})$ c'est-\`{a}-dire $\forall i\in 1,..,h$ 
\[
\{V_{1},...,V_{i}\}=\{W_{1},...,W_{i}\} 
\]%
o\`{u} $\{U_{i}\}$ d\'{e}signe l'espace engendr\'{e} par les vecteurs $U_{i}$%
.
\end{theo}

\bigskip Pour la d\'{e}monstration, on pourra se \ r\'{e}ferrer \`{a} \cite%
{G.R1}. Nous allons appliquer cette d\'{e}composition \`{a} une d\'{e}formation valu%
\'{e}e.

\subsubsection{D\'{e}composition d'une d\'{e}formation valu\'{e}e}

Ceci \'{e}tant, soit $\mathfrak{g}$ une $\mathbb{K}$-alg\`{e}bre de Lie et $A$
une $\mathbb{K}$-alg\`{e}bre commutative unitaire de valuation. Supposons
que son corps r\'{e}siduel $A/\mathfrak{m}$ soit isomorphe \`{a} $\mathbb{K}$
Dans ce cas on a une augmentation naturelle.\ Soit $\mathfrak{g}_{A}$= $%
\mathfrak{g}\otimes A$ une $A$-d\'{e}formation de $\mathfrak{g}$. Si $dim_{%
\mathbb{K}}(\mathfrak{g})$ est fini, alors 
\[
dim_{A}(\mathfrak{g}_{A})=dim_{\mathbb{K}}(\mathfrak{g}). 
\]%
Comme $A$ est aussi une $\mathbb{K}$-alg\`{e}bre, on identifiera $\mathfrak{g%
}$ au sous espace $\mathfrak{g}\otimes e$ de $\mathfrak{g}_{A}$. \ En
particulier on a si $\mu _{\mathfrak{g}_{A}}$ d\'{e}signe la multiplication
de $\mathfrak{g}_{A}$ et $\mu _{\mathfrak{g}}$ celle de $\mathfrak{g,}$
alors $\mu _{\mathfrak{g}_{A}}(X,Y)-\mu _{\mathfrak{g}}(X,Y)$ appartient au
quasi-module $\mathfrak{g}\otimes \mathfrak{m}$ pour tout $X,Y\in \mathbb{K}%
^{n}.$ Les d\'{e}formations formelles de Gerstenhaber et les
perturbations sont des d\'{e}formations valu\'{e}es. Supposons que $%
\mathfrak{g}$ soit de dimension finie et soit $\{X_{1},...,X_{n}\}$ une base
de $\mathfrak{g}$. On a alors 
\[
\mu _{\mathfrak{g}_{A}}(X_{i},X_{j})-\mu _{\mathfrak{g}}(X_{i},X_{j})=%
\sum_{k}C_{ij}^{k}X_{k} 
\]%
avec $C_{ij}^{k}\in \mathfrak{m}$. Cette diff\'{e}rence appara\^{\i}t donc
comme un vecteur de $\mathfrak{m}^{n^{2}(n-1)/2}$ ayant pour composantes les 
$C_{ij}^{k}$.\ La d\'{e}composition canonique s'\'{e}crit alors 
\[
\begin{array}{lll}
\mu _{\mathfrak{g}_{A}}(X_{i},X_{j})-\mu _{\mathfrak{g}}(X_{i},X_{j}) & = & 
\epsilon _{1}\phi _{1}(X_{i},X_{j})+\epsilon _{1}\epsilon _{2}\phi
_{2}(X_{i},X_{j}) \\ 
&  & +...+\epsilon _{1}\epsilon _{2}...\epsilon _{k}\phi _{k}(X_{i},X_{j})%
\end{array}%
\]%
avec $\epsilon _{s}\in \mathfrak{m}$ et $\phi _{1},...,\phi _{l}$ :$%
\mathfrak{g}\otimes \mathfrak{g}\rightarrow \mathfrak{g}$ lin\'{e}airement
ind\'{e}pendantes.\ Cette d\'ecomposition est en
particulier valable si $A$=$\mathbb{C}[[t]].$ Nous pouvons donc r\'{e}duire
le syst\`{e}me infini d'int\'{e}grabilit\'{e} de Gerstenhaber \`{a} un syst%
\`{e}me fini. Rappelons tout d'abord que le complexe de Chevalley-Eilenberg
est un complexe diff\'{e}rentiel gradu\'{e}.\ Le crochet gradu\'{e} est d%
\'{e}fini \`{a} partir des produits $\circ $ donn\'{e}s par 
\[
(g_{q}\circ f_{p})(X_{1},...,X_{p+q})=\sum (-1)^{\epsilon (\sigma
)}g_{q}(f_{p}(X_{\sigma (1)},...,X_{\sigma (p)}),X_{\sigma
(p+1)},...,X_{\sigma (q)}) 
\]%
o\`{u} $\sigma $ est une permutation de ${1,...,p+q}$ telle que $\sigma
(1)<...<\sigma (p)$ et $\sigma (p+1)<...<\sigma (p+q)$ (c'est un $(p,q)$%
-schuffle) et $g_{q}\in \mathcal{C}^{q}(\mathfrak{g},\mathfrak{g})$ et $%
f_{p}\in \mathcal{C}^{p}(\mathfrak{g},\mathfrak{g})$. Le crochet gradu\'{e}
est alors donn\'{e} par%
\[
\lbrack f,g]=f\circ g-(-1)^{p-1}g\circ f. 
\]
La condition de Jacobi relative \`{a} $\mu _{\mathfrak{g}_{A}}$ se r\'{e}%
sume \`{a} $\mu _{\mathfrak{g}_{A}^{\prime }}\circ \mu _{\mathfrak{g}%
_{A}^{\prime }}=0$. Ceci donne 
\[
(\mu _{\mathfrak{g}}+\sum_{i\in I}\epsilon _{1}\epsilon _{2}...\epsilon
_{i}\phi _{i})\circ (\mu _{\mathfrak{g}}+\sum_{i\in I}\epsilon _{1}\epsilon
_{2}...\epsilon _{i}\phi _{i})=0.\quad \quad (1) 
\]%
Comme $\mu _{\mathfrak{g}}\circ \mu _{\mathfrak{g}}=0$, cette \'{e}quation
se r\'{e}duit \`{a} : 
\[
\epsilon _{1}(\mu _{\mathfrak{g}}\circ \phi _{1}+\phi _{1}\circ \mu _{%
\mathfrak{g}})+\epsilon _{1}U=0 
\]%
o\`{u} $U$ $\in $ $\mathcal{C}^{3}(\mathfrak{g},\mathfrak{g})\otimes 
\mathfrak{m}$. Simplifions par $\epsilon _{1}$ qui est suppos\'{e} non nul,
sinon la d\'{e}formation est triviale : 
\[
(\mu _{\mathfrak{g}}\circ \phi _{1}+\phi _{1}\circ \mu _{\mathfrak{g}%
})(X,Y,Z)+U(X,Y,Z)=0 
\]%
pour tout $X,Y,Z\in \mathfrak{g}$. Comme $U(X,Y,Z)$ $\in $ $\mathfrak{g}%
\otimes \mathfrak{m}$ et que le premier terme est dans $\mathfrak{g}$, on en
d\'{e}duit 
\[
(\mu _{\mathfrak{g}}\circ \phi _{1}+\phi _{1}\circ \mu _{\mathfrak{g}%
})(X,Y,Z)=0. 
\]%
Or $\mu _{\mathfrak{g}}\circ \phi _{1}+\phi _{1}\circ \mu _{\mathfrak{g}}$
n'est rien d'autre que $\delta _{\mu }\phi _{1}$ o\`{u} $\delta _{\mu }$ est
l'op\'{e}rateur cobord de la cohomologie de Chevalley de l'alg\`{e}bre de
Lie $\mathfrak{g}$. On retrouve donc, dans le cadre des d\'{e}formations valu%
\'{e}es le r\'{e}sultat classique de Gerstenhaber :%
\[
\delta _{\mu }\phi _{1}=0. 
\]%
Regardons maintenant les \'{e}quations donn\'{e}es par $U=0$. Elles vont s'%
\'{e}crire \`{a} l'aide du crochet gradu\'{e} sur l'espace des cochaines
rappel\'{e} ci-dessus. En particulier, on a si $\phi _{i},\phi _{j}\in 
\mathcal{C}^{2}(\mathfrak{g},\mathfrak{g}):$ 
\[
\lbrack \phi _{i},\phi _{j}]=\phi _{i}\circ \phi _{j}+\phi _{j}\circ \phi
_{i} 
\]%
et $[\phi _{i},\phi _{j}]\in \mathcal{C}^{3}(\mathfrak{g},\mathfrak{g})$.

\medskip

\noindent 

\begin{theo} Soit
$$\mu _{\frak{g}_A}=\mu _{\frak{g}}+\sum _{i \in I} \epsilon _1 \epsilon _2 ...\epsilon _i \phi_i$$
une d\'eformation valu\'ee de longueur $k$. Alors les  3-cochaines $[\phi _i,\phi _j]$ et $[\mu, \phi_i]$, 
$1 \leq i,j \leq k-1$, engendrent un sous espace vectoriel 
 $V$ de  $\mathcal{C}^3(\frak{g},\frak{g})$
de dimension inf\'erieure ou \'egale \`a $k(k-1)/2$ et $\mu _{\frak{g}_A} \circ \mu _{\frak{g}_A}=0$ est \'equivalent \`a
$$
\left\{
\begin{array}{l}
\delta \phi _1 = 0 \\
\delta \phi _2 = a_{11}^2 [\phi _1,\phi _1] \\
\delta \phi _3 = a_{12}^3 [\phi _1,\phi _2]+a_{22}^3[\phi _1,\phi _1]\\
... \\
\delta \phi _k = \sum _{1 \leq i \leq j \leq k-1} a_{ij}^{k} [\phi _i,\phi _j] \\
\lbrack \phi _1,\phi _k]=\sum _{1 \leq i \leq j \leq k-1} b_{ij}^{1} [\phi _i,\phi _j] \\
.... \\
\lbrack \phi _{k-1},\phi _k]=\sum _{1 \leq i \leq j \leq k-1} b_{ij}^{k-1} [\phi _i,\phi _j]

\end{array}
\right.
$$
\end{theo}

\medskip 

\noindent{\it D\'{e}monstration.} Soit $V$ le sous-espace de $\mathcal{C}%
^{3}(\mathfrak{g},\mathfrak{g})$ engendr\'{e} par les applications $[\phi
_{i},\phi _{j}]$ et $[\mu ,\phi _{i}]$. Si $\omega $ est une forme lin\'{e}%
aire sur $V$ dont le noyau contient les vecteurs $[\phi _{i},\phi _{j}]$
pour $1\leq i,j\leq (k-1)$, alors l'\'{e}quation (1) donne: 
\[
\epsilon _{1}\epsilon _{2}...\epsilon _{k}\omega ([\phi _{1},\phi
_{k}])+\epsilon _{1}\epsilon _{2}^{2}...\epsilon _{k}\omega ([\phi _{2},\phi
_{k}])+...+\epsilon _{1}\epsilon _{2}^{2}...\epsilon _{k}^{2}\omega ([\phi
_{k},\phi _{k}])+\epsilon _{2}\omega ([\mu ,\phi _{2}]) 
\]%
\[
+\epsilon _{2}\epsilon _{3}\omega ([\mu ,\phi _{3}])...+\epsilon
_{2}\epsilon _{3}...\epsilon _{k}\omega ([\mu ,\phi _{k}])=0. 
\]%
Comme chacun des coefficients est dans l'id\'{e}al $\mathfrak{m}$, on a n%
\'{e}cessairement 
\[
\omega ([\phi _{1},\phi _{k}])=...=\omega ([\phi _{k},\phi _{k}])=\omega
([\mu ,\phi _{2}])=...=\omega ([\mu ,\phi _{k}])=0 
\]%
et ceci pour toute forme lin\'{e}aire dont le noyau contient $V$. Le syst%
\`{e}me du th\'{e}or\`{e}me correspond aux relations de d\'{e}pendances dans 
$V$ et au fait que 
\[
\mathfrak{m}\supset \mathfrak{m}^{(2)}\supset ...\supset \mathfrak{m}%
^{(p)}... 
\]%
o\`{u} $\mathfrak{m}^{(p)}$est l'id\'{e}al engendr\'{e} par les produits${\
a_{1}a_{2}...a_{p},\ a_{i}\in \mathfrak{m}}$ de longueur $p$.

\medskip

\medskip \textbf{Cas particulier : dim V=}$k(k-1)/2$

Supposons que la dimension de $V$ soit maximum et \'{e}gale \`{a} $(k-1)/2$.
Nous allons voir que dans ce cas la d\'{e}formation est isomorphe \`{a} une $%
\mathbb{C[}t]-$d\'{e}formation.

\medskip \noindent%

\begin{proposition}
Soit $\mu_{\frak{g}_A}$ une $A$-d\'eformation de $\mu_{\frak{g}}$ de longueur $k$ telle que $dim V=k(k-1)/2.$ Elle est
alors \'equivalente \`a une d\'eformation polynomiale v\'erifiant
$$\mu_t(X,Y)=\mu_{\frak{g}}(X,Y) + \sum_{i=1,..,k} t^i \phi _i (X,Y)$$
pour tout $X,Y \in \g.$
\end{proposition}

\medskip

\noindent \textit{D\'{e}monstration. } Consid\'{e}rons l'\'{e}quation 
\[
\mu _{\mathfrak{g}_{A}}\circ \mu _{\mathfrak{g}_{A}}=0. 
\]%
Comme dim$V$=$k(k-1)/2$, il existe des polyn\^{o}mes $P_{i}(X)\in \mathbb{K}%
[X]$ de degr\'{e} $i$ tels que 
\[
\epsilon _{i}=a_{i}\epsilon _{k}\frac{P_{k-i}(\epsilon _{k})}{%
P_{k-i+1}(\epsilon _{k})} 
\]%
avec $a_{i}\in \mathbb{K}$. On a alors 
\[
\mu _{\mathfrak{g}_{A}^{\prime }}=\mu _{\mathfrak{g}_{A}}+%
\sum_{i=1,...,k}a_{1}a_{2}...a_{i}(\epsilon _{k})^{i}\frac{P_{k-i}(\epsilon
_{k})}{P_{k}(\epsilon _{k})}\phi _{i}. 
\]%
Ainsi 
\[
P_{k}(\epsilon _{k})\mu _{\mathfrak{g}_{A}^{\prime }}=P_{k}(\epsilon
_{k})\mu _{\mathfrak{g}_{A}}+\sum_{i=1,...,k}a_{1}a_{2}...a_{i}(\epsilon
_{k})^{i}P_{k-i}(\epsilon _{k})\phi _{i}. 
\]%
Le r\'{e}sultat se d\'{e}duit en \'{e}crivant cette expression suivant les
puissances croissantes.$\quad \Box $

\medskip

\noindent Notons que pour une telle d\'{e}formation, on a 
\[
\left\{ 
\begin{array}{l}
\delta \varphi _{2}+[\varphi _{1},\varphi _{1}]=0 \\ 
\delta \varphi _{3}+[\varphi _{1},\varphi _{2}]=0 \\ 
... \\ 
\delta \varphi _{k}+\sum_{i+j=k}[\varphi _{i},\varphi _{j}]=0 \\ 
\sum_{i+j=k+s}[\varphi _{i},\varphi _{j}]=0. \\ 
\end{array}%
\right. 
\]

\noindent

\section{\protect\medskip Alg\`{e}bres de Lie rigides}

\bigskip La notion de rigidit\'{e} est une notion topologique.\ Nous
pourrons la relier naturellement \`{a} celle de d\'{e}formation lorsque ces d%
\'{e}formations seront g\'{e}n\'{e}riques.\ Soit $\mathfrak{g}$ une alg\`{e}%
bre de Lie complexe de dimension $\ n$.\ On la consid\`{e}re comme un
point $\mu $ de la vari\'{e}t\'{e} alg\'{e}brique $L^{n}$ munie de sa
topologie de Zariski.

\medskip \noindent%

\begin{defi}
L'alg\`ebre de Lie $\g$ est rigide si l'orbite $\mathcal{O}(\mu)$ est Zariski ouverte dans $L^n$.
\end{defi}

\bigskip

Dans ce cas, l'adh\'{e}rence de l'orbite $\overline{\mathcal{O}(\mu )}$ est
une composante alg\'{e}brique connexe de $L^{n}.$ On en d\'{e}duit imm\'{e}%
diatement, comme toute vari\'{e}t\'{e} alg\'{e}brique complexe est r\'{e}%
union d'un nombre fini de composantes alg\'{e}briques connexes, qu'il
n'existe qu'un nombre fini de classe d'isomorphie d'alg\`{e}bres de Lie
rigides de dimension $n.$ Le th\'{e}or\`{e}me suivant, que nous ne d\'{e}%
montrons pas ici, permet de mettre en \'{e}vidence certaines alg\`{e}bres de
Lie rigides.

\noindent%

\begin{theo}Th\'eor\`eme de Nijenhuis-Richardson. 
Soit $\g=(\C^n,\mu)$ une alg\`ebre de Lie de dimension  $n$. Si le deuxi\`eme groupe $H^2(\g,\g)$ de la cohomologie de Chevalley est nul, alors $\g$ 
est  rigide. 
\end{theo}

Ainsi, toute alg\`{e}bre de Lie semi-simple complexe est rigide. Mais la r%
\'{e}ciproque du th\'{e}or\`{e}me de Nijenhuis-Richardson est fausse.\ Consid%
\'{e}rons en effet l'alg\`{e}bre de Lie de dimension $11$ d\'{e}finie dans
la base $\{X,X_{0},X_{1},\ldots ,X_{9}\}$ par,%
\[
\left\{ 
\begin{array}{l}
\mu (X,X_{i})=iX_{i},\quad i=0,\ldots ,9 \\ 
\mu (X_{0},X_{i})=X_{i},\quad i=4,5,\ldots ,9 \\ 
\mu (X_{1},X_{i})=X_{i+1},\quad i=2,4,5,6,7,8 \\ 
\mu (X_{2},X_{i})=X_{i+2},\quad i=4,5,6,7.%
\end{array}%
\right. 
\]%
On montre, soit par un calcul direct utilisant la suite d'Hochschild-Serre,
soit en utilisant un calcul sur ordinateur (un programme est donn\'{e} dans
le livre \cite{G.K}, bas\'{e} sur Mathematica) que la dimension de $%
H^{2}(\mu ,\mu )$ est \'{e}gale \`{a} $1$.\ Quant \`{a} la rigidit\'{e},
elle est montr\'{e}e (voir toujours \cite{G.K}) en utilisant le r\'{e}sultat
suivant:

\noindent 

\begin{theo}
Soit $\g=(\C^n,\mu)$ une alg\`ebre de Lie complexe de dimension $n$. Alors $\g$ est rigide
si et seulement si toute perturbation lui est isomorphe.
\end{theo}

\medskip \textit{D\'{e}monstration. }En effet, toute perturbation de $%
\mathfrak{g=(}\mathbb{C}^{n},\mu )$ est un point g\'{e}n\'{e}rique de la
composante alg\'{e}brique passant par le point $\mu .$ R\'eciproquement, si toutes les perturbations de $\mathfrak{g}$
sont isomorphes \`a $\mathfrak{g}$, son orbite est ouverte et $\mathfrak{g}$ est rigide.
D'o\`{u} le r\'{e}%
sultat.

\bigskip

Le r\'{e}sultat pr\'{e}c\'{e}dent se g\'{e}n\'{e}ralise naturellement aux d%
\'{e}formations g\'{e}n\'{e}riques.\ Ainsi

\bigskip

\noindent%

\begin{theo}
Soit $\g=(\C^n,\mu)$ une alg\`ebre de Lie complexe de dimension $n$. Alors $\g$ est rigide
si et seulement si toute d\'eformation valu\'ee lui est isomorphe.
\end{theo}

\noindent

\medskip Remarque : L'existence d'alg\`{e}bres de Lie rigides de dimension $%
n $ dont le $H^{2}(\mathfrak{g,g)}$ est non nul implique que le sch\'{e}ma $%
\mathcal{L}^{n}$ associ\'{e} \`{a} la vari\'{e}t\'{e} $L^{n}$ n'est pas r%
\'{e}duit.

\bigskip

Nous allons nous int\'{e}resser \`{a} pr\'{e}sent \`{a} la classification
des alg\`{e}bres de Lie rigides.\ Cette classification est loin d'\^{e}tre
achev\'{e}e. En particulier on ne sait rien sur la rigidit\'{e} \'{e}%
ventuelle d'alg\`{e}bres de Lie nilpotentes.\ On peut raisonnablement
conjecturer le r\'{e}sultat suivant

\bigskip

\textbf{Conjecture :}\textit{\ Il n'existe pas d'alg\`{e}bres de Lie
nilpotentes rigides.}

\bigskip

Afin de dresser une \'{e}ventuelle classification des alg\`{e}bres de Lie
rigides, nous devons commencer par comprendre la structure de ces alg\`{e}%
bres.\ Pour cela nous devons revenir sur l'\'{e}tude des alg\`{e}bres de Lie
alg\'{e}briques.

\medskip

\noindent \noindent%

\begin{defi}
Une alg\`ebre de Lie r\'eelle ou complexe  $\g$ est appel\'ee alg\'ebrique si elle est isomorphe 
\`a l'alg\`ebre de Lie d'un groupe de Lie alg\'ebrique.
\end{defi}

\bigskip \medskip

Un groupe de Lie lin\'{e}aire alg\'{e}brique complexe est un sous-groupe de
Lie d'un groupe $GL(p,\mathbb{C})$ d\'{e}fini comme l'ensemble des z\'{e}ros
d'un syst\`{e}me fini de relations polynomiales.\ C'est donc une vari\'{e}t%
\'{e} alg\'{e}brique mais avec la propri\'{e}t\'{e} surprenante que dans ce
cas il n'existe aucun point singulier.\ C'est donc bien une vari\'{e}t\'{e}
diff\'{e}rentielle.\medskip

\noindent \textbf{Exemples.}

1. Toute alg\`{e}bre de Lie simple complexe est alg\'{e}brique.

2. Toute alg\`{e}bre de Lie nilpotente est alg\'{e}brique.\ Plus g\'{e}n\'{e}%
ralement, toute alg\`{e}bre de Lie dont le radical est nilpotent est alg\'{e}%
brique.

3. Toute alg\`{e}bre de Lie complexe de dimension $n$ dont l'alg\`{e}bre de
Lie des d\'{e}rivations est aussi de dimension $n$ est alg\'{e}brique.

4. Toute alg\`{e}bre de Lie v\'{e}rifiant $\mathcal{D}^{1}(\mathfrak{g})=%
\mathfrak{g}$ est alg\'{e}brique.

\medskip

La famille des alg\`{e}bres de Lie alg\'{e}brique est donc vaste et nous
n'avons pas de classification pr\'{e}cise de cette classe. Par contre, le r%
\'{e}sultat suivant donne la structure pr\'{e}cise de ces alg\`{e}bres:

\medskip

\noindent \noindent%

\begin{proposition}
Les propositions suivantes sont \'equivalentes:

1. $\g$ est alg\'ebrique.

2. $ad(\g)=\{ad_\mu X, X \in \g\}$ est alg\'ebrique.

3. $\g=\s \oplus \n \oplus \m$ o\`u $\s$ est une sous-alg\`ebre semi-simple de Levi, $\n$ 
le nilradical c'est-\`a-dire le plus grand id\'eal nilpotent et $t$ un tore de Malcev tel que $ad_\mu \m$ soit alg\'ebrique.
\end{proposition}

\medskip

Cette proposition fait appara\^{\i}tre le tore externe de Malcev. Par d\'{e}%
finition, un tore externe de Malcev $\mathfrak{t}$ de $\mathfrak{g}$ est une
sous-alg\`{e}bre de Lie ab\'{e}lienne telle que tous les endomorphismes
$adX$, $X\in \mathfrak{t}$ soient semi-simple (simultan\'{e}ment
diagonalisables). Tous les tores de Malcev maximaux (pour l'inclusion) sont
conjugu\'{e}s.\ Leur dimension commune est appel\'{e}e le rang de $\mathfrak{%
g.}$\medskip

\bigskip

Revenons aux alg\`{e}bres de Lie rigides.\ Soit $\mathfrak{g}$ une alg\`{e}%
bre de Lie rigide dans $L^{n}.$ Son orbite est ouverte.\ On en d\'{e}duit
que tout d\'{e}formation g\'{e}n\'{e}rique lui est isomorphe.\ Si $\mathfrak{%
g}^{\prime }$ est une telle d\'{e}formation (par exemple une perturbation),
alors c'est un point g\'{e}n\'{e}rique de la composante d\'{e}finie par
l'orbite de $\mathfrak{g}$ (rappelons que l'adh\'{e}rence de l'orbite d'une
alg\`{e}bre rigide est une composante alg\'{e}brique de la vari\'{e}t\'{e} $%
L^{n}).$ Or cette composante est d\'{e}finie par l'action du groupe alg\'{e}%
brique $GL(n,\mathbb{C)}$.\ C'est l'adh\'{e}rence de l'image de $\mathfrak{g}%
^{\prime }$ par cette action. Comme cette action est d\'{e}finie par des 
\'{e}quations polynomiales et que la composante est \'{e}galement donn\'{e}e
par un syst\`{e}me d'\'{e}quations polynomiales, on en d\'{e}duit imm\'{e}%
diatement que le point g\'{e}n\'{e}rique est alg\'{e}brique.\ Comme il est
isomorphe \`{a} $\mathfrak{g,}$ cette alg\`{e}bre est aussi alg\'{e}brique.

\noindent \noindent%

\begin{proposition}
Toute alg\`ebre de Lie complexe rigide dans $L^n$ est alg\'ebrique.
\end{proposition}

\medskip

\bigskip

Supposons que $\mathfrak{g}$ soit r\'{e}soluble.\ Si elle est rigide, elle
se d\'{e}compose donc sous la forme

\[
\mathfrak{g=t\oplus n} 
\]%
o\`{u} $\mathfrak{t}$ est un tore de Malcev et $\mathfrak{n}$ le nilradical.
Reprenons ici l'\'{e}tude faite dans \cite{A.G} et \cite{G.A} permettant de
donner une description pr\'{e}cise des alg\`{e}bres r\'{e}solubles rigides
et leur classification lorsque le  nilradical est filiforme.

\bigskip

\noindent%

\begin{defi}
Soit $\g=\frak{t} \oplus \frak{n}$ une alg\`ebre de Lie r\'esoluble rigide de multiplication $\mu$. Un vecteur $X \in \frak{t}$ est dit r\'egulier
si la dimension de l'espace
$$V_0(X)=\{Y \in \g, \  \mu(X,Y)=0\}$$
est minimale c'est-\`a-dire $dimV_0(X) \leq dim V_0(Z)$ pour tout $Z \in \frak{t}$.
\end{defi}

\bigskip

Supposons que $\mathfrak{g}$ ne soit pas nilpotente.\ Dans ce cas $\mathfrak{%
t}$ n'est pas trivial.\ Soit $X$ un vecteur r\'{e}gulier et posons $p=\dim
V_{0}(X)$. Consid\'{e}rons une base $\left\{ X,Y_{1},\ldots
,Y_{n-p},X_{1},\ldots ,X_{p-1}\right\} $ de vecteurs propres de $adX$ (qui
par hypoth\`{e}se est diagonalisable) telle que $\left\{ X,X_{1},\ldots
,X_{p-1}\right\} $ soit une base de $V_{0}(X)$ et $\left\{ Y_{1},\ldots
,Y_{n-p}\right\} $ une base du nilradical $\mathfrak{n}$. On suppose \'{e}%
galement que $\left\{ X,X_{k_{0}+1},\ldots ,X_{p-1}\right\} $ soit une base
de $\mathfrak{t.}$ On notera par $(S)$ le syst\`{e}me de racines associ\'{e} 
\`{a} $adX$.\ Il est d\'{e}fini par%
\[
\left\{ 
\begin{array}{l}
x_{i}+x_{j}=x_{k}\quad \text{si la composante de }\mu (X_{i},X_{j})\text{
sur }X_{k}\text{ est non nulle} \\ 
y_{i}+y_{j}=y_{k}\quad \text{si la composante de }\mu (Y_{i},Y_{j})\text{
sur }Y_{k}\text{ est non nulle} \\ 
x_{i}+y_{j}=y_{k}\quad \text{si la composante de }\mu (X_{i},Y_{j})\text{
sur }Y_{k}\text{ est non nulle} \\ 
y_{i}+y_{j}=x_{k}\quad \text{si la composante de }\mu (Y_{i},Y_{j})\text{
sur }X_{k}\text{ est non nulle}%
\end{array}%
\right. 
\]

\bigskip

\bigskip

\noindent%

\begin{theo}
Si $rang(S) \neq dim(\frak{n})-1$, alors $\g$ n'est pas rigide.
\end{theo}

\bigskip

\bigskip La d\'{e}monstartion est donn\'{e}e dans \cite{A.G}.\ Les cons\'{e}%
quences sont nombreuses:

\begin{coro}
Si $\mathfrak{g=t\oplus n}$ est rigide, alors

\begin{itemize}
\item $\mathfrak{t}$ est un tore de Malcev maximal.

\item Il existe un vecteur r\'{e}gulier $X\in \mathfrak{t}$ tel que les
valeurs propres de $adX$ soit enti\`{e}res.

\item Le nilradical est d\'{e}fini par une solution isol\'{e}e du syst\`{e}%
me polynomial de Jacobi d\'{e}fini par les racines.
\end{itemize}
\end{coro}

Notons que la deuxi\`{e}me condition n'implique pas que l'alg\`{e}bre de Lie
rigide soit rationnelle.\ Il existe en effet des exemples d'alg\`{e}bres
rigides non rationnelles.\ La troisi\`{e}me propri\'{e}t\'{e} signifie que,
une fois donn\'{e}es les racines de $adX,$ les constantes de structure de $n$
sont donn\'{e}es par des conditions (r\'{e}duites) de Jacobi.\ L'alg\`{e}bre
est rigide si le nilradical correspond \`{a} une solution isol\'{e}e.\ 

\bigskip

Comme cons\'{e}quence, signalons la classification des alg\`{e}bres de Lie
rigides r\'{e}solubles de dimension inf\'{e}rieure ou \'{e}gale \`{a} $8$
(voir \cite{G.A}) et la classification g\'{e}n\'{e}rale des alg\`{e}bres de
Lie r\'{e}solubles dont le nilradical est filiforme (voir aussi \cite{G.A}).

\newpage

\bigskip

\bigskip

\bigskip

\part{Structures g\'{e}om\'{e}triques sur les alg\`{e}bres de Lie}

Soit $\mathfrak{g}$ une alg\`{e}bre de Lie r\'{e}elle ou complexe de
dimension finie. Soit $G$ un groupe de Lie connexe d'alg\`{e}bre de Lie $%
\mathfrak{g}$.\ Les \'{e}l\'{e}ments de $\mathfrak{g}$ sont donc les champs
de vecteurs invariants \`{a} gauche sur $G$ et le dual vectoriel $\Lambda
^{p}(\mathfrak{g}^{\ast })$ de $\mathfrak{g}$ est l'espace des formes de
Pfaff invariantes \`{a} gauche sur $G$. Rappelons que si $\{X_{1},\ldots
X_{n}\}$ est une base de $\mathfrak{g}$ et si $\{\omega _{1},\ldots ,\omega
_{n}\}$ en est la base duale de $\mathfrak{g}^{\ast },$ alors les \'{e}%
quations de Maurer-Cartan sont donn\'{e}es \`{a} partir de 
\[
\lbrack X_{i},X_{j}]=\sum_{k=1}^{n}C_{ij}^{k}X_{k} 
\]%
par%
\[
d\omega _{k}=\sum_{1\leq i<j\leq n}^{{}}C_{ij}^{k}\omega _{i}\wedge
\omega _{j} 
\]%
o\`{u} $d$ est la diff\'{e}rentielle ext\'{e}rieure des formes invariantes 
\`{a} gauche sur $G$. Si $\Lambda (\mathfrak{g}^{\ast })=\oplus _{p\in 
\mathbb{N}}\Lambda ^{p}(\mathfrak{g}^{\ast })$ est l'alg\`{e}bre ext\'{e}%
rieure sur $\mathfrak{g}^{\ast }$, cette diff\'{e}rentielle est donc un
morphisme gradu\'{e} de degr\'{e} $1$ :%
\[
d:\Lambda ^{p}(\mathfrak{g}^{\ast })\rightarrow \Lambda ^{p+1}(\mathfrak{g}%
^{\ast }). 
\]%
Nous la consid\`{e}rerons donc de mani\`{e}re \'{e}quivalente, soit comme
une diff\'{e}rentielle soit comme un morphisme lin\'{e}aire gradu\'{e}.

\bigskip

Les structures g\'{e}om\'{e}triques que nous allons aborder dans ce qui suit
seront d\'{e}finies sur les alg\`{e}bres de Lie.\ Elles correspondent en
fait \`{a} des structures g\'{e}om\'{e}triques invariantes \`{a} gauche sur $%
G$. On s'int\'{e}erssera souvent au cas o\`{u} l'alg\`{e}bre de Lie est
nilpotente.\ En effet dans ce cas, si cette alg\`{e}bre est rationnelle,
c'est-\`{a}-dire si elle admet une base par rapport \`{a} laquelle les
constantes de structure sont rationnelles, alors le groupe de Lie nilpotent
simplement connexe et connexe associ\'{e} admet un sous-groupe discret $%
\Gamma $ tel que le quotient $G/\Gamma $ soit une vari\'{e}t\'{e} diff\'{e}%
rentielle compacte appel\'{e}e nilvari\'{e}t\'{e}.\ Les structures d\'{e}%
finies sur $\mathfrak{g}$ correspondant \`{a} des structures invariantes 
\`{a} gauche sur $G$, donnent des structures diff\'{e}rentiables analogues
sur la vari\'{e}t\'{e} compacte $G/\Gamma $ d\`{e}s que ces structures sont
invariantes \`{a} droite par $\Gamma .$

\section{Structures symplectiques}

Soit $\mathfrak{g}$ une alg\`{e}bre de Lie r\'{e}elle (ou complexe) de
dimnsion $2n$. Une structure symplectique sur $\mathfrak{g}$ est donn\'{e}e
par une $2-$forme $\omega $ v\'{e}rifiant 
\[
\left\{ 
\begin{array}{c}
d\omega =0 \\ 
\omega ^{n}\neq 0%
\end{array}%
\right. 
\]%
o\`{u} $\omega ^{n}=\omega \wedge \omega \wedge \ldots \wedge \omega $ ($n$
fois) et%
\[
d\omega (X,Y,Z)=\omega (X,[Y,Z])+\omega (Y,[Z,X])+\omega (Z,[X,Y]) 
\]%
la multiplication de $\mathfrak{g}$ est not\'{e}e ici comme en g\'{e}om\'{e}%
trie diff\'{e}rentielle par le crochet. Les alg\`{e}bres de Lie frob\'{e}%
niusiennes sont munies d'une structure symplectique.\ En effet si $\alpha
\in \mathfrak{g}^{\ast }$ v\'{e}rifie $(d\alpha )^{n}\neq 0$, alors $\omega
=d\alpha $ est une forme symplectique (dite exacte). Le probl\`{e}me
d'existence d'une structure symplectique sur une alg\`{e}bre de Lie est
toujours d'actualit\'{e}. Par exemple, d'apr\`{e}s \cite{Go3}, il n'existe
pas de structure frob\'{e}niusienne sur une alg\`{e}bre de Lie nilpotente.
En effet, dans ce cas le centre de $\mathfrak{g}$ est au moins de dimension $%
1$, et tout vecteur du centre est dans le noyau de $d\alpha .~$La forme $%
\omega =d\alpha $ est donc d\'{e}g\'{e}n\'{e}r\'{e}e et ne peut v\'{e}rifier 
$\omega ^{n}\neq 0.$ Toutefois la classification des alg\`{e}bres de Lie
nilpotentes de dimension $6$ munie d'une structure symplectique est connue.
On peut la consulter dans \cite{G.K}

Il existe un proc\'{e}d\'{e} de construction des alg\`{e}bres de Lie munies
d'une forme symplectique, appel\'{e} le proc\'{e}d\'{e} de double extension
et d\'{e}fini par Alberto Medina et Philippe\ Revoy. Soit $(\mathfrak{%
g,\omega )}$ une alg\`{e}bre de Lie munie d'une forme symplectique.\ On dira
que l'alg\`{e}bre est symplectique. Alors le produit $\circledast $ d\'{e}%
fini par%
\[
\omega (X\circledast Y,Z)=-\omega (Y,[X,Z]) 
\]%
correspondant \`{a} l'adjoint pour la forme non d\'{e}g\'{e}n\'{e}r\'{e}e $%
\omega $ de l'application lin\'{e}aire $adX$ est un produit sym\'{e}trique
gauche 
\[
(X\circledast Y)\circledast Z-X\circledast (Y\circledast Z)=(Y\circledast
X)\circledast Z-Y\circledast (X\circledast Z) 
\]%
tel que 
\[
X\circledast Y-Y\circledast X=[X,Y] 
\]%
pour tout $X,Y,Z\in \mathfrak{g.}$ On dit alors que $\mathfrak{g}$ est munie
d'une structure affine, structure que l'on va \'{e}tudier deux paragraphes
plus loin.\ Soit $D$ une d\'{e}rivation de $\mathfrak{g.}$ Alors
l'application bilin\'{e}aire $f$ sur $\mathfrak{g}$ donn\'{e}e par%
\[
f(X,Y)=\omega (D(X),Y)+\omega (X,D(Y)) 
\]%
est un $2$-cocycle pour la cohomologie scalaire de $\mathfrak{g}$.\ Elle
permet donc de d\'{e}finir une extension centrale $E=\mathfrak{%
g\oplus }\mathbb{K}e$ o\`{u} $\mathbb{K=R}$ ou $\mathbb{C}$ en posant%
\[
[X+\lambda e,Y+\beta e]_{E}=[X,Y]+f(X,Y)e. 
\]%
On v\'{e}rifie ais\'{e}ment que ceci est un crochet de Lie sur $E$, que $e$
est dans le centre de $E$ et $\dim E=\dim \mathfrak{g+}1.$ Consid\'{e}rons 
\`{a} pr\'{e}sent l'application bilin\'{e}aire sur $\mathfrak{g}$ donn\'{e}e
par%
\[
\Omega (X,Y)=\omega (((D+D^{\ast })\circ D+D^{\ast }\circ (D+D^{\ast
}))(X),Y) 
\]%
o\`u $D^*$ est l'application adjointe de $D$ par rapport \`a la forme non d\'eg\'en\'er\'ee $\omega$. Cette application appartient \`{a} $Z^{2}(\mathfrak{g,}\mathbb{K)}$.
Supposons que $\Omega \in B^{2}(\mathfrak{g,}\mathbb{K)}$. Il existe alors $%
Z_{\Omega }\in \mathfrak{g}$ tel que 
\[
\Omega (X,Y)=\omega (Z_{\Omega },[X,Y]) 
\]%
pour tout $X,Y\in \mathfrak{g.}$ D\'{e}finissons alors la d\'{e}rivation $%
D_{1}$ de l'alg\`{e}bre de Lie $E$ en posant%
\begin{eqnarray*}
D_{1}(X) &=&-D(X)-\omega (Z_{\Omega },X)e,~X\in \mathfrak{g} \\
D_{1}(e) &=&0.
\end{eqnarray*}%
Cette d\'{e}rivation permet de construire une extension par d\'{e}rivation
de $E$ de dimension $\dim \mathfrak{g+}2$, qui est un produit semi-direct%
\[
\mathfrak{g}^{\prime }=(\mathfrak{g\oplus }\mathbb{K}e)\ltimes _{D_{1}}%
\mathbb{K}d 
\]%
de $E$ par un espace de dimension $1$, not\'{e} $\mathbb{K}d.$ Son crochet
est donn\'{e} par 
\begin{eqnarray*}
\lbrack X,Y]_{\mathfrak{g}^{\prime }} &=&[X,Y]_{E}\,\ \ X,Y\in E \\
\lbrack d,X]_{\mathfrak{g}^{\prime }} &=&-D_{1}(X)-\omega (Z_{\Omega },X)e
\end{eqnarray*}%
pour tout $X,Y\in E.$ Alors l'application bilin\'{e}aire%
\[
\omega _{1}:\mathfrak{g}^{\prime }\times \mathfrak{g}^{\prime }\rightarrow 
\mathfrak{g}^{\prime } 
\]%
donn\'{e}e par%
\[
\left\{ 
\begin{array}{l}
\omega _{1}\mid _{\mathfrak{g\times g}}=\omega \\ 
\omega _{1}(e,d)=1%
\end{array}%
\right. 
\]%
les autres produits non d\'{e}finis \'{e}tant nuls, est une forme
symplectique sur $\mathfrak{g}^{\prime }~$\ L'alg\`{e}bre symplectique $(%
\mathfrak{g}^{\prime },\omega _{1})$ est appel\'{e}e la double extension
symplectique de $(\mathfrak{g},\omega _{{}})$ au moyen de $D$ et $Z_{\Omega
}.$ Cette construction ne fonctionne que sous l'hypoth\`{e}se $\Omega \in
B^{2}(\mathfrak{g,}\mathbb{K)}$. Par contre, on montre que toute alg\`{e}bre
de Lie nilpotente symplectique de dimension $2n+2$ est une double extension
symplectique d'une alg\`{e}bre symplectique nilpotente de dimension $2n.$

\section{Structures complexes}

Soit $\mathfrak{g}$ une alg\`{e}bre de Lie r\'{e}elle de dimension paire $2n$
dont la multiplication est not\'{e}e $\mu $. Une structure complexe sur $%
\mathfrak{g}$ est donn\'{e}e par un endomorphisme $J:\mathfrak{g\rightarrow g%
}$ v\'{e}rifiant%
\[
\left\{ 
\begin{array}{l}
J^{2}=-Id \\ 
\mu (J(X),J(Y)=\mu (X,Y)+J(\mu (J(X),Y)+\mu (X,J(Y)))%
\end{array}%
\right. 
\]%
pour tout $X,Y\in \mathfrak{g.}$ Etablir l'existence d'une telle structure
est un probl\`{e}me assez difficile. Dans le cas nilpotent, S.Salamon donne
dans \cite{Sa} la classification des alg\`{e}bres de Lie nilpotentes r\'{e}%
elles de dimension inf\'{e}rieure ou \'{e}gale \`{a} $6$ admettant une telle
structure.\ Cette classification montre qu'il existe des alg\`{e}bres
nilpotentes n'admettant aucune structure complexe.\ Dans \cite{G.R3} on
montre le r\'{e}sultat suivant :

\begin{proposition}
Soit $\mathfrak{g}$ une alg\`{e}bre de Lie filiforme r\'{e}elle de dimension 
$2n$. Alors il n'existe aucune structure complexe sur $\mathfrak{g.}$
\end{proposition}

Rappelons que la suite caract\'{e}ristique $c(\mathfrak{g)}$ d'une alg\`{e}%
bre de Lie nilpotente est l'invariant \`{a} isomorphisme pr\`{e}s donn\'{e}
par%
\[
c(\mathfrak{g)=}\max \mathfrak{\{}c(X\mathfrak{),}X\in \mathfrak{g-}\mathcal{%
D}^{1}(\mathfrak{g)\}} 
\]%
o\`{u} $c(X)$ est la suite ordonn\'{e}e d\'{e}croissante des dimensions des
blocs de Jordan de l'op\'{e}rateur nilpotent $adX$. En particulier la classe
des alg\`{e}bres filiformes de dimension $2n$ est la classe des alg\`{e}bres
de caract\'{e}ristique $(2n-1,1).$ Une alg\`{e}bre est dite quasifiliforme
si sa caract\'{e}ristique est $(2n-2,1,1).$ En dimension $6$, il n'existe
qu'une seule classe d'alg\`{e}bre quasi-filiforme admettant une structure
complexe. On peut lire ce travail dans \cite{Ga.R} qui repose sur la notion
de structures complexes g\'{e}n\'{e}ralis\'{e}es que l'on pr\'{e}sente au
paragraphe suivant.

\section{Structures complexes g\'{e}n\'{e}ralis\'{e}es}

Les structures complexes g\'{e}n\'{e}ralis\'{e}es sont une nouvelle esp\`{e}%
ce de structures g\'{e}om\'{e}triques, introduites par Nigel Hitchin, et qui
contiennent les structures symplectiques et les structures complexes comme
cas extr\^{e}mes. Elles sont en g\'{e}n\'{e}ralement d\'{e}finies sur des
vari\'{e}t\'{e}s diff\'{e}rentables.\ Nous particularisons cette \'{e}tude
aux alg\`{e}bres de Lie, ceci correspondant aux structures complexes g\'{e}n%
\'{e}ralis\'{e}es invariantes \`{a} gauche sur des Groupes de Lie. Une pr%
\'{e}sentation d\'{e}taill\'{e}e est faite dans \cite{Ca.Gu}.

Soit $\g$ une alg\`ebre de Lie r\'eelle de dimension $2n$. Notons par $[X,Y]$ le crochet de $\g$. Si $\g^*$ est
 le dual vectoriel de $\g$, on d\'efinit
sur la somme directe externe $\g \oplus \g ^*$ une structure d'alg\`ebre 
de Lie en posant

$$\mu(X+\alpha ,Y+\beta )=[X,Y]+i(X)d\gamma +i(Y)d\alpha $$
o\`u $X,Y \in\ g$,
 $\alpha ,\beta  \in \g^*$ et $i(X)$ d\'esignant le produit int\'erieur, c'est-\`a-dire
 
$$i(X)d\alpha (Z)=d\alpha (X,Z)=-\alpha [X,Z].$$

Munissons cette alg\`ebre de Lie $\g \oplus \g^*$ du
 produit scalaire invariant (on dit que c'est une alg\`ebre de Lie quadratique) 
donn\'e  par
$$<X+\alpha ,Y+\beta >=\displaystyle \frac{1}{2}(\beta (X)+\alpha (Y)).$$
Ce produit scalaire 
est de signature $(2n,2n)$ (rappelons que $\g \oplus \g^*$ est de dimension $4n$). 

\begin{defi} 
On dit 
qu'un endomorphisme lin\'eaire
$$J:\g \oplus \g^* \longrightarrow \g \oplus \g^*$$
est une structure
 complexe g\'en\'eralis\'ee si

1. $J$ est une isom\'etrie de $<,>$, c'est-\`a-dire

$$<J(X+\alpha ),J(Y+\beta )>=<X+\alpha ,Y+\beta >$$

pour tout $X,Y \in \g$ et $\alpha ,\beta \in \g^*$,

2. Si $L$ est l'espace propre associ\'e \`a la valeur
 propre $i$ de $J$ sur l'espace complexe $(\g \oplus \g ^*) \otimes \mathbb{C}$, alors $L$ est un
 sous-espace isotrope maximal de $<,>$ (donc de dimension $2n$), involutif pour $\mu$, c'est-\`a-dire 
$\mu(L,L) \subset L.$

\end{defi}

On peut caract\'eriser les sous-espaces isotropes et isotropes maximaux en utilisant les alg\`ebres de Clifford.
Rappelons bri\`evement la d\'efinition.
Une alg\`ebre de Clifford est une alg\`ebre unitaire associative
qui est engendr\'ee par un espace 
vectoriel $V$ muni d'une forme quadratique $Q$
soumise \`a la condition
$$v^2 = Q(v)\ \rm{pour~tout}\ v\in V.$$ 
Soit $\varphi $ un spineur (on peut prendre par exemple un \'el\'ement de $\Lambda \g ^*.)$ 
L'action de Clifford de $\g\oplus \g^*$ sur le spineur $\phi $ est donn\'ee par

$$(X+\alpha )\bullet \phi=i(X)\phi +\alpha \wedge \phi.$$
Cette action correspond \`a une repr\'esentation des alg\`ebres de Clifford car

$$(X+\alpha )^2\bullet \phi=<X+\alpha ,X+\alpha >\phi.$$
A tout spineur $\phi$, faisons correspondre le sous-espace vectoriel $L_{\phi}$ de $\g \oplus \g^*$ donn\'e par

$$L_{\phi}=\{(X+\alpha ) \in \g \oplus \g^*, (X+\alpha )\bullet \phi =0\}.$$

On dira que $\phi$ est un spineur pur lorsque $dim L_{\phi} = 2n.$ Dans ce cas $L_{\phi}$ est un 
sous-espace isotrope maximal pour $<,>$. Inversement, si $L$ est un sous-espace isotrope maximal, 
alors l'ensemble des spineurs

\[
U_{L}=\left\{ \varphi \in \Lambda \mathfrak{g}^{\ast },\quad L=L_{\varphi
}\right\} 
\]%
est une droite de spineurs (purs) engendr\'{e}e par un spineur du type%
\[
e^{B+i\omega }\theta _{1}\wedge \ldots \wedge \theta _{k}
\]%
o\`{u} $B$ et $\omega $ sont des formes r\'{e}elles et $\theta _{i}$ des
formes complexes de degr\'{e} $1$. Dans cette expression $e^{A}$ d\'{e}signe 
$Id+B+B\wedge B/2+$\ldots 

\begin{defi}
On dira que la structure complexe g\'{e}n\'{e}ralis\'{e}e $J$ est de type $k$
si $k$ est la codimension de la projection de $L$ sur $\mathfrak{g}$.
\end{defi}

Notons que si $J$ est de type $k$, alors le spineur pur engendrant $L$ s'%
\'{e}crit $e^{B+i\omega }$ $\theta _{1}\wedge \ldots \wedge \theta _{k}$ (le
m\^{e}me $k$).

\bigskip 

\textbf{Exemples. }

1.\ Soit $j$ une structure complexe sur l'alg\`{e}bre de Lie $\mathfrak{g}$ r%
\'{e}elle de dimension $2n.$ Alors l'endomorphisme de $\mathfrak{g\oplus g}%
^{\ast }$ donn\'{e} par%
\[
J_{j}(X+\alpha )=-j(X)+j^{\ast }(\alpha )
\]%
o\`{u} $j^{\ast }$ d\'{e}signe la transpos\'{e}e de $j$ est une structure
complexe g\'{e}n\'{e}ralis\'{e}e de type (maximal) $n.$ Dans ce cas, si $%
T_{0,1}$ est l'espace propre de $j$ associ\'{e} \`{a} la valeur propre $i$%
, on a%
\[
L=(T_{0,1})\oplus (T_{0,1})^{\ast }
\]%
et le spineur pur d\'{e}finissant $L$ est donn\'{e} par%
\[
\rho =e^{B}\theta _{1}\wedge \ldots \wedge \theta _{n}
\]%
R\'{e}ciproquement, toute structure symplectique g\'{e}n\'{e}ralis\'{e}e de
type $n$ correspond \`{a} une structure complexe sur $\mathfrak{g.}$

\bigskip 

2. Supposons que $\mathfrak{g}$ admette une forme symplectique $\omega .\ $%
Nous pouvons consid\'{e}rer $\omega $ comme un isomorphisme%
\[
\omega :\mathfrak{g\rightarrow g}^{\ast }
\]%
donn\'{e} par $\omega (X)=i(X)\omega .$ Soit l'endomorphisme de $\mathfrak{%
g\oplus g}^{\ast }$ donn\'{e} par%
\[
J_{\omega }(X+\alpha )=i(X)\omega -\omega ^{-1}(\alpha ).
\]%
Il d\'{e}finit une structure complexe g\'{e}n\'{e}ralis\'{e}e de type $0$.
Dans ce cas%
\[
L=\left\{ X+\alpha ,\omega ^{-1}(\alpha )=iX,\omega (X)=i\alpha \right\}
=\left\{ X-i\omega (X),X\in \mathfrak{g\otimes }\mathbb{C}\right\} .
\]%
R\'{e}ciproquement, toute  structure symplectique g\'{e}n\'{e}ralis\'{e}e de
type $0$ correspond \`{a} une structure symplectique sur $\mathfrak{g.}$

\bigskip 

\section{Structures affines}

Une structure affine sur une alg\`{e}bre de Lie $\mathfrak{g}$ dont la
multiplication est not\'{e}e $\mu $ est donn\'{e}e par une application bilin%
\'{e}aire%
\[
\nabla :\mathfrak{g\otimes g\rightarrow g} 
\]%
v\'{e}rifiant%
\[
\left\{ 
\begin{array}{l}
\nabla (X,Y)-\nabla (Y,X)=\mu (X,Y) \\ 
\nabla (X,\nabla (Y,Z))-\nabla (Y,\nabla (X,Z))=\nabla (\nabla
(X,Y),Z)-\nabla (\nabla (Y,X),Z)%
\end{array}%
\right. 
\]%
pour tout $X,Y,Z\in \mathfrak{g.}$ Cette op\'{e}ration correspond en fait 
\`{a} la donn\'{e}e d'une connexion affine sans courbure ni torsion
invariante \`{a} gauche sur un groupe de Lie d'alg\`{e}bre de Lie $\mathfrak{%
g.\ }$Les identit\'{e}s ci-dessus signifie que l'alg\`{e}bre $(\mathfrak{%
g,\nabla )}$ est une alg\`{e}bre sym\'{e}trique gauche encore appel\'{e}e alg%
\`{e}bre de Pr\'{e}-Lie. Toute alg\`{e}bre associative est donc de Pr\'{e}%
-Lie.\ On en d\'{e}duit que toute alg\`{e}bre de Lie associ\'{e}e \`{a} une
alg\`{e}bre associative admet une structure affine. Si $\nabla $ est
commutative, alors c'est une multiplication associative commutative et l'alg%
\`{e}bre de Lie associ\'{e}e est ab\'{e}lienne. Il y a donc une
correspondance bijective entre les alg\`{e}bres associatives commutatives de
dimension $n$ et les structures affines sur l'alg\`{e}bre ab\'{e}lienne de
dimension $n$.\ Par exemple, on pourra consulter cette classification pour $%
n=3$ dans \cite{G.R2}.\ Pour la dimension $4$ on pourra consulter \cite{B}.

Le probl\`{e}me d'existence d'une structure affine sur une alg\`{e}bre
nilpotente a longtemps \'{e}t\'{e} dict\'{e} par la conjecture de Milnor
stipulant que toute alg\`{e}bre de Lie nilpotente \'{e}tait munie d'une
structure affine.\ Ceci \ impliquait en particulier que toute nilvari\'{e}t%
\'{e} \'{e}tait affine (munie d'une connexion affine sans courbure ni
torsion).\ En fait Y.\ Besnoit a mis en \'{e}vidence une nilvari\'{e}t\'{e}
de dimension $11$ sans structure affine.\ L'alg\`{e}bre de Lie
correspondante est filiforme \cite{Bes} . Dans ce qui suit, on va donner des
exemples de structure affine construite en fonction de certaines propri\'{e}t%
\'{e}s de $\mathfrak{g.}$

\bigskip

\begin{itemize}
\item Supposons que $\mathfrak{g}$ soit de dimension $2n$ et munie d'une
forme symplectique $\omega .$ Pour tout $X\in \mathfrak{g,}$ soit $f_{X}$
l'adjoint de $adX$ pour la forme $\omega :$%
\[
\omega (\mu (Y,X),Z)=-\omega (Y,f_{X}(Z)) 
\]%
pour tout $Y,Z\in \mathfrak{g.}$ On v\'{e}rifie assez facilement que l'op%
\'{e}ration%
\[
\nabla (X,Y)=f_{X}(Y) 
\]%
munie $\mathfrak{g}$ d'une structure affine.

\item Supposons que $\mathfrak{g}$ soit munie d'une d\'{e}rivation $f$
inversible. Ceci implique n\'{e}cessairement que $\mathfrak{g}$ soit
nilpotente. Rappelons qu'une d\'{e}rivation est un endomorphisme lin\'{e}%
aire $f:\mathfrak{g\rightarrow g}$ v\'{e}rifiant%
\[
f(\mu (X,Y))=\mu (f(X),Y)+\mu (X,f(Y)) 
\]%
pour tout $X,Y\in \mathfrak{g}$. L'application bilin\'{e}aire d\'{e}finie par%
\[
\nabla (X,Y)=f^{-1}(\mu (f(X),Y)) 
\]%
est une structure affine sur $\mathfrak{g.\ }$La d\'{e}monstration est laiss%
\'{e}e ici aussi en exercice.

\item On peut \'{e}galement construire une telle structure en consid\'{e}%
rant une d\'{e}rivation $f$ inversible ou non mais en supposant que la
restriction de la d\'{e}rivation $f$ \`{a} l'alg\`{e}bre d\'{e}riv\'{e}e
soit inversible. Dans ce cas l'op\'{e}ration $\nabla $ est d\'{e}finie comme
ci-dessus.

\item Supposons que $\mathfrak{g}$ soit munie d'un op\'{e}rateur de
Baxter-Yang $R$, c'est-\`{a}-dire d'un endomorphisme v\'{e}rifiant%
\[
\mu (R(X),R(Y)=R(\mu (R(X),Y)+\mu (X,R(Y))). 
\]

Alors l'op\'{e}ration%
\[
\nabla (X,Y)=\mu (R(X),Y) 
\]%
est une multiplication d'alg\`{e}bre de Pr\'{e}-Lie. On en d\'{e}duit que
l'alg\`{e}bre de Lie $\mathfrak{g}^{\prime }$ de multiplication $\mu
^{\prime }$ donn\'{e}e par%
\[
\mu ^{\prime }(X,Y)=\nabla (X,Y)-\nabla (Y,X) 
\]%
est munie d'une structure affine. En g\'{e}n\'{e}ral $\mu ^{\prime }\neq \mu 
$ sauf si l'application $Id-R$ est \`{a} valeurs dans le centre de $%
\mathfrak{g.}$

\item Dans le m\^{e}me ordre d'id\'{e}e, si $\mathfrak{g}$ est munie d'un op%
\'{e}rateur de Baxter-Rota $R$, c'est-\`{a}-dire d'un endomorphisme v\'{e}%
rifiant%
\[
\mu (R(X),R(Y)+\mu (X,Y)=R(\mu (R(X),Y)+\mu (X,R(Y))) 
\]%
alors l'application%
\[
\nabla (X,Y)=\mu (R(X),Y) 
\]%
donne un r\'{e}sultat analogue au pr\'{e}dent.

\item Si $\mathfrak{g}$ est une alg\`{e}bre de Lie r\'{e}elle de dimension $%
2n$ munie d'une structure complexe int\'{e}grable, c'est-\`{a}-dire d'un
endomorphisme $J:\mathfrak{g\rightarrow g}$ v\'{e}rifiant $J^{2}=-Id$ et la
condition de Nijenhuis%
\[
\mu (J(X),J(Y)=\mu (X,Y)+J(\mu (J(X),Y)+\mu (X,J(Y))) 
\]%
alors le morphisme $R=-iJ$ est un op\'{e}rateur de Baxter-Rota.\ On en d\'{e}%
duit que%
\[
\nabla (X,Y)=\mu (J(X),Y) 
\]%
d\'{e}finit une structure affine sur l'alg\`{e}bre de Lie de multiplication $%
\mu ^{\prime }(X,Y)=\nabla (X,Y)-\nabla (Y,X).$
\end{itemize}

\section{Espaces homog\`{e}nes r\'{e}ductifs}

\noindent Soient $G$ un groupe de Lie connexe et $H$ un sous-groupe de Lie.\
On dit que l'espace homog\`{e}ne $M=G/H$ est r\'{e}ductif si l'alg\`{e}bre
de Lie $\mathfrak{g}$ de $G$ peut se d\'{e}composer en une somme directe
vectorielle%
\[
\mathfrak{g=h\oplus m} 
\]%
o\`{u} $\mathfrak{h}$ est l'alg\`{e}bre de Lie de $H$ et $\mathfrak{m}$ un
sous espace vectoriel v\'{e}rifiant%
\[
ad(H)\mathfrak{m\subset m} 
\]%
ce qui est \'{e}quivalent si $H$ est connexe \`{a}%
\[
\lbrack \mathfrak{h,m]\subset m.} 
\]%
Ici, on note par $[X,Y]$ la multiplication de $\mathfrak{g}$, pour rester
conforme aux \'{e}critures classiques en g\'{e}om\'{e}trie diff\'{e}%
rentielle. De telles vari\'{e}t\'{e}s admettent toujours des connexions
invariantes par $G$ et une unique connexion affine $\nabla $ sans torsion $G$%
-invariante qui est compl\`{e}te. En ce qui concerne les m\'{e}triques
(riemanniennes ou pseudo-riemanniennes), il y a une correspondance bijective
entre l'ensembles des m\'{e}triques sur l'espace homog\`{e}ne r\'{e}ductif $%
M=G/H$ invariante par \ $G$ et les formes bilin\'{e}aires sym\'{e}triques
non d\'{e}g\'{e}n\'{e}r\'{e}es $B$ sur $\mathfrak{m}$ qui sont $ad(H)-$%
invariantes, c'est-\`{a}-dire, si $H$ est connexe, qui v\'{e}rifient%
\[
B([Z,X],Y)+B(X,[Z,Y])=0 
\]%
pour $X,Y\in \mathfrak{m}$ et $Z\in \mathfrak{h.}$ La connexion riemannienne
associ\'{e}e co\"{\i}ncide avec la connexion affine canonique sans torsion
si et seulement si%
\[
B([Z,X]_{\mathfrak{m}},Y)+B(X,[Z,Y]_{\mathfrak{m}})=0 
\]%
pour $X,Y,Z\in \mathfrak{m}$ o\`{u} $[,]_{\mathfrak{m}}$ d\'{e}signe la
projection sur $\mathfrak{m}$ du crochet. Dans ce cas, $M$ est appel\'{e}
espace homog\`{e}ne riemannien naturellement r\'{e}ductif.

\subsection{Espaces sym\'{e}triques}

La classe la plus connue et la plus \'{e}tudi\'{e}e d'espaces homog\`{e}nes r%
\'{e}ductifs et dans le cas riemannien, naturellement r\'{e}ductifs est
celle des espaces sym\'{e}triques. Un espace sym\'{e}trique est la donn\'{e}%
e d'un espace homog\`{e}ne $M=G/H$ et d'un automorphisme involutif $\sigma $
de $G$ tel que $H$ soit compris entre le sous groupe $G_{\sigma }$ des
points fixes de $G$ par $\sigma $ et sa composante connexe passant par
l'identit\'{e} de $G$. Ceci revient \`{a} se donner, en tout point $x\in M$,
un diff\'{e}omorphisme involutif, appel\'{e}e sym\'{e}trie et not\'{e}e $%
s_{x},$ tel que $x$ soit un point fixe isol\'{e}. Dans ce cas, l'alg\`{e}bre
de Lie $\mathfrak{g}$ est sym\'{e}trique, c'est-\`{a}-dire s'\'{e}crit, d'apr%
\`{e}s la d\'{e}composition r\'{e}ductive%
\[
\mathfrak{g=h\oplus m} 
\]%
avec%
\[
\left\{ 
\begin{array}{c}
\lbrack \mathfrak{h,h]\subset h} \\ 
\lbrack \mathfrak{h,m]\subset m} \\ 
\lbrack \mathfrak{m,m]\subset h}%
\end{array}%
\right. 
\]%
Notons que dans ce cas, on a n\'{e}cessairement $ad(H)\mathfrak{m\subset m}$
et donc un espace sym\'{e}trique est bien r\'{e}ductif. De plus, la
connexion canonique associ\'{e}e \`{a} tout espace homog\`{e}ne r\'{e}ductif
est, dans le cas sym\'{e}trique, la seule connexion affine invariante par
les sym\'{e}tries. Consid\'{e}rons une alg\`{e}bre de Lie sym\'{e}trique.\
L'application $\rho $ d\'{e}finie par $\rho (X)=X$ pour tout $X\in \mathfrak{%
h}$ et $\rho (Y)=-Y$ pour tout $Y\in \mathfrak{m}$ est un automorphisme
involutif de $\mathfrak{g}$.\ Ainsi toute espace sym\'{e}trique d\'{e}finit
une alg\`{e}bre de Lie sym\'{e}trique naturellement munie d'un automorphisme
involutif.\ La r\'{e}ciproque n'est pas exacte.\ Toutefois, si on suppose
que le groupe de Lie $G$ est connexe et simplement connexe, alors
l'automorphisme $\rho $ de $\mathfrak{g}$ induit un automorphisme involutif $%
\sigma $ de $G$ et pour tout sous-groupe de Lie $H$ compris entre le
sous-groupe des points fixes de $\sigma $ et sa composante de l'ident\'{e},
son alg\`{e}bre de Lie \'{e}tant $\mathfrak{h}$, l'espace homog\`{e}ne $G/H$
est sym\'{e}trique et l'automorphisme associ\'{e} est $\sigma .$

Un espace sym\'{e}trique est dit riemannien s'il est muni d'une m\'{e}trique 
$G$-invariante pour laquelle les sym\'{e}tries $s_{x}$ sont des isom\'{e}%
tries pour tout $x\in M=G/H.$ En fait, dans ce cas, $G$ correspond au plus
grand groupe connexe des isom\'{e}tries de $M$ et $H$ le sous-groupe
d'isotropie en un point fix\'{e} de $M$. Si la m\'{e}trique est
riemannienne, alors $H$ est compact,~si elle est pseudo-riemannienne ceci
n'est pas toujours vrai. Dans tous les cas, on suppose que le sous-groupe $%
ad(H)$ des transformations lin\'{e}aires de $\mathfrak{g}$ est compact. On
en d\'{e}duit que $\mathfrak{g}$ admet un produit scalaire $ad(H)$ invariant
tel que $\mathfrak{h}$ et $\mathfrak{m}$ soient orthogonaux. En restriction 
\`{a} $\mathfrak{m}$, on retrouve la m\'{e}trique $G-$invariante munissant $%
G/H$ d'une structure d'espace sym\'{e}trique riemannien. La connexion
riemannien co\"{\i}ncide n\'{e}cessairement avec la connexion naturelle sans
torsion de l'espace homog\`{e}ne r\'{e}ductif $G/H$.\ L'espace riemannien $%
G/H$ est donc naturellement r\'{e}ductif. Les alg\`{e}bres de Lie sym\'{e}%
triques correspondantes aux espaces sym\'{e}triques riemanniens sont les alg%
\`{e}bres de Lie orthogonales sym\'{e}triques. Ce sont des alg\`{e}bres de
Lie sym\'{e}triques $\mathfrak{g=h\oplus m}$ telles que l'alg\`{e}bre de Lie 
$ad(\mathfrak{h)}$ soit compacte.\ Si $H$ n'a qu'un nombre fini de
composantes connexes, ceci est \'{e}quivalent \`{a} dire que le groupe $%
ad(H) $ est compact. Nous avons vu que $\mathfrak{g}$ admettait un produit
scalaire $B$ qui est donc $ad(\mathfrak{h)-}$invariant et pour lequel $%
\mathfrak{h}$ et $\mathfrak{m}$ sont orthogonaux. L'invariance de $B$ se
traduit par%
\[
B([X,Y],Z)+B(Y,[X,Z])=0 
\]%
pour tout $X\in \mathfrak{h}$ et $Y,Z\in \mathfrak{m}$. On en d\'{e}duit la
structure d'une alg\`{e}bre sym\'{e}trique orthogonale.\ Supposons l'espace
sym\'{e}trique riemannien (et non pseudo-riemannien).\ Si le centre de $%
\mathfrak{g}$ ne rencontre pas $\mathfrak{m}$, alors $\mathfrak{g}$ est la
somme directe de deux alg\`{e}bres sym\'{e}triques orthogonales $\mathfrak{g}%
_{1}=\mathfrak{h}_{1}\oplus \mathfrak{m}_{1}$ et $\mathfrak{g}_{2}=\mathfrak{%
h}_{2}\oplus \mathfrak{m}_{2}$ avec $[\mathfrak{m}_{1},\mathfrak{m}_{1}]=0$
et $\mathfrak{g}_{2}$ semi-simple. L'\'{e}tude des alg\`{e}bres de Lie sym%
\'{e}triques orthogonales se ram\`{e}ne donc \`{a} la classe des alg\`{e}%
bres semi-simples. \ Si l'alg\`{e}bre sym\'{e}trique orthogonale est simple,
la forme de Killing-Cartan $K$ de $\mathfrak{g}$ est d\'{e}finie (positive
ou n\'{e}gative) sur $\mathfrak{m.~}$Elle est dite de type compact si $K$
est d\'{e}finie n\'{e}gative sur $\mathfrak{m}$, et de type non-compact si
elle est d\'{e}finie positive. Dans le premier cas, l'espace $G/H$ est sym%
\'{e}trique riemannien compact. Leur classification se d\'{e}crit enti\`{e}%
rement par celle des alg\`{e}bres simples sym\'{e}triques orthogonales
compactes.\ Elle est due \`{a} Elie Cartan.\ Dans la liste suivante on donne
le couple $(\mathfrak{g},\mathfrak{h}).$

AI : $(su(n),so(n))$

AII : $(su(2n),sp(n))$

AIII : $(su(p+q),su(p)\oplus su(q)$

BDI : $(so(p+q),so(p)\oplus so(q)$

DIII : $(so(2n),u(n))$

CI : $(sp(n),u(n))$

CII : $(sp(p+q),sp(p)\oplus sp(q)$

EI : ($E(6),sp(4))$

EII ($E(6),su(6)\oplus su(2))$

EIII : ($E(6),so(10)\oplus so(2))$

EIV : ($E(6),F(4))$

EV : ($E(7),su(8))$

EVI : ($E(7),so(12)\oplus su(2))$

EVII : ($E(7),E(6)\oplus so(2))$

EVIII: ($E(8),so(16))$

EIX : ($E(8),E(7)\oplus su(2))$

FI : ($F(4),sp(3)\oplus su(2))$

FII: ($F(4),so(9))$

G : ($G(2),su(2)\oplus su(2))$

\subsection{Espaces homog\`{e}nes r\'{e}ductifs non-sym\'{e}triques}

Soit $\Gamma $ un groupe ab\'{e}lien fini. Un espace homog\`{e}ne $M=G/H$
est dit $\Gamma -$sym\'{e}trique s'il existe un sous-groupe $\Gamma _{G}$ de 
$Aut(G)$ isomorphe \`{a} $\Gamma $ tel que $H$ soit compris entre l'ensemble
des points fixes de tous les \'{e}l\'{e}ments $\sigma _{\gamma }$ de $\Gamma
_{G}$ et sa composante connexe passant par l'\'{e}l\'{e}ment neutre. Dans ce
cas, il existe pour tout $x\in M\,\ $un sous-groupe du groupe des diff\'{e}%
omorphismes de $M$ isomorphe \`{a} $\Gamma $. Les \'{e}l\'{e}ments sont appel%
\'{e}s les sym\'{e}tries et not\'{e}es $s_{\gamma ,x}$ et $x$ et un point
fixe isol\'{e} de ces sym\'{e}tries. L'alg\`{e}bre de Lie $\mathfrak{g}$ de $%
G$ admet une d\'{e}composition $\Gamma -$sym\'{e}trique, c'est-\`{a}-dire s'%
\'{e}crit%
\[
\mathfrak{g=\oplus }_{\gamma \in \Gamma }\text{ }\mathfrak{g}_{\gamma }
\]%
avec%
\[
\lbrack \mathfrak{g}_{\gamma _{1}},\mathfrak{g}_{\gamma _{2}}]\subset 
\mathfrak{g}_{\gamma _{1}\gamma _{2}}.
\]%
Une telle alg\`{e}bre est appel\'{e}e $\Gamma $-sym\'{e}trique. Cette
graduation est d\'{e}finie par un sous-groupe d'automor- -phisme de $\mathfrak{g%
}$ isomorphe \`{a} $\Gamma $ et si $G$ est connexe et simplement connexe,
cette graduation d\'{e}finit une structure $\Gamma -$sym\'{e}trique sur $G/H.
$ Si on pose $\mathfrak{m=\oplus }_{\gamma \neq 1_{\Gamma }}$ $\mathfrak{g}%
_{\gamma },$ on voit que $\mathfrak{g=g}_{1_{\Gamma }}\mathfrak{\oplus m}$
et cette d\'{e}composition donne une structure d'espace homog\`{e}ne r\'{e}%
ductif sur $G/H.$ Contrairement aux espaces sym\'{e}triques, la connexion
canonique de Nomizu et la connexion canonique sans torsion d'un espace homog%
\`{e}ne r\'{e}ductif ne co\"{\i}ncident pas. Dans \cite{Ba.G}, en utilisant
les graduations des alg\`{e}bres simples par des groupes finis, on donne la
classification des espaces homog\`{e}nes r\'{e}ductifs $\mathbb{Z}_{2}^{2}$%
-sym\'{e}triques lorsque $\mathfrak{g}$ est simple et non exceptionnelle. On
en d\'{e}duit en particulier la classification des alg\`{e}bres simples
compactes $\mathbb{Z}_{2}^{2}$-sym\'{e}triques non exceptionnelles :

$(su(2n),su(n))$

$(su(k_{1}+k_{2}),su(k_{1})\oplus su(k_{2})\oplus \mathbb{C})$

$(su(k_{1}+k_{2}+k_{3}),su(k_{1})\oplus su(k_{2})\oplus su(k_{3})\oplus 
\mathbb{C}^{2})$

$(su(k_{1}+k_{2}+k_{3}+k_{4}),su(k_{1})\oplus su(k_{2})\oplus
su(k_{3})\oplus su(k_{4})\oplus \mathbb{C}^{3})$

$(su(n),so(n))$

$(su(2m),sp(m))$

$(su(k_{1}+k_{2}),so(k_{1})\oplus so(k_{2}))$

$(su2(k_{1}+k_{2}),sp(2k_{1})\oplus sp(2k_{2}))$

$(so(k_{1}+k_{2}+k_{3}+k_{4}),so(k_{1})\oplus so(k_{2})\oplus
so(k_{3})\oplus so(k_{4}))$

$(so(4m),sp(2m))$

$(so(2m),so(m))$

$(sp(k_{1}+k_{2}+k_{3}+k_{4}),sp(k_{1})\oplus sp(k_{2})\oplus
sp(k_{3})\oplus sp(k_{4}))$

$(sp(4m),sp(2m))$

$(sp(2m),so(m))$

\subsection{Espaces riemanniens et pseudo-riemanniens $\Gamma $-sym\'{e}%
triques}

Soit $(M=G/H, \Gamma)$ un espace homog\`ene $\Gamma$-sym\'etrique.

\begin{defi}
Une m\'etrique riemannienne $g$ sur $M$ est dite adapt\'ee \`a la structure $%
\Gamma$-sym\'etrique si chacune des sym\'etries $s_{\gamma,x}$ est une
isom\'etrie.
\end{defi}

Si $\bigtriangledown_g $ est la connexion de Levi-Civita de $g$,
cette connexion ne co\"{\i}ncide pas en g\'en\'eral avec la connexion
canonique $\bigtriangledown $ de l'espace homog\`ene ($\Gamma$%
-sym\'etrique). Ces deux connexions co\"{\i}ncident si et seulement si $g$
est naturellement r\'eductive.

Par exemple dans le cas de la sph\`ere $S^3$ consid\'er\'ee comme espace $(%
\mathbb{Z}_2)^2$-sym\'etrique, les m\'etriques adapt\'ees \`a cette
structure sont les m\'etriques sur $SO(4)/Sp(2)$ invariantes par $SO(4)$
chacune \'etant d\'efinie par une forme bilin\'eaire sym\'etrique $B$ sur $%
so(4)$ qui est $ad(sp(2))$-invariante. Si $so(4)=sp(2) \oplus \mathfrak{g}_a
\oplus \mathfrak{g}_b \oplus \mathfrak{g}_c$ est la d\'ecomposition $(%
\mathbb{Z}_2)^2$-gradu\'ee correspondante, le fait de dire que sur $S^3$ les
sym\'etries $s_{\gamma,x}$ sont des isom\'etries est \'equivalent \`a dire
que les espaces $\mathfrak{g}_e,\mathfrak{g}_a,\mathfrak{g}_b, \mathfrak{g}%
_c $ sont deux \`a deux orthogonaux pour $B$. D\'ecrivons en d\'etail cette
graduation: 
\[
sp(2)=\left\{ \left( 
\begin{array}{llll}
0 & -a_2 & -a_3 & -a_4 \\ 
a_2 & 0 & -a_4 & a_3 \\ 
a_3 & a_4 & 0 & -a_2 \\ 
a_4 & -a_3 & a_2 & 0 \\ 
&  &  & 
\end{array}
\right) \right\}, \ \mathfrak{g}_a=\left\{ \left( 
\begin{array}{llll}
0 & 0 & 0 & x \\ 
0 & 0 & -x & 0 \\ 
0 & x & 0 & 0 \\ 
-x & 0 & 0 & 0 \\ 
&  &  & 
\end{array}
\right) \right\}, 
\]
\[
\mathfrak{g}_b=\left\{ \left( 
\begin{array}{llll}
0 & y & 0 & 0 \\ 
-y & 0 & 0 & 0 \\ 
0 & 0 & 0 & -y \\ 
0 & 0 & y & 0 \\ 
&  &  & 
\end{array}
\right) \right\} \ \mbox{\rm et} \ \mathfrak{g}_c=\left\{ \left( 
\begin{array}{llll}
0 & 0 & z & 0 \\ 
0 & 0 & 0 & z \\ 
-z & 0 & 0 & 0 \\ 
0 & -z & 0 & 0 \\ 
&  &  & 
\end{array}
\right) \right\}. 
\]
Si $\left\{ A_1 ,A_2, A_3, X,Y,Z \right\}$ est une base adapt\'ee \`a cette
graduation et si $\left\{ \alpha _1 ,\alpha _2, \alpha _3, \omega_1,
\omega_2, \omega_3 \right\}$ en est la base duale alors 
\[
B\mid_{ \mathfrak{g}_a \oplus \mathfrak{g}_b \oplus \mathfrak{g}_c } =
\lambda _1^2\omega _1^2+\lambda _2^2\omega _2^2+\lambda _3^2\omega _3^2. 
\]
La m\'etrique correspondante sera naturellement r\'eductive si et seulement
si $\lambda_1= \lambda_2 = \lambda_3$ et dans ce cas-l\`a elle correspond
\`a la restriction de la forme de Killing Cartan.

\subsection{Exemple : l'espace pseudo-riemannien $SO(2m)/Sp(m)$}

\subsubsection{La graduation $(\mathbb{Z}_{2})^{2}$-sym\'{e}trique}

Consid\'erons les matrices 
\[
S_m= \left( 
\begin{array}{ll}
0 & I_m \\ 
-I_n & 0%
\end{array}
\right), \ X_a=\left( 
\begin{array}{ll}
-1 & 0 \\ 
0 & 1%
\end{array}
\right),X_b=\left( 
\begin{array}{ll}
0 & 1 \\ 
1 & 0%
\end{array}
\right),X_c=\left( 
\begin{array}{ll}
0 & -1 \\ 
1 & 0%
\end{array}
\right). 
\]
Soit $M \in so(2m).$ Les applications 
\[
\begin{array}{l}
\tau_a(M)=J_a^{-1}MJ_a \\ 
\tau_b(M)=J_b^{-1}MJ_b \\ 
\tau_a(M)=J_c^{-1}MJ_c%
\end{array}%
\]
o\`u $J_a=S_m \otimes X_a, \ J_b=S_m \otimes X_b, \ J_c=S_m \otimes X_c$
sont des automorphismes involutifs de $so(2m)$ qui commutent deux \`a deux.
Ainsi $\left\{ Id,\tau _a,\tau _b,\tau _c \right\}$ est un sous groupe de $%
Aut (so(2m))$ isomorphe \`a $(\mathbb{Z}_2)^2$. Il d\'efinit donc une $(%
\mathbb{Z}_2)^2$-graduation 
\[
so(2m)= \mathfrak{g}_e \oplus \mathfrak{g}_a \oplus \mathfrak{g}_b \oplus 
\mathfrak{g}_c 
\]
o\`u 
\[
\begin{array}{l}
\mathfrak{g}_e= \left\{ M \in so(2m)\, / \, \tau _a(M)=\tau _b(M)=\tau _c
(M)=M \right\} \\ 
\mathfrak{g}_a= \left\{ M \in so(2m)\, / \, \tau _a(M) =\tau _c (M)=-M, \tau
_b(M)=M \right\} \\ 
\mathfrak{g}_b= \left\{ M \in so(2m)\, / \, \tau _b(M)=\tau _c(M)=-M, \tau
_a (M)=M \right\} \\ 
\mathfrak{g}_c= \left\{ M \in so(2m)\, / \, \tau _a(M)=\tau _b(M)=-M, \tau
_c (M)=M \right\}%
\end{array}%
\]
Ainsi 
\[
\mathfrak{g}_e= \left\{ \left( 
\begin{array}{rr|rr}
A_1 & B_1 & A_2 & B_2 \\ 
-B_1 & A_1 & B_2 & -A_2 \\ \hline
-^t \! A_2 & -^t \! B_2 & A_1 & B_1 \\ 
-^t \! B_2 & ^t \! A_2 & -B_1 & A_1%
\end{array}
\right) \ \mbox{\rm avec} \ 
\begin{array}{cc}
^t \! A_1 = -A_1, & ^t \! B_1 = B_1 \\ 
^t \! A_2 = A_2, & ^t \! B_2 = B_2 \\ 
& 
\end{array}
\right\} 
\]
\[
\mathfrak{g}_a= \left\{ \left( 
\begin{array}{rr|rr}
X_1 & Y_1 & Z_1 & T_1 \\ 
Y_1 & -X_1 & -T_1 & Z_1 \\ \hline
-^t \! Z_1 & ^t \! T_1 & -X_1 & -Y_1 \\ 
-^t \! T_1 & -^t \! Z_1 & -Y_1 & X_1%
\end{array}
\right) \ \mbox{\rm avec} \ 
\begin{array}{cc}
^t \! X_1 = -X_1, & ^t \! Y_1 = -Y_1 \\ 
^t \! Z_1 = -Z_1, & ^t \! T_1 = T_1 \\ 
& 
\end{array}
\right\} \label{g_a} 
\]
\[
\mathfrak{g}_b= \left\{ \left( 
\begin{array}{rr|rr}
X_2 & Y_2 & Z_2 & T_2 \\ 
-Y_2 & X_2 & T_2 & Z_2 \\ \hline
-^t \! Z_2 & -^t \! T_2 & -X_2 & -Y_2 \\ 
-^t \! T_2 & -^t \! Z_2 & Y_2 & -X_2%
\end{array}
\right) \ \mbox{\rm avec} \ 
\begin{array}{cc}
^t \! X_2 = -X_2, & ^t \! Y_2 = Y_2 \\ 
^t \! Z_2 = -Z_2, & ^t \! T_2 = -T_2 \\ 
& 
\end{array}
\right\} 
\]
\[
\mathfrak{g}_c= \left\{ \left( 
\begin{array}{rr|rr}
X_3 & Y_3 & Z_3 & T_3 \\ 
Y_3 & -X_3 & -T_3 & Z_3 \\ \hline
-^t \! Z_3 & ^t \! T_3 & X_3 & Y_3 \\ 
-^t \! T_3 & -^t \! Z_3 & Y_3 & -X_3%
\end{array}
\right) \ \mbox{\rm avec} \ 
\begin{array}{cc}
^t \! X_3= -X_3, & ^t \! Y_3 = -Y_3 \\ 
^t \! Z_3 = Z_3, & ^t \! T_3 = -T_3 \\ 
& 
\end{array}
\right\} 
\]
Notons que $dim \mathfrak{g}_e=m(2m+1), \, dim \mathfrak{g}_a=dim \mathfrak{g%
}_b=dim \mathfrak{g}_c=m(2m-1)$.

\begin{proposition}
Dans cette graduation $\mathfrak{g}_e$ est isomorphe \`a $sp(m)$ et toute $(%
\mathbb{Z}_2)^2$-graduation de $so(2m)$ telle que $\mathfrak{g}_e$ soit
isomorphe \`a $sp(m)$ est \'equivalente \`a la graduation ci-dessus.
\end{proposition}

En effet $\mathfrak{g}_e$ est simple de rang $m$ et de dimension $m(2m+1)$.
La deuxi\`eme partie r\'esulte de la classification donn\'ee dans \cite{B.G}
et \cite{B.G.R}.

\begin{coro}
Il n'existe, \`a \'equivalence pr\`es, qu'une seule structure d'espace
homog\`ene $(\mathbb{Z}_2)^2$-sym\'etrique sur l'espace homog\`ene compact $%
SO(2m)/Sp(m)$.
\end{coro}

Cette structure est associ\'{e}e \`{a} l'existence en tout point $x$ de $%
SO(2m)/Sp(m)$ d'un sous-groupe de $\mathcal{D}iff(M)$ isomorphe \`{a} $(%
\mathbb{Z}_{2})^{2}$. Notons $\Gamma _{x}$ ce sous-groupe. Il est enti\`{e}%
rement d\'{e}fini d\`{e}s que l'on connait $\Gamma _{\bar{1}}$ o\`{u} $\bar{1%
}$ est la classe dans $SO(2m)/Sp(m)$ de l'\'{e}l\'{e}ment neutre $1$ de $%
SO(2m)$. Notons 
\[
\Gamma _{\bar{1}}=\left\{ s_{e,\bar{1}},s_{a,\bar{1}},s_{b,\bar{1}},s_{c,%
\bar{1}}\right\} , 
\]%
les sym\'{e}tries $s_{\gamma ,\bar{1}}(x)=\pi (\rho _{\gamma }(A))$ o\`{u} $%
\pi :SO(2m)\rightarrow SO(2m)/Sp(m)$ est la submersion canonique, $x=\pi (A)$
et $\rho _{\gamma }$ est un automorphisme de $SO(2m)$ dont l'application
tangente en $1$ coincide avec $\tau _{\gamma }$. Ainsi 
\[
\left\{ 
\begin{array}{l}
\rho _{a}(A)=J_{a}^{-1}AJ_{a} \\ 
\rho _{b}(A)=J_{b}^{-1}AJ_{b} \\ 
\rho _{c}(A)=J_{c}^{-1}AJ_{c} \\ 
\end{array}%
\right. . 
\]%
Si $B\in \pi (A)$ alors il existe $P\in Sp(m)$ tel que $B=AP$. On a $%
J_{a}^{-1}BJ_{a}=J_{a}^{-1}AJ_{a}J_{a}^{-1}PJ_{a}=J_{a}^{-1}AJ_{a}$ car $P$
est invariante pour tous les automorphismes $\rho _{a},\rho _{b},\rho _{c}.$

\noindent Une m\'etrique non d\'eg\'en\'er\'ee $g$ invariante par $SO(2m)$
sur $SO(2m)/Sp(m)$ est adapt\'ee \`a la $(\mathbb{Z}_2)^2$-structure si les
sym\'etries $s_{x,\gamma}$ sont des isom\'etries c'est-\`a-dire si les
automorphismes $\rho _\gamma$ induisent des isom\'etries lin\'eaires.

\noindent Ceci implique que $g$ soit d\'efinie par une forme bilin\'eaire
sym\'etrique non d\'eg\'en\'er\'ee $B$ sur $\mathfrak{g}_a \oplus \mathfrak{g%
}_b \oplus \mathfrak{g}_c $ telle que les espaces $\mathfrak{g}_a,\, 
\mathfrak{g}_b,\, \mathfrak{g}_c$ soient deux \`a deux orthogonaux.
D\'eterminons toutes les formes bilin\'eaires $B$ v\'erifiant les
hypoth\`eses ci-dessus. Une telle forme s'\'ecrit donc 
\[
B=B _a +B _b+B _c 
\]
o\`u $B _a$(resp. $B _b$, resp. $B _c$ ) est une forme bilin\'eaire
sym\'etrique non d\'eg\'en\'er\'ee invariante par $\mathfrak{g}_e$ dont le
noyau contient $\mathfrak{g}_b \oplus \mathfrak{g}_c$ (resp. $\mathfrak{g}_a
\oplus \mathfrak{g}_c$, resp. $\mathfrak{g}_a \oplus \mathfrak{g}_b$).

\subsubsection{Exemples}

1) Dans le cas de la sph\`ere $SO(4)/Sp(2)$ la m\'etrique adapt\'ee \`a la
structure $(\mathbb{Z}_2)^2$-sym\'etrique est d\'efinie par la forme
bilin\'eaire $B$ sur $\mathfrak{g}_a \oplus \mathfrak{g}_b \oplus \mathfrak{g%
}_c$ qui est $ad(sp(2))$-invariante. Nous avons vu qu'une telle forme
s'\'ecrivait 
\[
B=\lambda _1^2\omega _1^2+\lambda _2^2\omega _2^2+\lambda _3^2\omega _3^2. 
\]
Elle est d\'efinie positive si et seulement si les coefficients $\lambda_i$
sont positifs ou nuls.

\medskip

2) Consid\'erons l'espace $(\mathbb{Z}_2)^2$-sym\'etrique compact $%
SO(8)/Sp(4)$. Afin de fixer les notations \'ecrivons la $(\mathbb{Z}_2)^2$%
-graduation de $so(8)$ ainsi: 
\[
\mathfrak{g}_a= \left\{ \left( 
\begin{array}{rr|rr}
X_1 & Y_1 & Z_1 & T_1 \\ 
Y_1 & -X_1 & -T_1 & Z_1 \\ \hline
-^t \! Z_1 & ^t \! T_1 & -X_1 & -Y_1 \\ 
-^t \! T_1 & -^t \! Z_1 & -Y_1 & X_1%
\end{array}
\right) \ \mbox{\rm avec} \ 
\begin{array}{cc}
^t \! X_1 = -X_1, & ^t \! Y_1 = -Y_1 \\ 
^t \! Z_1 = -Z_1, & ^t \! T_1 = T_1 \\ 
& 
\end{array}
\right\} 
\]
\[
\mathfrak{g}_b= \left\{ \left( 
\begin{array}{rr|rr}
X_2 & Y_2 & Z_2 & T_2 \\ 
-Y_2 & X_2 & T_2 & Z_2 \\ \hline
-^t \! Z_2 & -^t \! T_2 & -X_2 & -Y_2 \\ 
-^t \! T_2 & ^t \! Z_2 & Y_2 & -X_2%
\end{array}
\right) \ \mbox{\rm avec} \ 
\begin{array}{cc}
^t \! X_2 = -X_2, & ^t \! Y_2 = Y_2 \\ 
^t \! Z_2 = -Z_2, & ^t \! T_2 = -T_2 \\ 
& 
\end{array}
\right\} 
\]
\[
\mathfrak{g}_c= \left\{ \left( 
\begin{array}{rr|rr}
X_3 & Y_3 & Z_3 & T_3 \\ 
Y_3 & -X_3 & -T_3 & Z_3 \\ \hline
-^t \! Z_3 & ^t \! T_3 & X_3 & Y_3 \\ 
-^t \! T_3 & -^t \! Z_3 & Y_3 & -X_3%
\end{array}
\right) \ \mbox{\rm avec} \ 
\begin{array}{cc}
^t \! X_3= -X_3, & ^t \! Y_3 = -Y_3 \\ 
^t \! Z_3 = Z_3, & ^t \! T_3 = -T_3 \\ 
& 
\end{array}
\right\} 
\]
et pour la matrice $X_i$ (resp. $Y_i,Z_i,T_i$) on notera $X_i=\left( 
\begin{array}{cc}
0 & x_i \\ 
-x_i & 0%
\end{array}
\right)$ si elle est antisym\'etrique ou $X_i=\left( 
\begin{array}{cc}
x_i^1 & x_i^2 \\ 
x_i^2 & x_i^3%
\end{array}
\right)$ si elle est sym\'etrique, c'est \`a dire $X_i=\sum_j x_i^jX_i^j .$

Enfin on notera par les lettres $\alpha_i,\beta_i, \gamma_i, \delta_i$ les
formes lin\'eaires duales des vecteurs d\'efinis respectivement par les
matrices $X_i,Y_i,Z_i,T_i$. Ainsi si $X_i$ est antisym\'etrique, la forme
duale correspondante sera not\'ee $\alpha_i$, et si $X_i$ est sym\'etrique,
les formes duales $\alpha_i^1,\alpha_i^2,\alpha_i^3$ correspondent aux
vecteurs $\left( 
\begin{array}{cc}
1 & 0 \\ 
0 & 0%
\end{array}
\right),\left( 
\begin{array}{cc}
0 & 1 \\ 
1 & 0%
\end{array}
\right),\left( 
\begin{array}{cc}
0 & 0 \\ 
0 & 1%
\end{array}
\right) .$ Ceci \'etant la forme $B$ s'\'ecrit $B_a+B_b+B_c$ o\`u la forme $%
B_{\gamma}$ a pour noyau $\mathfrak{g}_{\gamma_1} \oplus \mathfrak{g}%
_{\gamma_2}$ avec $\mathfrak{g}_{\gamma} \neq \mathfrak{g}_{\gamma_1}$ et $%
\mathfrak{g}_{\gamma} \neq \mathfrak{g}_{\gamma_2}$. D\'eterminons $B_a$.
Comme elle est invariante par $ad(sp(2))$ on obtient: 
\[
\left\{ 
\begin{array}{l}
\medskip B_a(X_1,Y_1)=B_a(X_1,Z_1)=B_a(X_1,T_1^i)=0 \\ 
\medskip B_a(Y_1,Z_1)=B_a(Y_1,T_1^i)=0 \\ 
\medskip B_a(Z_1,T_1^i)=B_a(T_1^i,T_1^2)=0 \mbox{\rm \ pour } \, i=1,3 \\ 
\medskip B_a(X_1,X_1)=B_a(Y_1,Y_1)=B_a(Z_1,Z_1)=B_a(T_1^2,T_1^2)=0 \\ 
\medskip B_a(T_1^1,T_1^1)=B_a(T_1^3,T_1^3) \\ 
\medskip B_a(X_1,X_1)=2B_a(T_1^1,T_1^1)-2B_a(T_1^1,T_1^3)%
\end{array}
\right. 
\]
Ainsi la forme quadratique associ\'ee s'\'ecrit 
\[
q_{\mathfrak{g}_a}=\lambda _1(\alpha _1^2+\beta _1^2+\gamma _1^2+(\delta
_1^2)^2)+ \lambda _2((\delta _1^1)^2)+(\delta _1^3)^2)+(\lambda_2-\frac{
\lambda_1}{2} )((\delta _1^1)(\delta _1^3)) 
\]
soit 
\[
q_{\mathfrak{g}_a}=\lambda _1(\alpha _1^2+\beta _1^2+\gamma _1^2+(\delta
_1^2)^2)+ (\frac{3\lambda _2}{4}-\frac{\lambda _1}{8})(\delta _1^1+\delta
_1^3)^2+ (\frac{\lambda _2}{4}+\frac{\lambda _1}{8})(\delta _1^1-\delta
_1^3)^2. 
\]
De m\^eme nous aurons 
\[
q_{\mathfrak{g}_b}=\lambda _3(\alpha _2^2+(\beta _2^2)^2+\gamma _2^2+\delta
_2^2)+ (\frac{3\lambda _4}{4}-\frac{\lambda _3}{8})(\beta _2^1+\beta_2^3)^2+
(\frac{\lambda _4}{4}+\frac{\lambda _3}{8})(\beta_2^1-\beta_2^3)^2 
\]
et 
\[
q_{\mathfrak{g}_c}=\lambda _5(\alpha _3^2+\beta _3^2+(\gamma _3^2)^2+\delta
_3^2)+ (\frac{3\lambda _6}{4}-\frac{\lambda _5}{8})(\gamma
_3^1+\gamma_3^3)^2+ (\frac{\lambda _6}{4}+\frac{\lambda _5}{8}%
)(\gamma_3^1-\gamma_3^3)^2. 
\]

\noindent \textbf{Remarques.} 1. La forme $B$ d\'efinit une m\'etrique
riemannienne si et seulement si 
\[
\lambda_{2p}> \frac{\lambda _{2p-1}}{6}>0 
\]
pour $p=1,2,3$. Si cette contrainte est relach\'ee, la forme $B$, suppos\'ee
non d\'eg\'en\'er\'ee, peut d\'efinir une m\'etrique pseudo-riemannienne sur
l'espace $(\mathbb{Z}_2)^2$-sym\'etrique. Nous verrons cela dans le dernier
paragraphe.

2. Consid\'erons le sous-espace $\mathfrak{g}_e \oplus \mathfrak{g}_a$.
Comme $[\mathfrak{g}_a,\mathfrak{g}_a] \subset \mathfrak{g}_e$, c'est une
sous-alg\`ebre de $so(8)$ (ou plus g\'en\'eralement de $\mathfrak{g}$)
admettant une stucture sym\'etrique. La forme $B_a$ induit donc une
structure riemannienne ou pseudo-riemannienne sur l'espace sym\'etrique
associ\'e \`a l'espace sym\'etrique local $(\mathfrak{g}_e , \mathfrak{g}_a)$%
. Dans l'exemple pr\'ec\'edent $\mathfrak{g}_e \oplus \mathfrak{g}_a$ est la
sous-alg\`ebre de $so(8)$ donn\'ee par les matrices: 
\[
\left( 
\begin{array}{rr|rr}
X_1 & X_3 & X_4 & X_5 \\ 
-^t \! X_3 & X_2 & X_6 & -^t \! X_4 \\ \hline
-^t \! X_4 & -^t \! X_6 & X_2 & ^t \! X_3 \\ 
-^t \! X_5 & X_4 & -^t \! X_3 & X_2%
\end{array}
\right) \ \mbox{\rm avec} \ 
\begin{array}{cc}
^t \! X_1 = -X_1, & ^t \! X_2 = -X_2 \\ 
^t \! X_5 = X_5, & ^t \! X_6 = X_6 \\ 
& 
\end{array}
\]
Dans \cite{B}, on d\'etermine les espaces r\'eels en \'etudiant ces
structures sym\'etriques $\mathfrak{g}_e \oplus \mathfrak{g}_a$ donn\'ees
par deux automorphismes commutant de $\mathfrak{g}$. En effet si $\mathfrak{g%
}$ est simple r\'eelle et si $\sigma $ est un automorphisme involutif de $%
\mathfrak{g} $, il existe une sous-alg\`ebre compacte maximale $\mathfrak{g}%
_1$ de $\mathfrak{g}$ qui est invariante par $\sigma $ et l'\'etude des
espaces locaux sym\'etriques $(\mathfrak{g} ,\mathfrak{g}_e)$ se ram\`ene
\`a l'\'etude des espaces locaux sym\'etriques $(\mathfrak{g}_1 ,\mathfrak{g}%
_{11})$ o\`u $\mathfrak{g}_1$ est compacte. Dans ce cas $\mathfrak{g} $ est
d\'efinie \`a partir de $\mathfrak{g}_1$ par un automorphisme involutif $%
\tau $ commutant avec l'automorphisme $\sigma$. Ici notre approche est en
partie similaire mais le but est de regarder la structure des espaces non
sym\'etrique associ\'es aux paires $( \mathfrak{g},\mathfrak{g}_e ).$

\smallskip

Dans le cas particulier de l'espace $(\mathbb{Z}_2)^2$-sym\'etrique compact $%
SO(8)/Sp(4)$, l'alg\`ebre de Lie $\mathfrak{g}_e \oplus \mathfrak{g}_a $ est
isomorphe \`a $so(4) \oplus \mathbb{R}$ o\`u $\mathbb{R}$ d\'esigne
l'alg\`ebre ab\'elienne de dimension $1$. Notons \'egalement que chacun des
espaces sym\'etriques $\mathfrak{g}_e \oplus \mathfrak{g}_a, \mathfrak{g}_e
\oplus \mathfrak{g}_b, \mathfrak{g}_e \oplus \mathfrak{g}_c$ est isomorphe
\`a $so(4) \oplus \mathbb{R}$. Mais ceci n'est pas g\'en\'eral, les
alg\`ebres sym\'etriques peuvent ne pas \^etre isomorphes ni m\^eme de
m\^eme dimension. L'espace sym\'etrique compact connexe associ\'e est
l'espace homog\`ene ${Su(4)}/{Sp(2)} \times \mathbb{T}$ o\`u $\mathbb{T}$
est le tore \`a une dimension. C'est un espace riemannien sym\'etrique
compact non irr\'eductible. La m\'etrique $q_{\mathfrak{g}_a}$ d\'efinie
pr\'ec\'edemment correspond \`a une m\'etrique riemannienne ou pseudo-riemannienne sur cet espace. La restriction au premier facteur correspond
\`a la m\'etrique associ\'ee \`a la forme de Killing Cartan sur $su(4)$.
Elle correspond \`a $\lambda_2=\frac{ \lambda_1}{2}$.

\subsubsection{Cas g\'{e}n\'{e}ral: m\'{e}triques adapt\'{e}es sur $%
SO(2m)/Sp(m)$}

\noindent {\bf Notations.} Nous avons \'ecrit une matrice g\'en\'erale de $\mathfrak{g}_a$
sous la forme (\ref{g_a}). Si on note $(X_1,Y_1,Z_1,T_1)$ un \'el\'ement de $%
\mathfrak{g}_a$, on consid\`ere la base de $\mathfrak{g}_a$, $%
\{X_{1,ij},Y_{1,ij},Z_{1,ij},T_{1,ij}\}$ correspondant aux matrices
\'el\'ementaires . La base duale sera not\'ee $(\alpha _{a,ij},\beta_{a,ij},
\gamma_{a,ij}, \delta_{a,ij} )$. Rappelons que $X_1,Y_1,Z_1$ sont
antisym\'etriques alors que $T_1$ est sym\'etrique. Les crochets
correspondent aux repr\'esentations de $so(\frac{m}{2})$ sur lui-m\^eme ou
de $so(\frac{m}{2})$ sur l'espace des matrices sym\'etriques. On aura donc 
\[
q_{\mathfrak{g}_a}=\lambda _1^a(\sum (\alpha _{a,ij}^2+\beta
_{a,ij}^2+\gamma _{a,ij}^2)+\sum_{i \neq j}\delta _{a,ij}^2)+ \lambda
_2^a(\delta _{a,ii}^2)+(\lambda_2^a-\frac{ \lambda_1^a}{2}
)(\sum_{i<j}(\delta _{a,ii}\delta _{a,jj}). 
\]
Les formes $q_{\mathfrak{g}_b}$ et $q_{\mathfrak{g}_c}$ admettent une
d\'ecomposition analogue, en tenant compte du fait que dans $\mathfrak{g}_b$
ce sont les matrices $Y_2$ qui sont sym\'etriques et pour $\mathfrak{g}_a$
les matrices $Z_1$ (\ref{g_a}). On note $(\alpha _{b,ij},\beta_{b,ij},
\gamma_{b,ij}, \delta_{b,ij} )$ la base duale de $%
\{X_{2,ij},Y_{2,ij},Z_{2,ij},T_{2,ij}\}$ et par $(\alpha
_{c,ij},\beta_{c,ij}, \gamma_{c,ij}, \delta_{c,ij} )$ la base duale de $%
\{X_{3,ij},Y_{3,ij},Z_{3,ij}, T_{3,ij}\}$.

\begin{proposition}
Toute m\'{e}trique non d\'{e}g\'{e}n\'{e}r\'{e}e adapt\'{e}e \`{a} la
structure $(\mathbb{Z}_{2})^{2}$-sym\'{e}trique de l'espace homog\`{e}ne $%
SO(2m)/Sp(m)$ est d\'{e}finie \`{a} partir de la forme bilin\'{e}aire $ad(%
\mathfrak{g}_{e})$-invariante sur $\mathfrak{g}_{a}\oplus \mathfrak{g}%
_{b}\oplus \mathfrak{g}_{c}$ $B=q_{\mathfrak{g}_{a}}+q_{\mathfrak{g}_{b}}+q_{%
\mathfrak{g}_{b}}$ avec 
\[
\left\{ 
\begin{array}{l}
\medskip q_{\mathfrak{g}_{a}}=\lambda _{1}^{a}\left( \sum (\alpha
_{a,ij}^{2}+\beta _{a,ij}^{2}+\gamma _{a,ij}^{2}\right) +\sum_{i\neq
j}\delta _{a,ij}^{2})+\lambda _{2}^{a}(\delta _{a,ii}^{2})+(\lambda _{2}^{a}-%
\frac{\lambda _{1}^{a}}{2})(\sum_{i<j}(\delta _{a,ii}\delta _{a,jj}) \\ 
\medskip q_{\mathfrak{g}_{b}}=\lambda _{1}^{b}(\sum (\alpha
_{b,ij}^{2}+\gamma _{ij}^{2})+\delta _{b,ij}^{2}+\sum_{i\neq j}\beta
_{b,ij}^{2})+\lambda _{2}^{b}(\beta _{b,ii}^{2})+(\lambda _{2}^{b}-\frac{%
\lambda _{1}^{b}}{2})(\sum_{i<j}(\beta _{b,ii}\beta _{b,jj}) \\ 
\medskip q_{\mathfrak{g}_{c}}=\lambda _{1}^{c}(\sum (\beta
_{c,ij}^{2}+\gamma _{c,ij}^{2})+\delta _{c,ij}^{2}+\sum_{i\neq j}\alpha
_{c,ij}^{2})+\lambda _{2}^{c}(\alpha _{c,ii}^{2})+(\lambda _{2}^{c}-\frac{%
\lambda _{1}^{c}}{2})(\sum_{i<j}(\alpha _{c,ii}\alpha _{c,jj})%
\end{array}%
\right. 
\]
\end{proposition}

Soit $\gamma \in \{a,b,c\}$. Les valeurs propres de la forme $q_{\mathfrak{g}%
_\gamma }$ sont 
\[
\mu_{1,\gamma }=\lambda_1^\gamma , \ \ \mu _{2,\gamma }=\lambda _2^\gamma
/2+\lambda _{1}^\gamma /4, \ \ \mu_{3,\gamma }=\lambda _2^\gamma\frac{r+1}{2}%
-\lambda _{1}^\gamma \frac{r-1}{4} 
\]
o\`u $r$ est l'ordre commun des matrices sym\'etriques $X_4,Y_2,Z_1$. Ces
valeurs propres sont respectivement de multiplicit\'e $dim \mathfrak{g}%
_{\gamma }-r,r-1,1$. Le signe des valeurs propres $\mu _{2,\gamma }$ et $\mu
_{3,\gamma }$ est donc 
\[
\left\{ 
\begin{array}{l}
\mu _{2,\gamma } > 0 \Longleftrightarrow \lambda_2^\gamma >-\lambda_1^\gamma
/2 \\ 
\mu _{3,\gamma } > 0 \Longleftrightarrow \lambda_2^\gamma >
-\lambda_1^\gamma \frac{r-1}{2(r+1)}%
\end{array}
\right. 
\]
On en d\'eduit, si $s(q)$ d\'esigne la signature de la forme quadratique $q$
: 
\[
\left\{ 
\begin{array}{lll}
\medskip s(q_{\mathfrak{g}_\gamma }) & = (dim \mathfrak{g}_\gamma ,0) & 
\Leftrightarrow ( \lambda_1^\gamma >0, \ \lambda_2^\gamma > \lambda_1^\gamma 
\frac{r-1}{2(r+1)}) \\ 
\medskip & = (dim \mathfrak{g}_\gamma -1 ,1) & \Leftrightarrow (
\lambda_1^\gamma >0, \ -\lambda_1^\gamma /2< \lambda_2^\gamma
<\lambda_1^\gamma \frac{r-1}{2(r+1)}) \\ 
\medskip & = (dim \mathfrak{g}_\gamma -r ,r) & \Leftrightarrow (
\lambda_1^\gamma >0, \ \lambda_2^\gamma <-\lambda_1^\gamma /2) \\ 
\medskip & = (r,dim \mathfrak{g}_\gamma -r) & \Leftrightarrow (
\lambda_1^\gamma <0, \ \lambda_2^\gamma >-\lambda_1^\gamma /2) \\ 
\medskip & = (1,dim \mathfrak{g}_\gamma -1) & \Leftrightarrow (
\lambda_1^\gamma <0, \ \lambda_1^\gamma \frac{r-1}{2(r+1)} <
\lambda_2^\gamma <-\lambda_1^\gamma /2) \\ 
\medskip & = (0, dim \mathfrak{g}_\gamma ) & \Leftrightarrow (
\lambda_1^\gamma <0, \ \lambda_2^\gamma < \lambda_1^\gamma \frac{r-1}{2(r+1)}
\\ 
&  & 
\end{array}
\right. 
\]
Notons que $\mu_{2,\gamma }= \mu_{3,\gamma }$ si et seulement si $%
\lambda_1^\gamma= 2\lambda_2^\gamma.$

\subsubsection{Classification des m\'{e}triques riemanniennes adapt\'{e}es
sur $SO(2m)/Sp(m)$}

Comme $r=\frac{m^2+m-2}{m^2+m+2}$ on a le r\'esultat suivant

\begin{theo}
Toute m\'etrique riemannienne sur $SO(2m)/Sp(m)$ adapt\'ee \`a la structure $%
(\mathbb{Z}_2)^2$-sym\'etrique est d\'efinie \`a partir de la forme
bilin\'eaire sur $\mathfrak{g}_a\oplus \mathfrak{g}_b\oplus \mathfrak{g}_c$
\[
B= q_{\mathfrak{g}_a }( \lambda_1^a,\lambda_2^a)+q_{\mathfrak{g}_b }(
\lambda_1^b,\lambda_2^b)+q_{\mathfrak{g}_b }( \lambda_1^b,\lambda_2^b) 
\]
avec 
\[
\left\{ 
\begin{array}{l}
\lambda_1^\gamma >0 \\ 
\lambda_2^\gamma > \lambda_1^\gamma \frac{m^2+m-2}{2(m^2+m+2}) \\ 
\end{array}
\right. 
\]
pour tout $\gamma \in \{a,b,c\}$.
\end{theo}

Pour une telle m\'etrique, la connexion de Levi-Civita ne co\"{\i}ncide pas
en g\'en\'eral avec la connexion canonique associ\'ee \`a la structure $(%
\mathbb{Z}_2)^2$-sym\'etrique (\cite{B.G.R}). Ces deux connexions sont les
m\^emes si et seulement si la m\'trique riemannienne est naturellement
r\'eductive. Elle correspond donc \`a la restriction de la forme de Killing
(au signe pr\`es) de $SO(2m)$. Cette m\'etrique correspond \`a la forme
bilin\'eaire $B$ d\'efinie par les param\`etres 
\[
\lambda_1^a= \lambda_1^b=\lambda_1^c=
2\lambda_2^a=2\lambda_2^b=2\lambda_2^c. 
\]

\subsubsection{Classification des m\'{e}triques lorentziennes adapt\'{e}es
sur l'espace $SO(2m)/Sp(m)$}

Les m\'etriques lorentziennes adapt\'ees \`a la structure $(\mathbb{Z}_2)^2$%
-sym\'etrique sont d\'efinies par les formes bilin\'eaires $B$, d\'efinies
dans la section pr\'ec\'edente, dont la signature est $(dim(\mathfrak{g}%
_a)+dim(\mathfrak{g}_b)+dim(\mathfrak{g}_c)-1,1)$. On a donc

\begin{theo}
Toute m\'etrique lorentzienne sur $SO(2m)/Sp(m)$ adapt\'ee \`a la structure $%
(\mathbb{Z}_2)^2$-sym\'etrique est d\'efinie par l'une des formes
bilin\'eaires 
\[
B= q_{\mathfrak{g}_a }( \lambda_1^a,\lambda_2^a)+q_{\mathfrak{g}_b }(
\lambda_1^b,\lambda_2^b)+q_{\mathfrak{g}_b }( \lambda_1^b,\lambda_2^b) 
\]
avec 
\[
\left\{ 
\begin{array}{l}
\medskip \forall \gamma \in \{a,b,c\}, \ \lambda_1^\gamma >0 \\ 
\medskip \exists \gamma _0 \in \{a,b,c\} \ \mbox{\rm tel \ que} \
-\lambda_1^{\gamma_0} /2< \lambda_2^{\gamma_0} <\lambda_1^{\gamma_0} \frac{%
r-1}{2(r+1)} \\ 
\medskip \forall \gamma \neq \gamma _0, \ \ \lambda_2^\gamma >
\lambda_1^\gamma \frac{r-1}{2(r+1)}. \  \\ 
\end{array}
\right. 
\]
\end{theo}

\end{document}